\documentclass[11pt,letter]{article}

\usepackage{booktabs, multirow}
\usepackage{fancyhdr}
\usepackage{tikz}
\usetikzlibrary{arrows.meta}
\usetikzlibrary{positioning}
\usepackage{bbm}
\usepackage{natbib}
\usepackage{amstext}
\usepackage{amsthm}
\usepackage{epsfig}
\usepackage{pdfpages}
\usepackage{amsmath,amssymb,amsthm}
\usepackage{graphicx}
\usepackage{calc}
\usepackage{dsfont}
\usepackage{rotating}
\usepackage{multirow}
\usepackage{amsmath}
\usepackage{lscape}
\usepackage{tabularx}
\usepackage{bm}
\usepackage{color}
\usepackage{enumerate,array}
\usepackage{verbatim}
\usepackage{stackengine}
\usepackage{enumitem}
\usepackage{mathrsfs}

\newcommand{\indep}{\!\perp\!\!\!\perp\!}

\usepackage{hyperref}
\definecolor{cornellred}{rgb}{0.7, 0.11, 0.11}
\hypersetup{    
	colorlinks=true,
	linkcolor=blue,
	citecolor=blue,
	urlcolor=black
}

\renewcommand{\baselinestretch}{1.5}

\numberwithin{equation}{section}

\renewcommand{\Re}{{\rm Re}}

\allowdisplaybreaks

\newtheorem{theorem}{Theorem}[section]
\newtheorem{lemma}[theorem]{Lemma}
\newtheorem{proposition}[theorem]{Proposition}

\newtheorem{corollary}[theorem]{Corollary}

\theoremstyle{definition}
\newtheorem{definition}[theorem]{Definition}
\newtheorem{remark}[theorem]{Remark}
\newtheorem{algorithm}[theorem]{Algorithm}
\newtheorem{example}[theorem]{Example}

\newtheorem{assumption}{Assumption}

\newcommand{\E}{{\mathbb E}}
\newcommand{\ft}{{\rm ft\!}}

\renewcommand{\H}{\mathcal{H}}

\newcommand{\as}{{\rm a.s.}}

\newcommand{\one}{\mathbbm{1}}
\newcommand{\dvar}{{\rm dVar}}
\newcommand{\dcov}{{\rm dCov}}
\newcommand{\dcor}{{\rm dCor}}
\newcommand{\hsic}{{\rm HSIC}}
\newcommand{\ed}{{\rm Ed}}
\newcommand{\mmd}{{\rm MMD}}

\def\T{{\mathcal T}}
\def\d{{\rm d}}

\def\trans{^{\mkern-1mu\mathsf{T}\mkern-1mu}}
\def\tr{{\rm tr}}

\def\B{\mathcal B}

\def\K{\mathcal K}

\newcommand{\var}{{\rm Var}}

\newcommand{\e}{\epsilon}

\renewcommand{\S}{\mathcal{S}}

\newcommand{\F}{\mathcal{F}}
\newcommand{\V}{\mathcal{V}}

\newcommand{\N}{\mathcal{N}}

\newcommand{\converged}{\overset{d.}{\longrightarrow}}

\newcommand{\iidsim}{\overset{\small{\text{iid}}}{\sim}}
\newcommand{\diag}{{\rm diag}}
\newcommand{\supp}{{\rm supp}}

\renewcommand{\l}{\langle}
\renewcommand{\r}{\rangle}
\renewcommand{\i}{{\rm i}}
\renewcommand{\phi}{\varphi}

\renewcommand{\tilde}{\widetilde}
\renewcommand{\hat}{\widehat}
\renewcommand{\epsilon}{\varepsilon}
\renewcommand{\P}{{\rm P}}

\def\boxit#1{\vbox{\hrule\hbox{\vrule\kern6pt  \vbox{\kern6pt#1\kern6pt}\kern6pt\vrule}\hrule}}

\newtheorem{innercustomgeneric}{\customgenericname}
\providecommand{\customgenericname}{}
\newcommand{\newcustomtheorem}[2]{%
  \newenvironment{#1}[1]
  {%
   \renewcommand\customgenericname{#2}%
   \renewcommand\theinnercustomgeneric{##1}%
   \innercustomgeneric
  }
  {\endinnercustomgeneric}
}
\newcustomtheorem{customthm}{Assumption}
\newcustomtheorem{customcor}{Corollary}

\renewcommand{\baselinestretch}{1.4}

\makeatletter
\def\singlespace{\deltaf\baselinestretch{1}\@normalsize}

  \renewenvironment{thebibliography}[1]{%
    \begin{oldthebibliography}{#1}%
      \setlength{\parskip}{0.3ex}%
      \setlength{\itemsep}{0ex}%
  }%
  {%
    \end{oldthebibliography}%
  }

\title{New energy distances 
for statistical inference on infinite dimensional Hilbert spaces without moment conditions }
\author{
  { Holger Dette} \\
{ Ruhr-Universit\"at Bochum} \\
{ Fakult\"at f\"ur Mathematik} \\
{ 44780 Bochum, Germany} \\
 \and 
Jiajun Tang \\
{Harvard University} \\
{Department of Statistics} \\
{ 02138 Cambridge, MA, USA} \\
}
\date{}

\begin{document}

\maketitle

\baselineskip=16pt

\begin{abstract}
For statistical inference on an infinite-dimensional Hilbert space $\H $ with no moment conditions we 
introduce a new class of energy distances on the 
space of probability measures on $\H$.  The proposed  distances  consist of the integrated squared modulus of the corresponding difference of the characteristic functionals with respect to a reference probability measure on the Hilbert space. Necessary and sufficient conditions are established for the reference probability measure to be {\em characteristic}, the property that guarantees that  the distance defines a metric 
on the space of   probability measures on $\H$. 
We also use these results to define new distance covariances, which can be used to measure  the dependence between the marginals of a two dimensional distribution of $\H^2$ without existing moments.

On the basis of the new  distances we  develop  statistical inference for Hilbert space valued data, which does not
require any moment assumptions. As a consequence,  our methods are robust with respect to heavy tails in finite dimensional data. In particular, we consider the problem of 
  comparing the distributions of two samples and    the problem of testing for independence and  construct  new minimax optimal tests for the corresponding hypotheses. We also develop 
aggregated (with respect to the reference measure) procedures for power enhancement  and investigate the  finite-sample properties  by means of  a simulation study. 
\end{abstract}

\vspace{0.5em}



\parindent 0cm 

\section{Introduction}\label{sec:intro}

Comparing the distributions of two  samples 
$X_1,\ldots,X_{n_1}\iidsim \mathbb P^X$ and $Y_1,\ldots,Y_{n_2}\iidsim \mathbb P^Y$
or investigating the  dependence 
between the components of  the sample 
$(X_1,Y_1),\ldots,(X_n,Y_n)\iidsim \mathbb P^{X,Y}$
are fundamental
problems in statistics and an enormous amount of literature has been published on these subjects. Many authors address these problems by testing  the equality of the marginal distribution, that is $\mathbb P^X=\mathbb P^Y$ (\textit{testing homogeneity}), or  by testing for a product structure, that is  $\mathbb P^{X,Y}=\mathbb P^X\otimes\mathbb P^Y$ (\textit{testing independence}). These (and other testing problems) usually require  the definition   of a  distance between probability distributions  which vanishes, if and only if  the  null hypothesis is satisfied, and which is estimated from the data. While  
classical work  considers distances between distributions on Euclidean spaces, in the last decades  advanced  sampling technology produces large data volumes with complex structures, which cannot be modeled 
by  random objects on  finite dimensional vector  spaces. As a consequence, new tools  to model and analyze data in spaces of such  complex structure have been developed. These include functional data analysis \citep{bosq2000,horvath2012,hsing2015,wang2016}, topological data analysis \citep{wasserman}, statistical shape analysis \citep{dryden2016} and manifold learning \citep{bhapat2003,bhapat2005,maandfu}, where various distances between random objects (more precisely their distributions) have been considered.

Recently, there has been also some interest  on  statistical analysis  of heavy-tailed  data in these fields. For example \cite{Gervini,
Boente2021} and \cite{shao2023} investigate heavy tailed  functional data   and propose  robust methods to cope with the ubiquitous cases where moment conditions on the underlying distribution of the random functions are violated. 
Robust inference tools are also proposed in topological data analysis \citep{Buchet2016,Bobrowski,Fasy}, in manifold learning \citep{chen2006,feng}, in statistical shape analysis \citep{schmid2011robust} among others. The research presented in this paper is 
motivated by these works, and develops statistical tools  for homogeneity and independence testing 
for data  on infinite dimensional spaces, which do not require the existence of moments of the underlying population distribution. 
\medskip

\textbf{Main contributions.}  
Our approach is based on the definition of  
a class of energy distances, say  $\ed_\nu$,  on the space of  probability distributions $\mathcal{P} (\H) $
on an infinite-dimensional Hilbert space $\H $, which do not require any moment conditions. This class of metrics is indexed by Borel probability measures $\nu$  on $\H$, which are  called reference measures.  We  establish sufficient and necessary conditions on the reference measure 
such that the corresponding  distance metric is characteristic,  that is   $\ed_\nu^2(\mathbb P,\mathbb Q)=0$ if and only if $\mathbb P=\mathbb Q$ ($\mathbb P,\mathbb Q\in\mathcal P(\H)$). We demonstrate that  these conditions 
are  satisfied   by the Gaussian and Laplace reference measure (under mild assumptions). 
From a more statistical point of  view  we will use these distances to construct tests for the hypotheses  of  homogeneity and independence and prove  their (minimax)  optimality  in terms of uniform separation rates. 
Our results are in particular applicable to Hilbert-space valued data in  the presence of additive measurement errors under mild conditions.

\medskip

\textbf{Related works.} 
The energy distance and maximum mean discrepancy  (e.g.~\citealp{smola2007,gretton2012,szekely2013}) are  common metrics 
for comparing distributions. 
For example,
\cite{zhu2022}  use the energy distance to develop  a homogeneity test for  sparsely observed functional data,   and 
\cite{wynne2022} propose a maximum mean discrepancy in a
reproducing kernel Hilbert space (RKHS)  for this purpose. 
In this context it is of importance 
that the null hypothesis of equal distributions is characterized by a vanishing distance (that is: the corresponding RKHS embedding maps all distributions uniquely). This property has been established   by 
\cite{sriperumbudur2010} for  integrally strictly positive definite kernels. We also refer to \cite{hlavka2022} who propose Cramér-von Mises and  ANOVA tests for functional data.

Similarly, the distance covariance and Hilbert-Schmidt independence criterion  are frequently used to measure dependence (see, for example,~\citealp{gretton,fukumizu2007,szekely2007}). For random variables in infinite dimensional Hilbert spaces, independence tests have been proposed using basis projection  \citep{miao2022,krzysko2022} and 
using the sup-norm  approaches \citep{bhar2023}. \cite{meintanis2022} suggest  to measure independence via integrating the squared modulus of difference in point-wise characteristic functions. However, it is unclear if their distance covariance is well-defined since their approach is based on marginal independence. A conditional mean independence test for functional data has been proposed, using the martingale difference divergence approach in \cite{lee2020} and the dimensional reduction approach in \cite{patilea2016}. \cite{dehling2020} investigate a distance covariance for discretized stochastic processes, while  \cite{lyons2013} considers general  metric spaces and  proves necessary and sufficient conditions on the metric space such that distance covariance is a well-defined measure of independence.


A common feature of these works
is that they usually require at least some condition on the underlying distribution of interest. For example, \cite{patilea2016} and \cite{lee2020}
assume 
the existence of the conditional mean  while others authors use  moment conditions, such as 
$\int d(x,y)\mathbb P(\d y)<\infty$
for some semimetric $d$ and  element $x$ in the metric space (see, for example,~\citealp{lyons2013,miao2022}).
  Since probability distributions on infinite-dimensional spaces and general metric spaces are usually difficult to capture, the validity of these conditions might not be easy to justify in applications.
These aspects motivate us to define energy distances and distance covariances  for distributions of infinite dimensional spaces without any moment conditions. Statistical methodology based on appropriate estimates of these distances is more suitable for the analysis of  heavy-tailed data, and we demonstrate the advantages of our approach in the construction of homogeneity and independence tests with very good finite sample properties for this type of data.

\medskip

\textbf{Organization.}
 In Sections~\ref{sec:ed} and \ref{sec:dcov}, we define a class of energy distances and distance covariances for probability measures on infinite dimensional Hilbert spaces with no  moments and investigate several properties of these metrics. In Section~\ref{sec:permu:indep} and \ref{sec:twosample} we use these results to develop  independence tests and two-sample tests, respectively, where we also focus on optimality  properties of the procedures. 
In Section~\ref{sec:me}, we consider the inference problem where the data are measured with error. An aggregated procedure is proposed in Section~\ref{sec:agg}, and the finite-sample properties of the new methods are investigated in 
Section~\ref{sec:simu} by means of a  simulation study.

\section{Energy distances in Hilbert space}\label{sec:ed}

Let $\H$ be a real separable infinite-dimensional Hilbert space with inner product $\l\cdot,\cdot\r$, norm $\|\cdot\|$, and Borel $\sigma$-algebra $\mathcal B(\H)$.
For an $\H$-valued random variable on  a probability space $(\Omega,\mathcal F,{\rm P})$, that is a $(\mathcal F,\mathcal B(\H))$ measurable
mapping $X:\Omega\to\H$, we denote by $\mathbb P^X={\rm P}\circ X^{-1}$ the distribution of $X$.
%
Let $\mathcal P(\H)$ and $\mathcal P(\H^2)$ denote the sets of all Borel probability measures on $\H$ and $\H^2$, respectively. 
For functionals $f_1,f_2:\H\to\mathbb C$, we define   their tensor product  $f_1\otimes f_2:\H^2\to\mathbb C$ by $f_1\otimes f_2(w_1,w_2)=f_1(w_1)f_2(w_2)$  ($w_1,w_2\in\H$). 
A compact operator $C: \H \to \H$ is a trace-class operator if  ${\rm tr}(C)=\sum_{j=1}^\infty|\l Ce_j,e_j\r|<\infty$, where $\{e_j\}_{j=1}^\infty$  is a complete orthogonal basis of $\H$. A self-adjoint, positive definite and trace-class operator on $\H $ is called $\mathcal S$-operator (see Section~4.5 of \citealp{hsing2015}). Throughout this paper we denote by  $\mathcal S$  the set of $\mathcal S$-operators with a slight abuse of notation. 
In the simulation study we will  be interested in the Hilbert space $\H=L^2([0,1])$ of square integrable functions with respect to  the Lebesgue measure on $[0,1]$. A mean-square continuous process $\{X(t)\}_{t\in[0,1]}$ such that $X(t,\omega)$ is jointly measurable is a random element in $L^2([0,1])$ (see Section~7.4 in \citealp{hsing2015}).



\subsection{An energy distance for measures with no moments}  \label{sec21}

For a Borel probability measure $\nu$ on $\H $ 
we denote by 
\begin{align}\label{phinu}
\phi_{\nu}(w)=\int_\H \exp(\i\l w,z\r)\,\nu(\d z)\,,\qquad w\in \H \,  
\end{align}
its characteristic functional (see, for example, \citealp{prokhorov1956}). If 
 $X$ and $Y$ are  $\H$-valued random variables with distributions $\mathbb P^X$ and $\mathbb P^Y$,   their characteristic functionals are defined by
\begin{align}\label{phix}
\phi_X(w)&:= \phi_{\mathbb P^X} (w)  = \E \exp(\i\l X,w\r)\,,\quad \phi_Y(w) := \phi_{\mathbb P^Y} (w)  = \E \exp(\i\l Y,w\r)\,,\quad w\in\H\,, 
\end{align}
respectively, and the joint characteristic functional of $(X,Y)$ is given by
\begin{align}\label{phixy}
\phi_{X,Y}(w_1,w_2) := \phi_{\mathbb P^{X,Y}} (w_1,w_2) 
=\E \exp(\i\l X,w_1\r+\i\l Y,w_2\r)\,,  \qquad w_1,w_2\in\H\,.
\end{align}

The following theorem states that probability measures on a separable Hilbert space and its characteristic functional are mutually uniquely identified. The proof of sufficiency is straightforward in view of the definition of characteristic functionals in \eqref{phinu}, and the necessity follows by directly applying Corollary~4.1 in \cite{kukush2020}.

\begin{theorem}
\label{thm:bijection}
Let $\H$ be a separable Hilbert space and $\mathbb P_1,\mathbb P_2\in\mathcal P(\H)$ with characteristic functionals $\phi_1$ and $\phi_2$, respectively. Then, $\mathbb P_1=\mathbb P_2$ if and only if $\phi_1=\phi_2$.
\end{theorem}



By Theorem~\ref{thm:bijection}, for $\mathbb P^X,\mathbb P^Y\in\mathcal P(\H)$, $\mathbb P^X=\mathbb P^Y$ is equivalent to $\phi_X=\phi_Y$. Therefore
we will use  the characteristic functionals to  define a class of distances between probability distributions on the Hilbert space $\H$ that do not require moment conditions. More precisely,  we define the energy distance between $\mathbb P^X$ and $\mathbb P^Y$ by integrating the squared modulus of the difference between their corresponding characteristic functionals, with respective to a Borel probability measure $\nu$ on $\H$.

\begin{definition}
\label{def:ed}

Let $\nu$ be a Borel probability measure  on $\H$. Define the (squared) energy distance between $\mathbb{P}^X$ and $\mathbb{P}^Y$ with respect to $\nu$ by
\begin{align}\label{ed}
\ed_\nu^2(\mathbb{P}^X,\mathbb{P}^Y)=\int_\H|\phi_X(w)-\phi_Y(w)|^2\,\nu(\d w)\,,\qquad \mathbb P^X,\mathbb P^Y\in\mathcal P(\H),
\end{align}
where $\phi_X$ and $\phi_Y$ are characteristic functionals defined  in \eqref{phix}.

\end{definition}

In its definition, we use the notation $\ed_\nu$ to emphasize the dependence of the energy distance on $\nu$.
Since $\nu$ is a probability measure on $\H$ and $|\phi_X|,|\phi_Y|\leq1$ it follows that 
$\ed_\nu^2(\mathbb{P}^X,\mathbb{P}^Y)\leq 4$ for  all $\mathbb P^X,\mathbb P^Y\in\mathcal P(\H)$, so the energy distance in \eqref{ed} does not require any moment conditions on $\mathbb P^X$ or $\mathbb P^Y$. 
Throughout this article we use the notation 
 $-B=\{w\in\mathcal B(\H):-w\in B\}$ for a set  $B\in\mathcal B(\H)$ and  make the following crucial assumptions on the probability measure $\nu$.

\begin{assumption}
\label{a:sym}
The probability measure $\nu$ is symmetric on $\H$, that is: for any $B\in\mathcal B(\H)$, it holds that $\nu(B)=\nu(-B)$.
\end{assumption}

\begin{assumption}
\label{a:supp}
The probability measure $\nu$ has full support on $\H$, i.e. $\supp(\nu)=\H$.
\end{assumption}

A direct consequence of the symmetry condition in Assumption~\ref{a:sym} is that the characteristic functional $\phi_\nu$ is a real-valued functional, that is
\begin{align}\label{lem:sym}
\phi_{\nu}(x)=\int_\H \cos\l x,y\r\,\nu(\d y)\,,\qquad \text{for }\,x\in \H \,.
\end{align}
In addition, $\phi_{\nu}(x)=\phi_{\nu}(-x)$ for any $x\in\H $.
We will show later  (see Theorem~\ref{thm:iff:ed} and Corollary~\ref{cor:4.11} below) that 
Assumption~\ref{a:supp} guarantees that the corresponding energy distance is a well-defined metric  on $\mathcal P(\H)$.  Throughout this article, we say that the probability measure $\nu$ is a {\em reference probability measure}, if $\nu$ satisfies Assumptions~\ref{a:sym} and \ref{a:supp}.

For an $\H$-valued random variable $Z$ and a reference probability measure $\nu$ with characteristic functional $\phi_{\nu}$, it follows from Fubini's theorem that
\begin{align}\label{phiz}
\int_\H\phi_Z(w)\,\nu(\d w)=\E\int_\H\exp(\i \l Z,w\r)\,\nu(\d w)=\E\phi_{\nu}(Z)\,.
\end{align}
This observation enables us to derive the following proposition which expresses the energy distance in terms of $\phi_{\nu}$; see Section~\ref{app:prop:ed} of the supplementary materialfor a proof.

\begin{proposition}\label{prop:ed}
Suppose  that $\nu$ is a symmetric probability measure on $\H$ 
and let $(X',Y')$ be an independent copy of $(X,Y)$. Then the energy distance $\ed_\nu$  in \eqref{ed} satisfies 
\begin{align}\label{ed:dnu}
\ed_\nu^2(\mathbb P^X,\mathbb P^Y)&=\E\phi_{\nu}(X-X')+\E \phi_{\nu}(Y-Y')-2\E\phi_{\nu}(X-Y')\,.
%
\end{align}
\end{proposition}

We are now in a position to state  a necessary and sufficient condition for  the  probability measure $\nu$, such that the property   
\begin{align}\label{statement2}
\text{$\ed_\nu^2(\mathbb P,\mathbb Q)=0$ if and only if $\mathbb P=\mathbb Q$\,.}
\end{align}
holds for all
probability distributions $\mathbb P,\mathbb Q\in\mathcal P(\H)$ on $\H$.
We call  a probability measure $\nu$  with this property   {\em characteristic}.
This term 
is borrowed from  the literature on characteristic kernels in reproducing kernel Hilbert spaces (RKHS), see, for example, \cite{fukumizu2007} and \cite{sriperumbudur2010}.

\begin{theorem}\label{thm:iff:ed}
Suppose that the probability measure $\nu$ is symmetric.
Then, $\nu$ is characteristic (i.e.~statement \eqref{statement2} is true) if and only if $\supp(\nu)=\H$.
\end{theorem}

A detailed proof of  Theorem~\ref{thm:iff:ed}
can be found in Section~\ref{app:thm:iff} of the supplementary material. Here we indicate the main steps 
to prove necessity. Suppose  that the symmetric reference probability measure $\nu$ is characteristic, but is not fully supported on $\H$. We aim to construct two probability measures  $\mathbb P_1,\mathbb P_2\in\mathcal P(\H)$ with characteristic functionals $\phi_1,\phi_2$ such that $\ed_\nu^2(\mathbb P_1,\mathbb P_2)=0$, but $\mathbb P_1\neq\mathbb P_2$. Since $\nu$ is symmetric and $\supp(\nu)\subsetneq \H$, there exists $w\in\H$ and $r>0$ such that $\nu[B_r(w)\cup\{-B_r(w)\}]=0$ and $0\notin B_r(w)$, where $B_r(w)=\{x\in\H:\|x-w\|\leq r\}$ denotes the $L^2$-ball centered at $w$ with radius $r$. We construct two well-defined characteristic functionals $\phi_1$ and $\phi_2$, such that $\phi_1-\phi_2$ is supported on $B_r(w)\cup\{-B_r(w)\}$. 
To achieve this, suppose $\{e_k\}_{k=1}^\infty$ is a complete orthogonal basis of $\H$. We take
$\phi_1(x)=\exp\big(-2^{-1}\sum_{k=1}^\infty k^{-2}\l x,e_k\r^2\big)$, for $x\in\H$, and
\begin{align*}
\phi_2(x)=\phi_1(x)+a\prod_{k=1}^\infty\bigg[1&-\frac{1}{3}\bigg(1-\frac{\pi k}{\sqrt{6}r}|\l x-w,e_k\r|\bigg)_+-\frac{1}{3}\bigg(1-\frac{\pi k}{\sqrt{6}r}|\l x+w,e_k\r|\bigg)_+\bigg]^{-k^{-2}}-a\,,
\end{align*}
where $(c)_{+}=\max\{c,0\}$ for $c\in\mathbb R$, and $a>0$ is a constant to be specified. 
We use the Minlos--Sazonov theorem (see Theorem~\ref{thm:ms} in the supplementary material) to verify that $\phi_2$ is a well-defined characteristic functional. We show that $\phi_1-\phi_2\neq0$ is supported on $B_r(w)\cup\{-B_r(w)\}$, so that $\int_\H|\phi_1(w)-\phi_2(w)|^2\nu(\d w)=0$,
which contradicts the fact that $\nu$ is characteristic.
\medskip

A direct consequence of Theorem~\ref{thm:iff:ed} is that 
the energy distance $\ed_\nu$ in~\eqref{ed}  is a semimetric on the space $\mathcal P(\H)$ of all Borel  probability measures on the Hilbert space $\H$, if $\nu$ is symmetric and has full support. It can further be easily  verified that $\ed_\nu$ satisfies the triangle inequality, so it is a well-defined metric on the space $\mathcal P(\H)$. The following corollary provides a formal statement; see Section~\ref{app:cor:4.11} of the supplementary material for more details  of the  proof.

\begin{corollary}\label{cor:4.11}
Suppose $\nu$ is a reference probability measure. 
Then, equipped with the metric $\ed_\nu$ in Definition~\ref{def:ed}, $\mathcal P(\H)$ is a metric space.

\end{corollary}

\subsection{Gaussian and Laplace reference probability measure}\label{sec:gaussian}


We provide two examples of characteristic reference probability measure. The first example consists of a class of Gaussian measures. For the second example, we define the Laplace measure  on an infinite-dimensional Hilbert space, and show that it is characteristic under certain conditions. 

Recall that a Borel measure $\nu $ on $\H$ is Gaussian if for each 
$x\in\H$, the function $\nu\{y\in\H:\l x,y\r\leq\cdot\}$ is the cumulative distribution function of a normal distribution with mean mean 
$\mu_x \in \mathbb{R} $ and variance $\sigma_x^2 > 0$
\cite[see, for example, Definition~2.6 in][]{kuo1975}.
The characteristic functional of a  Gaussian measure $\nu$  is given by 
\begin{align}\label{key}
\phi(x)=\exp\bigg(\i\l \mu,x\r-\frac{1}{2}\l C x,x\r\bigg)\,
\end{align}
 \citep[see, for example,][]{prokhorov1956},
where $\mu$ and $C$ are its mean and covariance operators:
\begin{equation}
   \label{key1} 
   \begin{split}
       \l \mu,x\r &=\int_\H  \l w,x\r\nu(\d w) ~~~~~~~~~~ ( x\in \H )~, \\
       \l C x,y\r& =\int_\H \l x,w-\mu\r\l y,w-\mu\r\nu(\d w)~~~~~~~~ ( x,y\in \H )~,
   \end{split}
\end{equation}
respectively 
\citep[Definitions~2.1 and 2.3 in][]{kuo1975}. Conversely, a functional of the form \eqref{key}  is a characteristic functional 
of Gaussian process if 
$C$ is an $\mathcal S$-operator \citep[see, for example, Theorem~2.3 in][]{kuo1975}.

A mean-zero Gaussian measure is symmetric. In addition, if a Gaussian measure $\nu$ on $\H$ is non-degenerate, i.e., its covariance operator is strictly positive definite, then $\supp(\nu)=\H$
\citep[see, for example, Proposition~1.25 in][]{daprato2006}. Therefore, as a consequence of Theorem~\ref{thm:iff:ed}, we have the following corollary, which states that for mean-zero Gaussian measures on $\H$, non-degeneracy and being characteristic are equivalent.

\begin{corollary}\label{lem:gaussian}

A mean-zero Gaussian measure on $\H$ is characteristic, if and only if it is non-degenerate.

\end{corollary}

It is known that, for finite dimensional Euclidean space $\mathbb R^p$, both Gaussian and Laplace measures are characteristic (see, for example, \citealp{fukumizu2007,sriperumbudur2010}). Due to the heavy-tails of  the characteristic function compared to the exponential decay of the Gaussian characteristic function, Laplace distributions have gained better performance in many statistical problems (see, for example, \citealp{delaigle2008,vandermaaten2008}). In the following, we define a Laplace probability measure on an infinite dimensional Hilbert space $\H$, which to our knowledge is new in the literature. 

To be precise, we call 
 a Borel measure  $\nu$    Laplace probability measure on $\H$,  if 
 for all  $x\in\H $, the function $\nu\{y\in\H:\l x,y\r\leq\cdot\}$ is the cumulative distribution function of a Laplace distribution on $\mathbb R$ with mean $\mu_x\in\mathbb R$ and variance $2b_x^2>0$ (depending on $x$). The mean $\mu$ and covariance operator $C$ of $\nu$ are defined by
 \eqref{key1}.
The following theorem establishes the characteristic functional of the Laplace measure $\nu$ on $\H$ in terms of its mean and covariance operator. Its proof uses the Minlos–Sazonov theorem 
and 
is given in Section~\ref{app:thm:laplace} of the supplementary material.

\begin{theorem}
\label{thm:laplace}
Suppose $\nu$ is a Laplace measure on $\H$ with mean operator $\mu$ and covariance operator $C$. Then, the functional
\begin{align}\label{phiLap}
\phi(x)=\exp(\i\l \mu,x\r)\bigg(1+\frac{1}{2}\l C x,x\r\bigg)^{-1}\,,\qquad x\in\H
\end{align}
is the characteristic functional of $\nu$ if and only if $C$ is an $\mathcal S$-operator. In addition, $\l \mu,x\r=\mu_x$ and $\l Cx,x\r=2b_x^2$, for $x\in\H$.

\end{theorem}

Recall that the characteristic function of multivariate Laplace distribution on $\mathbb R^p$ with mean $\mu\in\mathbb R^p$ and covariance matrix $\Sigma\in\mathbb R^{p\times p}$ is given by $\phi_p(x)=\exp(i\mu\trans x)\big(1+x\trans\Sigma x/2\big)^{-1}$, for $x\in\mathbb R^p$; see, for example, \cite{kotz2001}. Therefore, $\phi_p$ can be viewed as the $p$-dimensional restriction of the characteristic functional $\phi$ in \eqref{phiLap}.

The following proposition provides a sufficient condition for a Laplace measure to have full support on $\H$, which  implies that it is characteristic. A proof is given in Section~\ref{app:prop:support} of the supplementary material.


\begin{proposition}\label{prop:support}
Suppose $\nu$ is a mean-zero symmetric Laplace measure on $\H$ with covariance operator $C$. Then, $\supp(\nu)=\H$ if
the eigenvalues $\{\lambda_k\}_{k\geq1}$ of $C$ satisfy that, for any fixed $r\in\mathbb R$ and real-valued sequence $\{w_k\}_{k\geq1}$ such that $\sum_{k=1}^\infty w_k^2<\infty$,
\begin{align}\label{pp2}
\limsup_{M\to\infty}
\frac{\P\big[\sum_{k=1}^M \{\sqrt{\lambda_k}(L_{k}+L'_{k})-w_k\}^2\leq r^2\big]}{\P\big\{\sum_{k=1}^M (\sqrt{\lambda_k}L_{k}-w_k)^2\leq r^2\big\}} <\infty \,,
\end{align}
where $(L_{1},\ldots,L_{M})\trans,(L_{1}',\ldots,L_{M}')\trans\in\mathbb R^M$ are i.i.d.~mean-zero Laplace random vectors with identity covariance matrices.
\end{proposition}


\begin{remark} ~\\
(a)   If $w_k=0$ for all $k\geq0$, condition \eqref{pp2} is satisfied. In fact the ratio of probabilities in 
\eqref{pp2} is less or equal than $1$ for any $M\geq1$, This follows from  Anderson's inequality (see, for example, \cite{gardner}, p.~377), and the fact that the multivariate Laplace distribution has a unimodal probability density function, which is symmetric around the origin. However, in general condition \eqref{pp2} is not so easy to prove.
\medskip

(b) In  Section~\ref{sec:simu} we develop tests 
for homogeneity and independence using the new energy distances. It turns out that procedures using  Laplace reference measures have a better   finite-sample performance 
compared to the ones using Gaussian measures.
\end{remark}

\begin{example}\label{exp:2.10}
Consider $\H =L^2([0,1])$. Let $c:[0,1]^2\to\mathbb R$ be a continuous positive definite kernel. Mercer's theorem \citep{mercer} guarantees that the integral operator $C(x)=\int_0^1c(s,\cdot) x(s)\d s$, $x\in L^2([0,1])$, is an $\S$-operator, with finite trace $\tr(C)=\int_0^1c(s,s)\d s$. In this case, the characteristic functionals of the mean-zero Gaussian and Laplace measure with covariance operator $C$ are respectively given by 
\begin{align*}
\phi_{G}(w)
&=\exp\bigg[-\frac{1}{2}\int_0^1\int_0^1c(s,t)w(s)w(t) \d s\d t\bigg];\
~~\phi_{L}(w)=\bigg[1+\frac{1}{2}\int_0^1\int_0^1c(s,t)w(s)w(t)\d s\d t\bigg]^{-1}.
\end{align*}


\end{example}

\subsection{Energy distances and  maximum mean discrepancy}\label{sec:rkhs}

We conclude this section with a discussion on the connection between the energy distance in Definition~\ref{def:ed} and the maximum mean discrepancy (e.g.~\cite{sriperumbudur2010,gretton2012}) over a functional space in a reproducing kernel Hilbert space (RKHS).

Recall that for a general class $\mathcal F$ of functionals $f:\H\to\mathbb R$,  the maximum mean discrepancy between $\mathbb P,\mathbb Q\in\mathcal P(\H)$ over the space $\mathcal F$ is given by 
\begin{align}\label{mmdnu}
\mmd(\mathbb P,\mathbb Q;\F)=\sup_{f\in\F}\big|\E_{\mathbb P} f(Z)-\E_{\mathbb Q}f(Z)\big|=\sup_{f\in\F}\bigg|\int_{\H} f\, \d(\mathbb P-\mathbb Q)\bigg|\,.
\end{align}

Let $\nu$ be the reference probability measure of the energy distance $\ed_\nu$ in Definition~\ref{def:ed}. For $\mathbb P,\mathbb Q\in\H$, we will establish the equivalence between $\ed_\nu(\mathbb P,\mathbb Q)$ and $\mmd(\mathbb P,\mathbb Q;\mathcal F)$ defining an appropriate  functional class $\mathcal F$ that depends on the reference probability measure $\nu$.
For this purpose we define
\begin{align}\label{d}
d_\nu(x,y)=1-\phi_{\nu}(x-y)=1-\int_{\H}\cos\l x-y,w\r\,\nu(\d w)\,,\qquad x,y\in\H\,,
\end{align}
and note that  $d_\nu$ is translation-invariant.
In addition, $\supp(\nu)=\H$ if and only if $d_\nu$ is such that $(\H,d_\nu)$ is a semimetric space, that is: (a) $d_\nu(x,y)=d_\nu(y,x)$ for $x,y\in\H$; (b) $d_\nu(x,y)=0$ if and only if $x=y$. Furthermore, $d_\nu$ is {\em conditionally negative definite}: for any $\ell\geq2$, $z_1,\ldots,z_\ell\in\H$ and $ a _1,\ldots, a _\ell\in\mathbb R$ such that $\sum_{j=1}^\ell a _j=0$,
\begin{align}\label{nt}
\sum_{j=1}^\ell\sum_{k=1}^\ell a _j a _kd_\nu(z_j,z_k)
=-\int_\H\Big|\sum_{j=1}^\ell a_j\exp(\i\l z_j,w\r)\Big|^2\,\nu(\d w)\leq0\,.
\end{align}
This shows that $(\H,d_\nu)$ is a {\em semimetric space of negative-type}. By Lemma~3.2.1 in \cite{berg1984}
the kernel 
\begin{align}\label{knu2}
K_\nu(x,y)=d_\nu(x,0)+d_\nu(y,0)-d_\nu(x,y)=1-\phi_{\nu}(x)-\phi_{\nu}(y)+\phi_{\nu}(x-y)\,,\qquad x,y\in\H\,
\end{align}
is  symmetric and  positive definite.
By the Moore-Aronszajn theorem \citep{aronszajn1950}, there exists a unique RKHS $\H_{K_\nu}$ with  a corresponding norm $\|\cdot\|_{\H_{K_\nu}}$ such that $K_\nu$ in \eqref{knu2} is its reproducing kernel. Furthermore,
\begin{align*}
d_\nu(x,y)=\|K_\nu(x,\cdot)-K_\nu(y,\cdot)\|_{\H_{K_\nu}}^2=K_\nu(x,x)+K_\nu(y,y)-2K_\nu(x,y)\,,\qquad x,y\in\H\,.
\end{align*}

The following proposition proved in Section~\ref{proof:prop:mmded} of the supplementary material shows that $\ed_\nu$ in \eqref{ed} is equivalent to the maximum mean discrepancy \eqref{mmdnu} over the unit ball in the RKHS $\H_{K_\nu}$ with reproducing kernel $K_\nu$ in \eqref{knu2}. 

\begin{proposition}\label{prop:mmded}

Suppose that $\nu$ is a reference measure. Then, for any $\mathbb P,\mathbb Q\in\mathcal P(\H)$,
$$\mmd(\mathbb P,\mathbb Q;\H_{K_\nu,1})=\ed_\nu(\mathbb P,\mathbb Q),
$$
where $\H_{K_\nu,1}:=\{f\in\H_{K_\nu}:\|f\|_{\H_{K_\nu}}\leq1\}$.

\end{proposition}

\section{Distance covariance and measures of independence}\label{sec:dcov}


In this section we  define a class of distance covariances between $\H$-valued random variables that do not require any moment conditions. Let $X$ and $Y$ be $\H$-valued random variables, with joint distribution $\mathbb P^{X,Y}\in\mathcal P(\H^2)$, and let $\mathbb P^X,\mathbb P^Y\in\mathcal P(\H)$ be their corresponding marginal distributions. 
A direct extension of Theorem~\ref{thm:bijection} implies that 
the independence of $X$ and $Y$, that is $\mathbb P^{X,Y}=\mathbb P^X \otimes\mathbb P^Y$,  is equivalent to $\phi_{X,Y}=\phi_X\otimes\phi_Y$, that is  the joint characteristic functional is the tensor product of its  marginals. 

For a Borel probability measure $\nu$, we define the distance covariance between $X$ and $Y$, by integrating the squared modulus of the difference between the joint characteristic functional $\phi_{X,Y}$ and the tensor product marginals $\phi_X\otimes\phi_Y$ with respect to the product measure  
$\nu\otimes\nu$ on $\H^2$.

\begin{definition}
\label{def:dcov}
The (squared) distance covariance between the  $\H$-valued random variables $X$ and $Y$ with respect to 
$\nu\in\mathcal P(\H)$ is defined by
\begin{align}\label{dcov_def}
\dcov^2_\nu(X,Y)
&=\int_{\H^2}|\phi_{X,Y}(w_1,w_2)-\phi_{X}(w_1)\phi_{Y}(w_2)|^2\,\nu(\d w_1)\,\nu(\d w_2)\,.
\end{align}
\end{definition}

The distance covariance  can be viewed as the energy distance on $\mathcal P(\H^2)$,
between the  measure $\mathbb P^{X,Y}$ and the product measure $\mathbb P^X\otimes\mathbb P^Y$, where the reference probability measure is taken to be the product measure $\nu\otimes\nu$, that is 
$$
\dcov_\nu^2(X,Y)=\ed_{\nu\otimes\nu}^2(\mathbb P^{X,Y},\mathbb P^X\otimes\mathbb P^Y).
$$
Our next result  states an alternative representation  for the  distance covariance in terms of $\phi_{\nu}$ and $d_\nu$, and is proved in Section~\ref{app:thm:dcov} of the supplementary material. 

\begin{proposition}\label{thm:dcov}
Suppose $\nu\in\mathcal P(\H)$ is symmetric.
Let $(X',Y'),(X'',Y'')$ be  independent copies of $(X,Y)$. Then, the distance covariance in \eqref{dcov_def} satisfies
\begin{align}\label{dcov}
\dcov^2_\nu(X,Y)&=\E[\phi_{\nu}(X-X')\,\phi_{\nu}(Y-Y')]+\E\phi_{\nu}(X-X')\,\E\phi_{\nu}(Y-Y')-2\E[\phi_{\nu}(X-X')\phi_{\nu}(Y-Y'')]\notag\\
&=\E\big\{[\phi_{\nu}(X-X')+\E\phi_{\nu}(X-X')-2\E_X\phi_{\nu}(X-X')]\notag\\
&\qquad\times[\phi_{\nu}(Y-Y')+\E\phi_{\nu}(Y-Y')-2\E_Y\phi_{\nu}(Y-Y')]\big\}\,.
\end{align}
In addition, the equations hold by replacing $\phi_\nu(\square-\triangle)$ by $d_\nu(\square,\triangle)$, for $d_\nu$ in \eqref{d}.

\end{proposition}

Note that the definition of the  distance covariance  does not require any moment condition on $\mathbb P^{X,Y}$, in contrast previous works such as \cite{lyons2013}. 
The second formula in \eqref{dcov} is useful for deriving properties of distance variance and correlation; see the discussion below.

Parallel to Theorem~\ref{thm:iff:ed}, the following result  states a sufficient and necessary condition on $\nu$ such that $\dcov_\nu^2(X,Y)=0$ is equivalent to $X$ and $Y$ being independent, and is proved in Section~\ref{app:thm:sufficient} of the supplementary material.

\begin{theorem}\label{thm:iff:dcov}
Suppose  that the probability measure $\nu$ is symmetric.
Then, $\supp(\nu)=\H$ 
if and only if  the following statement is true for all  $\H$-valued random variables $X$ and $Y$:
\begin{align}\label{characteristic}
\text{$\dcov_\nu^2(X,Y)=0$ if and only if $X$ and $Y$ are independent.}
\end{align}

\end{theorem}


\begin{remark}
\label{prop:hsic}
For the kernel $K_\nu$ in \eqref{knu2} and the corresponding RKHS $\H_{K_\nu}$, let $\H_{K_\nu}^{\otimes2}=\H_{K_\nu}\otimes \H_{K_\nu}$ be the tensor product RKHS with reproducing kernel $K_\nu(x,x')K_\nu(y,y')$, for $x,x',y,y'\in\H$ and  corresponding norm $\|\cdot\|_{\H_{K_\nu}^{\otimes2}}$.
The  
Hilbert-Schmidt independence criterion (HSIC) is defined by
\begin{align*}
\hsic_\nu(X,Y)
=\sup_{f\in\H_{K_\nu,1}^{\otimes2}}\bigg|\int_{\H^2} f \d(\mathbb P^{X,Y}-\mathbb P^X\otimes \mathbb P^Y)\bigg|\,,
\end{align*}
where $\H_{K_\nu,1}^{\otimes2}=\{f\in\H_{K_\nu}^{\otimes2}:\|f\|_{\H_{K_\nu}^{\otimes2}}\leq 1\}$ \citep[see, for example,][]{gretton,gretton2007,smola2007}. It is shown in Section~\ref{app:proof:hsic} of the supplementary material that 
$$
\hsic_\nu^2(X,Y)=\dcov_\nu^2(X,Y).
$$
\end{remark}


We conclude this section by defining a class of distance variance and correlation, which do not require any moment conditions.
For a reference probability measure $\nu$ on $\H$, and two $\H$-valued random variables $X$ and $Y$, the distance variance of an $\H$-valued random variable $X$ is defined by $\dvar_\nu(X)=\dcov_\nu(X,X)$, and the distance correlation between  $X$ and $Y$ is defined by
\begin{align}\label{dcor}
\dcor_\nu^2(X,Y)=\bBigg@{3}\{\begin{array}{cl}
\displaystyle\frac{\dcov_\nu^2(X,Y)}{\dvar_\nu(X)\,\dvar_\nu(Y)}&\  \text{ if }\,\dvar_\nu(X)\dvar_\nu(Y)\neq0\,;\\
0 &\  \text{ if }\,\dvar_\nu(X)\dvar_\nu(Y)=0\,.\\
\end{array}
\end{align}
The following two propositions state important properties of distance variance and correlation, and are  proved in Section~\ref{app:thm:dvar} and \ref{app:thm:prop_dcor} of the supplementary material, respectively.

\begin{proposition}
~\label{thm:dvar}
Suppose that $\nu$  is an reference probability measure, and $X,X',X''$ are independent and identically distributed $\H$-valued random variables. Then,
\begin{align*}
\dvar_\nu^2(X)&=\E [\phi_{\nu}(X-X')]^2+[\E \phi_{\nu}(X-X')]^2-2\E[\phi_{\nu}(X-X')\phi_{\nu}(X-X'')]\\
&=\E[\phi_{\nu}(X-X')+\E \phi_{\nu}(X-X')-2\E_X\phi_{\nu}(X-X')]^2.
\end{align*}
Furthermore, $\dvar_\nu^2(X)=0$ if and only if $X  $ is  constant.

\end{proposition}

Note that the first equation for $\dvar_\nu(X)$ in Proposition~\ref{thm:dvar} implies that $\dvar_\nu^2(X)\leq4$ for any $\mathbb P^X\in\mathcal P(\H)$ and reference probability measure $\nu$. Therefore, the definition of $\dcov_\nu^2(X)$ does not require any moment condition on $\mathbb P^X$.

\begin{proposition}
\label{thm:prop_dcor}

Suppose $\nu$  is a reference probability measure.
Then,
\begin{enumerate}[nolistsep,label={\rm(\roman*)}]
\item $0\leq\dcor_\nu^2(X,Y)\leq 1$; 

\item $\dcor_\nu^2(X,Y)=0$ if and only if $X$ and $Y$ are independent.

\item $\dcor_\nu^2(X,Y)=1$ if and only if $\phi_{\nu}(X-X')=a \phi_{\nu}(Y-Y')$ a.s.~for some constant $a\in\mathbb R$, where $(X,Y)$ and $(X',Y')$ are i.i.d.
\end{enumerate}
\end{proposition}

\section{Statistical inference}\label{sec:minimax}

In this section we use the new distances for statistical inference 
 of $\H$-valued data. We start by introducing an empirical  analog of the characteristic functional introduced in \eqref{phixy}.
For  a random sample of $\H$-valued random variables $\T_n= \{ (X_1,Y_1),\ldots,(X_n,Y_n) \} $ with distribution $\mathbb{P}^{X,Y}$, the empirical joint characteristic functional is defined by 
\begin{align}\label{hatphixy}
\hat\phi_{X,Y}(w_1,w_2)&=\frac{1}{n}\sum_{j=1}^n\exp(\i\l X_j,w_1\r+\i\l Y_j,w_2\r)\,,\qquad w_1,w_2\in\H\,.
\end{align}
The corresponding empirical marginal characteristic functionals are, respectively,
\begin{align}\label{empcha}
\hat\phi_X(w)=\frac{1}{n}\sum_{j=1}^n\exp(\i\l X_j,w\r)\,,\qquad \hat\phi_Y(w)=\frac{1}{n}\sum_{j=1}^n\exp(\i\l Y_j,w\r)\,,\qquad w\in\H\,.
\end{align}
In view of Theorem~\ref{thm:iff:ed}, let 
\begin{align}
    \label{det51}
\V=\{\nu\in\mathcal P(\H):\nu\text{ is symmetric},\,\supp(\nu)=\H\}
\end{align}
denote the set of reference  probability measures on $\H$.
Consider the problem of testing a null hypothesis $H_0:  \mathbb{P}^{X,Y} \in \mathcal P_0 $ against an alternative hypothesis $H_1: \mathbb{P}^{X,Y} \in \mathcal P_1$  for some disjoint subsets  $\mathcal P_0,\mathcal P_1\subset\mathcal P(\H^2)$.
In the following discussion the definition of $\mathcal P_0,\mathcal P_1$ will always be clear from the context.
Let $\Psi$ denote the set of all non-randomized tests and define for a 
given a nominal level $\alpha \in (0,1)$  
\begin{align}\label{Psialpha}
\Psi(\alpha)=\Big\{\psi\in\Psi:\sup_{\mathbb{P}^{X,Y} \in \mathcal P_0}\mathbb{P}^{X,Y} (\psi=1)\leq\alpha\Big\}\,
\end{align}
as the set of all tests with a uniformly controlled
 type-I error.

\subsection{Test of independence}\label{sec:permu:indep}

For inference on independence, the data $\T_n=\{ (X_1,Y_1),\ldots,(X_n,Y_n) \} $ consists of i.i.d.~copies of the $\H^2$-valued random variable $(X,Y)$. The goal is to test whether $X$ and $Y$ are independent, that is
\begin{align}
    \label{independ}
H_0:\mathbb P^{X,Y}= \mathbb P^X\otimes \mathbb P^Y
\text{~~~versus ~~~~}  H_1:\mathbb P^{X,Y} \not = \mathbb P^X\otimes \mathbb P^Y.
\end{align}
In view of Theorem~\ref{thm:iff:dcov}, this is equivalent to testing the hypotheses $H_0: \dcov_\nu^2(\mathbb P^X,\mathbb P^Y)=0$ 
versus $H_1: \dcov_\nu^2(\mathbb P^X,\mathbb P^Y) \not =0$ for some $\nu\in\V$. Then, a reasonable test statistic is a consistent estimator of the distance covariance $\dcov_\nu^2(\mathbb P^X,\mathbb P^Y)$, defined by the integrated squared modulus of the difference between the empirical characteristic functionals 
$\hat\phi_{X,Y}$ and $\hat\phi_X\otimes\hat\phi_Y$, that is 
\begin{align}\label{dcovn}
\hat\dcov_\nu^2(X,Y)=\int_{\H^2}\big|\hat\phi_{X,Y}(w_1,w_2)-\hat\phi_{X}(w_1)\hat\phi_{Y}(w_2)\big|^2\,\nu(\d w_1)\,\nu(\d w_2)\,.
\end{align}
The first result in this section gives a useful formula for computing the empirical distance covariance in \eqref{dcovn}, and  is proved in Section~\ref{app:prop:dcovn} of the supplementary material.

\begin{proposition}\label{prop:dcovn}
Suppose Assumption~\ref{a:sym} is satisfied. The empirical distance covariance in \eqref{dcovn} is given by 
$$
\hat\dcov_\nu^2(X,Y)=n^{-2}\sum_{k,\ell=1}^n\hat V_{k,\ell}^{(X,\nu)}\hat V_{k,\ell}^{(Y,\nu)},
$$
where
\begin{align}\label{hatv}
\hat V_{k,\ell}^{(X,\nu)}&=\phi_{\nu}(X_k-X_\ell)-\frac{1}{n}\sum_{i=1}^n\phi_{\nu}(X_{i}-X_\ell)-\frac{1}{n}\sum_{j=1}^n\phi_{\nu}(X_k-X_{j})+\frac{1}{n^2}\sum_{i,j=1}^n\phi_{\nu}(X_{i}-X_{j})\,;\notag\\
\hat V_{k,\ell}^{(Y,\nu)}&=\phi_{\nu}(Y_k-Y_\ell)-\frac{1}{n}\sum_{i=1}^n\phi_{\nu}(Y_{i}-Y_\ell)-\frac{1}{n}\sum_{j=1}^n\phi_{\nu}(Y_k-Y_{j})+\frac{1}{n^2}\sum_{i,j=1}^n\phi_{\nu}(Y_{i}-Y_{j})\,.
\end{align}
\end{proposition}

The following proposition states the asymptotic properties of the empirical distance covariance defined in \eqref{dcovn}, and is proved in Section~\ref{app:thm:converge} of the supplementary material.

\begin{proposition}
\label{thm:converge}

The empirical distance covariance $\hat\dcov_\nu^2(X,Y)$ in \eqref{dcovn} satisfies 

\begin{enumerate}[nolistsep,label={\rm(\roman*)},wide]
\item $\hat\dcov_\nu^2(X,Y)\to\dcov_\nu^2(X,Y)$ almost surely as $n\to\infty$;

\item If $\H$-valued random variables $X$ and $Y$ are independent, then as $n\to\infty$,
\begin{align}\label{ev}
n\,\hat\dcov_\nu^2(X,Y)\longrightarrow\E \phi_{\nu}(X-X')\,\E \phi_{\nu}(Y-Y')+\sum_{j=1}^\infty\lambda_{\nu,j}(\zeta_j^2-1)\qquad\text{in law}\,,
\end{align}
where $\{\zeta_j\}_{j\geq1}$ 
is a sequence of independent standard normal distributed random variables,  $\{\lambda_{\nu,j}\}_{j\geq1}$ are  the eigenvalues of the operator $S_\nu$, which is defined  for $\eta:\H^2\to\mathbb R$ with  $\E\{\eta^2(X,Y)\}<\infty$ by 
\begin{align*}
( S_\nu\eta) (x,y)=\E\big[\{\phi_{\nu}(x-X)-\E \phi_{\nu}(x-X)\}\{\phi_{\nu}(y-Y)-\E \phi_{\nu}(y-Y)\}\eta(X,Y)\big]\,   ~~( x,y\in\H ) .
\end{align*}

\end{enumerate}

\end{proposition}

Given the data $\T_n=\{ (X_1,Y_1),\ldots,(X_n,Y_n)\}$, for the purpose of testing whether $X$ and $Y$ are independent, the null hypothesis $\mathbb P^{X,Y}=\mathbb P^X\otimes\mathbb P^Y$ should be rejected if
$\hat\dcov_\nu^2(X,Y)$
exceeds a threshold, which can be the $(1-\alpha)$-quantile of the limit distribution in \eqref{ev}. However, this distribution  involves unknown quantities such as the eigenvalues $\lambda_{\nu,j}$, 
which are often intractable and difficult to estimate.

This difficulty can be resolved by applying permutation tests, whose validity of controlled type-I error depends solely on the assumption of exchangeability under the null hypothesis. To be specific, for a positive integer $B$, let $\Pi_1,\Pi_2,\ldots,\Pi_B$ be i.i.d.~uniform random permutations of $(1,2,\ldots,n)$ that are 
independent of the data.
For each $1\leq b\leq B$, we obtain the permuted samples $\T_n(\Pi_b)=\{(X_j,Y_{\Pi_b(j)})\}_{1\leq j\leq n}$, and the corresponding permuted test statistic is then given by
\begin{align}\label{hatdb}
\hat\dcov_\nu^2(X,Y,\Pi_b)=\frac{1}{n^2}\sum_{k,\ell=1}^n\hat V_{k,\ell}^{(X,\nu)}\hat V_{\Pi_b(k),\Pi_b(\ell)}^{(Y,\nu)}\,,
\end{align}
where $\hat V_{k,\ell}^{(X,\nu)}$ and $\hat V_{k,\ell}^{(Y,\nu)}$ are defined in \eqref{hatv}.
Then, at nominal level $\alpha$, the decision rule to reject the null hypothesis of independence is defined by
\begin{align}\label{hatpsi}
\hat\psi_{\nu,n}(\alpha)=\one\bigg\{\frac{1}{B}\sum_{b=1}^B\one\big\{\hat\dcov_\nu^2(X,Y)\leq\hat\dcov_\nu^2(X,Y,\Pi_b)\big\}>1-\alpha\bigg\}\,.
\end{align}
Under the null hypothesis that $X$ and $Y$ are independent, it is seen that the original sample $\T_n$ and the permuted sample $\T_n(\Pi_1),\ldots,\T_n(\Pi_B)$ have the same distribution, such that the test statistics $\hat\dcov_\nu^2(X,Y)$ and $\{\hat\dcov_\nu^2(X,Y,\Pi_j)\}_{1\leq j\leq B}$ are identically distributed. This implies that $\hat\psi_{\nu,n}$ defined in \eqref{hatpsi} achieves uniform type-I error control at nominal level $\alpha$, that is  $\hat\psi_{\nu,n}(\alpha)\in\Psi(\alpha)$.

\begin{remark}
\label{rem:computation}

The test statistic $\hat\dcov_\nu^2(X,Y)$ and its permuted version $\hat\dcov_\nu^2(X,Y,\Pi_b)$ can be computed via the following steps, where the computation is efficient using matrix operations:
\begin{enumerate}[nolistsep,label={(\arabic*)}]
\item Compute matrices $M^{(X)},M^{(Y)}\in\mathbb R^{n,n}$, with entries 
$$
M^{(X)}_{k,\ell}=\phi_{\nu}(X_k-X_\ell) ~\text{  and }  ~M^{(Y)}_{k,\ell}=\phi_{\nu}(Y_k-Y_\ell)~~~\quad(1\leq k,\ell\leq n).
$$

\item For $Z=X$ and $Y$, define $1_n=(1, \ldots , 1)^\top \in \mathbb{R}^n$ and compute 
$$\hat V^{(Z)}=M^{(Z)}-n^{-1}(M^{(Z)}1_n1_n\trans+1_n1_n\trans M^{(Z)})+n^21_n1_n^\top M^{(Z)}1_n1_n\trans.
$$
\item Compute 
\begin{align*}
\hat\dcov_\nu^2(X,Y)& =n^{-2}1_n\trans(\hat V^{(X)}\circ\hat V^{(Y)})1_n ,  \\
\hat\dcov_\nu^2(X,Y,\Pi_b) &=n^{-2}1_n\trans\{\hat V^{(X,\nu)}\circ (P_{\Pi_b}\hat V^{(Y,\nu)} P_{\Pi_b}\trans)\}1_n,
\end{align*}
where $P_{\Pi_b}$ is the row-wise permutation matrix that corresponds to the permutation $\Pi_b$.

\end{enumerate}

\end{remark}

To analyze the power   of the independence test $\hat\psi_{\nu,n}$ in \eqref{hatpsi}, consider the set of local alternatives that are separated from the null by $\rho>0$, where the separation is measured via the squared distance covarianace between $X$ and $Y$,
that is
\begin{align}\label{fnurhoa}
\mathcal F_\nu(\rho)=\bigg\{\mathbb P^{X,Y}\in\mathcal P(\H^2):\int_{\H^2}|\phi_{X,Y}-\phi_X\otimes\phi_Y|^2\d(\nu\otimes\nu)\geq\rho\bigg\}\,.
\end{align}
In order to investigate the minimax optimality of the permutation test $\hat\psi_{\nu,n}$ in \eqref{hatpsi}, we follow 
\cite{baraud2002} and define for  $\psi\in\Psi(\alpha)$
and $\beta \in (0,1) $
\begin{align}\label{seprate}
\rho(\psi,\beta)=\inf\Big\{\rho>0:\sup_{\nu\in\V}\sup_{\mathcal F_\nu(\rho)}\P(\psi=0)\leq\beta\Big\}\,
\end{align}
as 
the uniform separation rate with respect to $\V$ and $\mathcal F_\nu(\rho)$.
Note that $\rho(\psi,\beta)$
is  the smallest value that separates the null and the local alternatives, such that the type-II error is  controlled at level $\beta$
uniformly over the set $\mathcal F_\nu(\rho)$ and the set $\V$ of all reference probability measures. 
The following result  establishes an  upper bound for the uniform separation rate of the test \eqref{hatpsi} and  is proved in Section~\ref{app:thm:upper:ind} of the supplementary material.
\begin{theorem}
\label{thm:upper:ind}
Suppose $0<\alpha+\beta<1$, $\beta<1/e$, $n\geq4$ and $B\geq2\alpha^{-2}\log(6/\beta)$.
For the  the test $\hat\psi_{\nu,n}(\alpha)$ defined in \eqref{hatpsi}, it holds  for an absolute constant $c>0$ that
\begin{align*}
\rho\{\hat\psi_{\nu,n}(\alpha),\beta\}\leq c n^{-1/2}\{\log(\alpha^{-1})+\log(\beta^{-1})\}\,.
\end{align*}
\end{theorem}

Theorem~\ref{thm:upper:ind} states that the uniform separation rate of the permutation independence test 
\eqref{hatpsi}
is given by the parametric rate $n^{-1/2}$ for given bounds
$\alpha$ and $\beta$ on the type I and II errors. In the following 
discussion we will show that 
this separation rate is 
the best  uniform separation rate
over $\Psi(\alpha)$,
 which is defined by 
\begin{align}\label{rho*}
\rho^*(\alpha,\beta)=\inf_{\psi_\alpha\in\Psi(\alpha)}\rho(\psi_\alpha,\beta)\,.
\end{align}
More precisely, the final result of this section
shows that for independence testing 
$\rho^*(\alpha,\beta)$ can  bounded from below  by a quantity of order  $n^{-1/2}$. As this  matches the uniform separation rate of the permutation test \eqref{hatpsi}
this test is minimax optimal.  Its proof leverages the Radon-Nikodym derivatives between Gaussian measures  and is given in Section~\ref{app:thm:lower:ind} of the supplementary material.

\begin{theorem}
\label{thm:lower:ind}
Suppose $\alpha,\beta\in(0,1)$ are fixed and satisfy $\alpha+\beta<1$. Then, there exists a constant $c>0$ independent of $n$ such that, for $n$ large enough, $\rho^*(\alpha,\beta)\geq cn^{-1/2}$.

\end{theorem}

\subsection{Two-sample test of equal distribution}\label{sec:twosample}

Suppose the data $\T_{n,2}=\{ X_1,\ldots,X_{n_1},Y_{1},\ldots,Y_{n_2} \}$ consist of independent $\H$-valued observations from two samples, where $X_1,\ldots,X_{n_1}\iidsim \mathbb P^X$ and  $Y_{1},\ldots,Y_{n_2}\iidsim\mathbb P^Y$. Our goal is to use the new energy distance to develop a test for  the hypotheses \begin{align} \label{twosample} 
H_0:\mathbb P^X=\mathbb P^Y \text{~~~versus ~~~~} H_1:\mathbb P^X \not =\mathbb P^Y.
\end{align}
By Theorem~\ref{thm:iff:ed}, for a characteristic reference probability measure $\nu\in\V$, the null hypothesis in  \eqref{twosample} is equivalent to  $\ed_\nu^2(\mathbb P^X,\mathbb P^Y)=0$. Therefore, a reasonable test statistic is obtained by 
replacing the characterstic functionals in \eqref{ed} by their empirical counterparts $\hat\phi_X$ and $\hat\phi_Y$ defined in \eqref{empcha}, that is 
\begin{align}\label{emed}
\ed_\nu^2(\hat{\mathbb P}^X,\hat{\mathbb P}^Y)=\int_\H|\hat\phi_X(w)-\hat\phi_Y(w)|^2\,\nu(\d w)\, .
\end{align}
Here, $ \hat{\mathbb P}^X$ and $\hat{\mathbb P}^Y $ denote the corresponding empirical measures.
The following proposition states two equivalent formulas for the empirical energy distance, and is proved in Section~\ref{app:prop:edhat} of the supplementary material.
\begin{proposition}\label{prop:edhat}
If  $\nu$ is symmetric, then the empirical energy distance \eqref{emed} 
is given by
\begin{align}\label{eed}
\ed^2_\nu(\hat{\mathbb P}^X,\hat{\mathbb P}^Y) & =\frac{1}{n_1^2}\sum_{j,k=1}^{n_1}\phi_{\nu}( X_{j}-X_k)+\frac{1}{n_2^2}\sum_{j,k=1}^{n_2}\phi_{\nu}( Y_j-Y_k)-\frac{2}{n_1n_2}\sum_{j=1}^{n_1}\sum_{k=1}^{n_2}\phi_{\nu}(X_j-Y_k)\\
& =\frac{n_1+n_2}{n_1n_2}\sum_{j,k=1}^{n_1+n_2}\phi_{\nu}(Z_j-Z_k)Q_{j,k}\,,
\label{new}
\end{align}
where 
$
(Z_1,\ldots,Z_{n_1},Z_{n_1+1},\ldots,Z_{n_1+n_2})=(X_1,\ldots,X_{n_1},Y_1,\ldots,Y_{n_2})$ and
\begin{align}
Q_{j,k}=\frac{\one\{1\leq j,k\leq n_1\}}{n_1}+\frac{\one\{n_1+1\leq j,k\leq n_1+n_2\}}{n_2}-\frac{1}{n_1+n_2}\,,\quad 1\leq j,k\leq n_1+n_2\,.\notag
\end{align}
\end{proposition}
The formula in \eqref{new} is useful for defining the permuted empirical energy distance. The next result  states the asymptotic properties of the empirical energy distance, and is proved in Section~\ref{app:thm:converge:ed} of the supplementary material.
\begin{proposition}
\label{thm:converge:ed}
The empirical energy distance $\ed_\nu^2(\hat{\mathbb P}^X,\hat{\mathbb P}^Y)$ defined in \eqref{ed} satisfies 
\begin{enumerate}[nolistsep,label={\rm(\roman*)},wide]
\item $\ed_\nu^2(\hat{\mathbb P}^X,\hat{\mathbb P}^Y)\to\ed_\nu^2 (\mathbb P^X,\mathbb P^Y)$ almost surely as $n\to\infty$;
\item Suppose $\dvar_\nu^2(X)>0$ and $n_1/n_2\to c_0>0$ as $n\to\infty$. If $\mathbb P^X=\mathbb P^Y$, then as $n_1,n_2\to\infty$,
\begin{align}\label{edlaw}
n_1\ed_\nu^2(\hat{\mathbb P}^X,\hat{\mathbb P}^Y)&
\longrightarrow-(1+c_0)\E \phi_{\nu}(X-X')+\sum_{\ell=1}^\infty\tau_{\nu,\ell}\{(\sqrt{c_0}\xi_\ell-\zeta_\ell)^2-1\}\qquad \text{in law}\,,
\end{align}
where $\{\xi_\ell\}_{\ell\geq1},\{\zeta_\ell\}_{\ell\geq1}$ 
are independent sequences of  independent standard normal distributed random variables, and the $\tau_{\nu,\ell}$ are the eigenvalues of the operator 
$R_\nu$ defined by 
$$
(R_\nu\eta )(x)=\E[\{\phi_{\nu}(x-X)-\E \phi_{\nu}(x-X)\}\eta(X)] \qquad(x\in\H )
$$
for $\eta:\H\to\mathbb R$ such that $\E\{\eta^2(X)\}<\infty$.
\end{enumerate}
\end{proposition}

For the purpose of defining a decision rule
for testing the hypotheses in \eqref{twosample}  we propose to reject the null hypothesis for large values of  $\ed_\nu^2(\hat{\mathbb P}^X,\hat{\mathbb P}^Y)$. Since the asymptotic distribution in \eqref{edlaw} is often intractable in practice, we apply permutation tests to find this critical value. To be specific, for a positive integer $B$, let $\Pi_1,\Pi_2,\ldots,\Pi_B$ denote $B$ i.i.d.~uniform random permutations of $(1,2,\ldots,n_1+n_2)$ independent of the original data $\T_{n,2}$. Recalling the definition of $\{Z_j\}$ in Proposition~\ref{prop:edhat}, for each $1\leq b\leq B$, let $\T_{n,2}(\Pi_b)=(Z_{\Pi_b(1)},\ldots,Z_{\Pi_b(n_1)},Z_{\Pi_b(n_1+1)}\ldots,Z_{\Pi_b(n_1+n_2)})$ denote the permuted sample under the permutation $\Pi_b$. In view of Proposition~\ref{prop:edhat}, define the permuted test statistic
\begin{align}\label{edtest}
&\ed_\nu^2(\hat{\mathbb P}^X,\hat{\mathbb P}^Y,\Pi_b)=\frac{n_1+n_2}{n_1n_2}\sum_{j,k=1}^{n_1+n_2}\phi_{\nu}(Z_{\Pi_b(j)}-Z_{\Pi_b(k)})Q_{\Pi_b(j),\Pi_b(k)}\,,
\end{align}
where $\{Q_{j,k}\}_{1\leq k,\ell\leq n_1+n_2}$ are defined in \eqref{new}.
Under the null hypothesis in \eqref{twosample}, the original test statistic $\ed_\nu^2(\hat{\mathbb P}^X,\hat{\mathbb P}^Y)$ and the permuted versions $\ed_\nu^2(\hat{\mathbb P}^X,\hat{\mathbb P}^Y,\Pi_1),\ldots,\ed_\nu^2(\hat{\mathbb P}^X,\hat{\mathbb P}^Y,\Pi_B)$ defined in \eqref{edtest} are identically distributed.
Then, the  test 
\begin{align}\label{hatpsi2}
\hat\psi_{\nu,n,2}(\alpha)=\one\bigg\{\frac{1}{B}\sum_{b=1}^B\one\big\{\ed_\nu^2(\hat{\mathbb P}^X,\hat{\mathbb P}^Y)\leq\ed_\nu^2(\hat{\mathbb P}^X,\hat{\mathbb P}^Y,\Pi_b)\big\}>1-\alpha\bigg\}\,
\end{align}
for the hypotheses \eqref{twosample}   satisfies $\hat\psi_{\nu,n,2}\in\Psi(\alpha)$, due to exchangeability.  
The final result of this section shows that the test $\hat\psi_{\nu,n,2}$  is minimax optimal. A proof can be found in 
Section~\ref{app:thm:upper:two} of the supplementary material.

\begin{theorem}\label{thm:upper:two}
Suppose $0<\alpha+\beta<1$, $c_1\leq n_1/ n_2\leq c_2$ for some constants $c_1,c_2>0$, and $\beta<1/e$, $B\geq2\alpha^{-2}\log(6/\beta)$. The uniform separation rate $\rho\{\hat\psi_{\nu,n,2}(\alpha),\beta\}$ of the two-sample test $\hat\psi_{\nu,n,2}(\alpha)$ in \eqref{hatpsi2} satisfies
\begin{align*}
\rho\{\hat\psi_{\nu,n,2}(\alpha),\beta\}
& =\inf\Big\{\rho>0:\sup_{\nu\in\V}\sup_{\mathcal F_\nu(\rho)}\P(\hat\psi_{\nu,n,2}(\alpha) =0) 
 \leq\beta\Big\} \\
& \leq c (n_1+n_2)^{-1/2}\{\log(\alpha^{-1})+{\log(\beta^{-1})} \}\,,
\end{align*}
where $\mathcal F_\nu(\rho)=\big\{\mathbb P^X,\mathbb P^Y\in\mathcal P(\H):\int_{\H}|\phi_{X}-\phi_{Y}|^2\d\nu\geq\rho\big\}$ and $c>0$ is an absolute constant.
Moreover, for given $\alpha,\beta$ such that $0<\alpha+\beta<1$, the minimax separation rate satisfies
$$
\rho^*(\alpha,\beta)=\inf_{\psi_\alpha\in\Psi(\alpha)}\rho(\psi_\alpha,\beta) \geq c'(n_1+n_2)^{-1/2}
$$
for some absolute constant $c'>0$.
\end{theorem}

\subsection{Inference in presence of measurement error}\label{sec:me}

In this section we develop a two-sample test  in the presence of measurement errors. Suppose the $\H$-valued random variables $X,Y$ are measured with errors:
\begin{align}\label{me}
\tilde X=X+U_{X},\qquad\tilde Y=Y+U_{Y}\,,
\end{align}
where $U_{X}$ and $U_Y$ are independent $\H$-valued random variables independent of $(X,Y)$ with characteristic functionals  $\phi_{U_X}$ and $\phi_{U_Y}$, respectively. Assume that $\phi_{U_X}$ and $\phi_{U_Y}$ do not vanish, that is, $\phi_{U_X}(w)\phi_{U_Y}(w)\neq0$, for any $w\in\H$. This assumption is satisfied when $\mathbb P^{U_X}$ and $\mathbb P^{U_Y}$ are mean-zero Gaussian measures (equation \eqref{key}) or mean-zero Laplace measures on $\H$ (Theorem~\ref{thm:laplace}), and is a common assumption for measurement errors in Euclidean space (see, for example, \citealp{carroll1988,meister2009}). Observing  that $\mathbb P^{\tilde X}$ (similar for $\mathbb P^{\tilde Y}$) is the convolution measure of $\mathbb P^X$ and $\mathbb P^{U_X}$, we have
\begin{align*}
&\ed_\nu^2(\mathbb P^{\tilde X},\mathbb P^{\tilde Y})=\int_\H|\phi_X(w)-\phi_Y(w)|^2|\phi_U(w)|^2\,\nu(\d w)\,,\qquad\text{when }\phi_{U_X}=\phi_{U_Y}=\phi_U\,;\\
&\dcov_\nu^2(\tilde X,\tilde Y)=\int_{\H^2}|\phi_{X,Y}(w_1,w_2)-\phi_{X}(w_1)\phi_{Y}(w_2)|^2|\phi_{U_X}(w_1)|^2|\phi_{U_Y}(w_2)|^2\,\nu(\d w_1)\,\nu(\d w_2)\,.
\end{align*}


Hence, following similar arguments as in the proof of Theorems~\ref{thm:iff:ed} and \ref{thm:iff:dcov}, we can prove the following proposition.

\begin{proposition}\label{prop:me}
Suppose  $\tilde X$ and $\tilde Y$ are defined by 
\eqref{me}, where the characteristic functionals $\phi_{U_X},\phi_{U_Y}$ of the measurement errors do not vanish. If $\nu$ is a characteristic reference measure, then, 
\begin{align*}
  \ed_\nu^2(\mathbb P^{\tilde X},\mathbb P^{\tilde Y})=0   & ~~\text{if and only if } \ed_\nu^2(\mathbb P^{X},\mathbb P^{Y})=0\,,~\text{when }\phi_{U_X}=\phi_{U_Y};\\
  \dcov_\nu^2(\tilde X,\tilde Y)=0 
  & ~~\text{if and only if }   \dcov_\nu^2(X,Y)=0 ~. 
\end{align*}
\end{proposition}

Proposition~\ref{prop:me} states that despite the measurement errors, the energy distance and distance covariance are still valid metrics for $\H$-valued random variables, so that statistical inference on equal distribution or independence can be carried out using the contaminated observations, even in presence of measurement errors.

For the sake of brevity we illustrate our approach for the two-sample problem (a method for independence testing can be developed in a similar way). In this case, the data $ \tilde \T_{n,2} = \{
\tilde Z_1,\ldots,\tilde Z_{n_1},\tilde Z_{n_1+1},\ldots,\tilde Z_{n_1+n_2} \} $ consist of $\H$-valued random variables that satisfy
\begin{align*}
\tilde Z_{j}=Z_j+U_j\,,\quad 1\leq j\leq n_1+n_2\,,\quad\text{where } Z_1,\ldots,Z_{n_1}\iidsim\mathbb P_1\,,\  Z_{n_1+1},\ldots,Z_{n_1+n_2}\iidsim \mathbb P_2\,,
\end{align*}
where $U_1,\ldots,U_{n_1+n_2}\iidsim\mathbb P^U$, with characteristic functional $\phi_U$, are independent of $\{Z_j\}_{1\leq j\leq n_1+n_2}$. 

In view of Proposition~\ref{prop:me}, the empirical energy distance $\tilde\ed_\nu^2(\hat{\mathbb P}_1,\hat{\mathbb P}_2)$ and its permuted version $\tilde \ed_\nu^2(\hat{\mathbb P}_1,\hat{\mathbb P}_2,\Pi_b)$ are respectively defined by, for the $Q_{j,k}$ in \eqref{new},
\begin{align*}
&\tilde\ed_\nu^2(\hat{\mathbb P}_1,\hat{\mathbb P}_2):=\frac{n_1+n_2}{n_1n_2}\sum_{j,k=1}^{n_1+n_2}\phi_{\nu}(\tilde Z_{j}-\tilde Z_{k})Q_{j,k}\,,\\ 
&\tilde\ed_\nu^2(\hat{\mathbb P}_1,\hat{\mathbb P}_2,\Pi_b):=\frac{n_1+n_2}{n_1n_2}\sum_{j,k=1}^{n_1+n_2}\phi_{\nu}(\tilde Z_{\Pi_b(j)}-\tilde Z_{\Pi_b(k)})Q_{\Pi_b(j),\Pi_b(k)}\,,
\end{align*}
are reasonable test statistics for the the two-sample problem in presence of measurement error. Therefore we define the test for the hypotheses \eqref{twosample} by
\begin{align}\label{tildepsime}
\tilde\psi_{\nu,n,2}(\alpha)=\one\bigg\{\frac{1}{B}\sum_{b=1}^B\one\{\tilde\ed_\nu^2(\hat{\mathbb P}_1,\hat{\mathbb P}_2)\leq\tilde\ed_\nu^2(\hat{\mathbb P}_1,\hat{\mathbb P}_2,\Pi_b)\}>1-\alpha\bigg\}\,.
\end{align}
Since under $H_0$, $\tilde Z_1,\ldots,\tilde Z_{n_1+n_2}$ are exchangeable, $\tilde\psi_{\nu,n,2}(\alpha)$ achieves uniform type-I error control at nominal level $\alpha$, that is $\tilde\psi_{\nu,n,2}(\alpha)\in\Psi(\alpha)$.
Next, we investigate the uniform separation rate of the decision rule $\tilde\psi_{\nu,n,2}(\alpha)$. For this purpose we take $\nu$ to be a (mean-zero) Gaussian reference probability measure and define 
for  a constant $c>0$ 
\begin{align*}
\V_G(c)=\big\{\nu\in\V:\nu\text{ is Gaussian on $\H$ with non-degenerate covariance operator }C_\nu;\, \tr(C_\nu)\leq c\big\}\,
\end{align*}
as the set of mean-zero Gaussian reference probability measures with covariance operator whose trace is bounded by $c$.

For $\nu\in\V_G(c)$, consider the local alternative that the squared $L_1$-norm of the difference in the characteristic functionals $\phi_1-\phi_2$ with respect to $\nu$ is at least $\rho>0$:
\begin{align}\label{f1nu}
\mathcal F_{1,\nu}(\rho)=\Big \{\mathbb P_1,\mathbb P_2\in\mathcal P(\H):\Big [\int_\H|\phi_{1}(w)-\phi_{2}(w)|\,\nu(\d w)\Big ]^2\geq\rho\Big \}\,.
\end{align}
In addition, for some constants $c_U>0$ and $s>0$, define   the function class
\begin{align}\label{Fu}
\mathcal F_U(c_U,s)=\Big\{\mathbb P^U\in\mathcal P(\H):|\phi_U(w)|\geq c_U(1+\|w\|)^{-s}\,,\text{for }\,w\in\H\Big\}\,
\end{align}
for the distribution   $\mathbb P^U$ of the measurement error. In other words,
 we assume that  $|\phi_U(w)|$ does not decay faster than a polynomial rate with respect to~$\|w\|$. Conditions of this type are  standard for measurement error problems in Euclidean spaces; see, for example, \cite{carroll1988,meister2009}.
For the mean-zero Laplace measure $\nu_L\in\mathcal F_U(c_U,s)$ on $\H$ with covariance operator $C_L$, by Theorem~\ref{thm:laplace}, its characteristic functional $\phi_L$ is such that $|\phi_L(w)|=(1+2^{-1}\l Cw,w\r)^{-1}\geq c_L(1+\|w\|)^{-1}$, where $c_L=\max\{2,\lambda_1(C_U)\}/2>0$, and $\lambda_1(C_L)$ denotes the leading eigenvalue of $C_L$.

The following proposition shows the upper bound for the uniform separation rate of the two sample test $\tilde\psi_{\nu,n,2}(\alpha)$ in presence of measurement error, and is proved in Section~\ref{app:prop:me:two} of the supplementary material.
\begin{proposition}\label{prop:me:two}
Assume $0<\alpha+\beta<1$, $c_1\leq n_1/ n_2\leq c_1$ for some constants $c_1,c_2>0$, and $\beta<1/e$, $B\geq2\alpha^{-2}\log(6/\beta)$. 
The uniform separation rate of  the two-sample test $\tilde\psi_{\nu,n,2}(\alpha)$ defined in \eqref{tildepsime}  with respect to~the set $\mathcal F_{1,\nu}(\rho)$ in \eqref{f1nu} satisfies
\begin{align*}
\inf\bigg\{\rho>0:\sup_{\nu\in\V_G(c)}\sup_{\substack{(\mathbb P_1,\mathbb P_2)\in\mathcal F_{1,\nu}(\rho)\\ \mathbb P^U\in\mathcal F_U(c_U,s)}}\P\{\tilde\psi_{\nu,n,2}(\alpha)=0\}\leq\beta\bigg\}\leq c_0(n_1+n_2)^{-1/2}\{\log(\alpha^{-1})+{\log(\beta^{-1})}\}\,,
\end{align*}
where the constant $c_0>0$ only depends on $c,s$ and $c_U$.
\end{proposition}

Note that a difference from the upper bound in Theorem~\ref{thm:upper:two} in the error-free case is that the separation from the null distributions is quantified with respect to~the $L_1$-norm in \eqref{f1nu}.
The challenge of deriving a lower bound in presence of measurement errors lies in computing the Radon-Nikodym derivatives of non-Gaussian probability measures on $\H$. Investigations of this type are  beyond the scope of the present  article and left  for future work.

\subsection{Aggregated permutation tests}\label{sec:agg}

Although our approach works with a wide range of reference probability measures, in practice, for a given sample, a theoretical optimal choice is unclear. To circumvent this difficulty, we follow the approach proposed in \cite{fromont2013} recently adopted by \cite{albert2022,schrab2022} and consider aggregated (permutation) tests. It uses, rather than a single, but several reference probability measures collected in a set $\V_0\subset\V$, to compute and aggregate the corresponding statistics.

To be precise, consider the problem of testing for  independence, that is $H_0: \mathbb{P}^{X,Y} = \mathbb{P}^{X} \otimes \mathbb{P}^{Y}$. 
For a given  sample $\T_n=\{(X_1,Y_1),\ldots,(X_n,Y_n)\}$ the  level-$\alpha$ aggregated permutation independence test consists of the following two steps. First, for each $\nu\in\V_0$ and $u\in (0,1)$, the $u$-quantile $q_\nu(u)$ of $\dcov_\nu^2(X,Y)$ is estimated via the permutation procedure, by taking the $\lceil (1-u)B_1\rceil$-th order statistic of  the sample $\{\hat\dcov_\nu^2(X,Y,\Pi_b)\}_{1\leq b\leq B_1}$ defined in \eqref{hatdb} computed from $B_1$ permuted samples $\T_n(\Pi_1),\ldots,\T_n(\Pi_{B_1})$, that is 
\begin{align}\label{hatq}
\hat q_{\nu}(u)=\big\{\hat\dcov_\nu^2(X,Y,\Pi_b):1\leq b\leq B_1\big\}_{(\lceil (1-u)B_1\rceil)}\,.
\end{align}

The second step involves estimating the adjusted nominal level.
For $B_2\in\mathbb N_+$, we generate $B_2$ independent permutations $\tilde\Pi_1,\ldots,\tilde\Pi_{B_2}$ that are independent of $\Pi_1,\ldots,\Pi_{B_1}$. For the permuted statistics $\{\hat\dcov_\nu^2(X,Y,\tilde\Pi_b)\}_{1\leq b\leq B_2}$ 
computed from the permuted samples $\T_n(\tilde\Pi_1),\ldots,\T_n(\tilde\Pi_{B_2})$, we estimate the adjusted nominal level by
\begin{align}\label{sup2}
\hat u_\alpha=\sup\Big \{u>0:\frac{1}{B_2}\sum_{b=1}^{B_2}\one\Big[\max_{\nu\in\V_0}\big\{\hat\dcov_\nu^2(X,Y,\tilde\Pi_b)-\hat q_{\nu}(1-u)\big\}>0\Big]>1-\alpha\Big \}\,,
\end{align}
where the empirical quantile $\hat q_{\nu}(u)$ is defined by \eqref{hatq}.
Finally, for the  sample $\{\hat\dcov_\nu^2(X,Y)\}_{\nu\in\V_0}$ computed  $\T_n$ using the reference probability measures in the set $\V_0$, the decision rule for the aggregated permutation test is defined by
\begin{align}
\label{hatpsi:ind}
\hat\psi_{n}(\alpha)=\one\Big[\max_{\nu\in\V_0}\{\hat\dcov_\nu^2(X,Y)-\hat q_{\nu}(1-\hat u_\alpha)\}>0\Big]\,.
\end{align}

We conclude this section with the following proposition that establishes the theoretical property of the aggregated permuted test, which is proved in Section~\ref{app:prop:agg} of the supplementary material.

\begin{proposition}\label{prop:agg} ~~\\
(i) For the type-I error of the aggregated test, we have $\hat\psi_n\in \Psi(\alpha)$. \\
(ii) If $\beta<1/e$, $B_1\geq 8\alpha^{-2}\log(6/\beta)|\V_0|^2$, $B_2\geq2\alpha^{-2}\log(2/\beta)$, the uniform separation rate satisfies
\begin{align*}
\rho\{\hat\psi_n(\alpha),\beta\}\leq c n^{-1/2}\{\log(\alpha^{-1})+{\log(\beta^{-1})}+\log(|\V_0|)\},
\end{align*}
where $c>0$ is an absolute constant.
\end{proposition}

In practice, the supremum in \eqref{sup2} can be estimated using the dichotomy method, and for a fixed nominal level $\alpha$, we may take $B_1=B_2=B$ for brevity. The following algorithm states the procedures for aggregated permutation test for inference on independence of $\H$-valued random variables.

\begin{algorithm}[Aggregated permutation independence test for $\H$-valued random variables]~\label{algo:indep}
\smallskip

{\bf Input:} The original data $\T_n=\{(X_j,Y_j)\}_{j=1}^n$, nominal level $\alpha$, a set of $M$ reference probability measures $\V_0=\{\nu_1,\ldots,\nu_M\}$  on $\H$, number of permutations $B$.
\begin{enumerate}[label={(\arabic*)},nolistsep]
\item Generate two sets of independent  random permutations of $(1,2,\ldots,n )$: $\mathfrak G_1=\{\Pi_1,\ldots,\Pi_{B}\}$ and $\mathfrak G_2=\{\tilde\Pi_1,\ldots,\tilde\Pi_{B}\}$, all independent of the original data $\T_n$.

\item Compute $\hat V^{(X,\nu)}_{k,\ell}$ and $\hat V^{(Y,\nu)}_{k,\ell}$ defined in \eqref{hatv}.

\item\label{(3)} Compute the empirical distance covariance test statistics for all $\nu\in\V_0$ and $1\leq b\leq B$:
\begin{align*}
&\hat\dcov_\nu^2(X,Y)=\frac{1}{n^2}\sum_{k,\ell=1}^n\hat V_{k,\ell}^{(X,\nu)}\hat V_{k,\ell}^{(Y,\nu)}\,,\qquad\\
&\hat\dcov_\nu^2(X,Y,\Pi_b)=\frac{1}{n^2}\sum_{k,\ell=1}^n\hat V_{k,\ell}^{(X,\nu)}\hat V_{\Pi_b(k),\Pi_b(\ell)}^{(Y,\nu)}\,,\qquad\hat\dcov_\nu^2(X,Y,\tilde\Pi_b)=\frac{1}{n^2}\sum_{k,\ell=1}^n\hat V_{k,\ell}^{(X,\nu)}\hat V_{\tilde\Pi_b(k),\tilde\Pi_b(\ell)}^{(Y,\nu)}\,.
\end{align*}

\item Compute the adjusted nominal level $\hat u_\alpha$ in \eqref{sup2}. For $\nu\in\V_0$, compute $\hat q_{\nu}(1-\hat u_\alpha)$ in \eqref{hatq}, by taking the $(\lceil (1-\hat u_\alpha)B\rceil)$ order statistic of $\{\hat\dcov_\nu^2(X,Y,\tilde\Pi_b)\}_{1\leq b\leq B}$.

\item Output the decision rule $\hat\psi_n(\alpha)$ in \eqref{hatpsi:ind}: reject $H_0$, if $\max_{\nu\in\V_0}\{\hat\dcov_\nu^2(X,Y)-\hat q_{\nu}(1-\hat u_\alpha)\}>0$.

\end{enumerate}

\end{algorithm}

\begin{remark}[Aggregated permutation two-sample test for $\H$-valued random variables]\label{rem:agg:two}
For the two-sample test, the algorithm is parallel to Algorithm~\ref{algo:indep}, where in step \ref{(3)} we compute the empirical energy distance $\ed_\nu^2(\hat{\mathbb P}^X,\hat{\mathbb P}^Y)$ and its permuted version $\ed_\nu^2(\hat{\mathbb P}^X,\hat{\mathbb P}^Y,\Pi_b)$ as follows: for $1\leq b\leq B$,
\begin{align*}
&\ed_\nu^2(\hat{\mathbb P}^X,\hat{\mathbb P}^Y)=\frac{n_1+n_2}{n_1n_2}\sum_{k,\ell=1}^{n_1+n_2}\phi_{\nu}(Z_k-Z_\ell)Q_{j,k};\\ &\ed_\nu^2(\hat{\mathbb P}^X,\hat{\mathbb P}^Y,\Pi_b)=\frac{n_1+n_2}{n_1n_2}\sum_{k,\ell=1}^{n_1+n_2}\phi_{\nu}(Z_{\Pi_b(k)}-Z_{\Pi_b(\ell)})Q_{\Pi_b(j),\Pi_b(k)};\\
&\ed_\nu^2(\hat{\mathbb P}^X,\hat{\mathbb P}^Y,\tilde\Pi_b)=\frac{n_1+n_2}{n_1n_2}\sum_{j,k=1}^{n_1+n_2}\phi_{\nu}(Z_{\tilde\Pi_b(k)}-Z_{\tilde\Pi_b(\ell)})Q_{\tilde\Pi_b(j),\tilde\Pi_b(k)},
\end{align*}
where, for $1\leq k,\ell\leq n$, the $Q_{j,k}$'s are defined in \eqref{new}. 
\end{remark}

\section{Finite sample properties}\label{sec:simu}

We illustrate 
the proposed inference methods for $L^2([0,1])$-random variables by means of a simulation study. The sample curves are evaluated on a equally spaced grid on the interval $[0,1]$ of $N=40$ and  $N=200$ time points, we use $B=300$ permutations  for the test and  the  nominal level is $\alpha=0.05$. For the aggregation tests, we use   mean-zero Gaussian and Laplace measures as  reference probability measures 
with various covariance kernels
\begin{equation} \label{det52}
\begin{split}
    {\rm (i)} & ~~C_1(s,t)=\min\{s,t\}  ~~~~~~~
     {\rm (ii)}  ~~ C_2(s,t)=\min\{s,t\}-st  ~~
      {\rm (iii)}  ~~ C_3(s,t)=\exp(-|s-t|^2) \\
       {\rm (iv)} & ~~  C_4(s,t)=\exp(-|s-t|) ~~  {\rm (v)}  ~~  C_5(s,t)=(1+|s-t|^2)^{-1}
\end{split}
\end{equation}
All results presented in this section are based on $500$ simulation runs.


\subsection{Independence tests} 
\begin{figure}
\centering
\includegraphics[width=4cm]{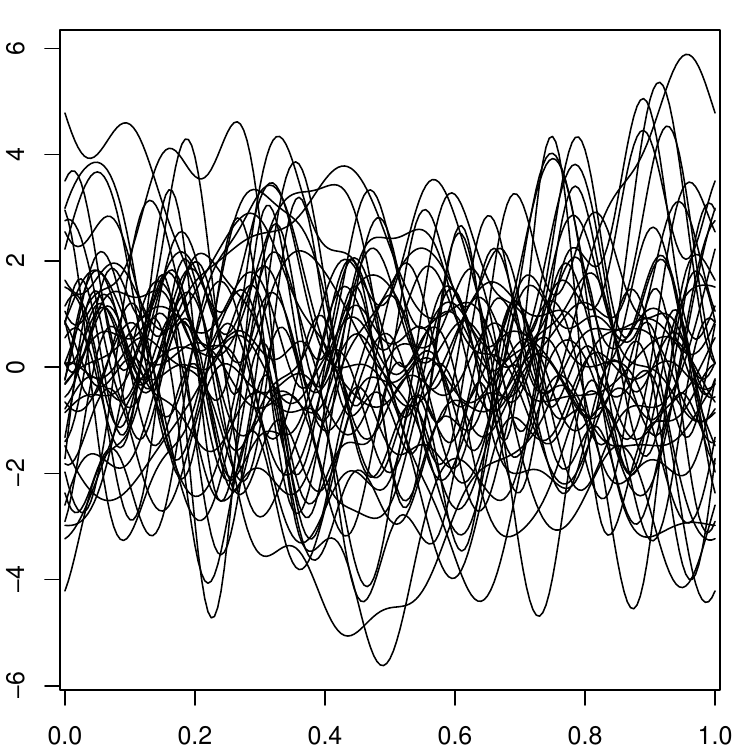}\includegraphics[width=4cm]{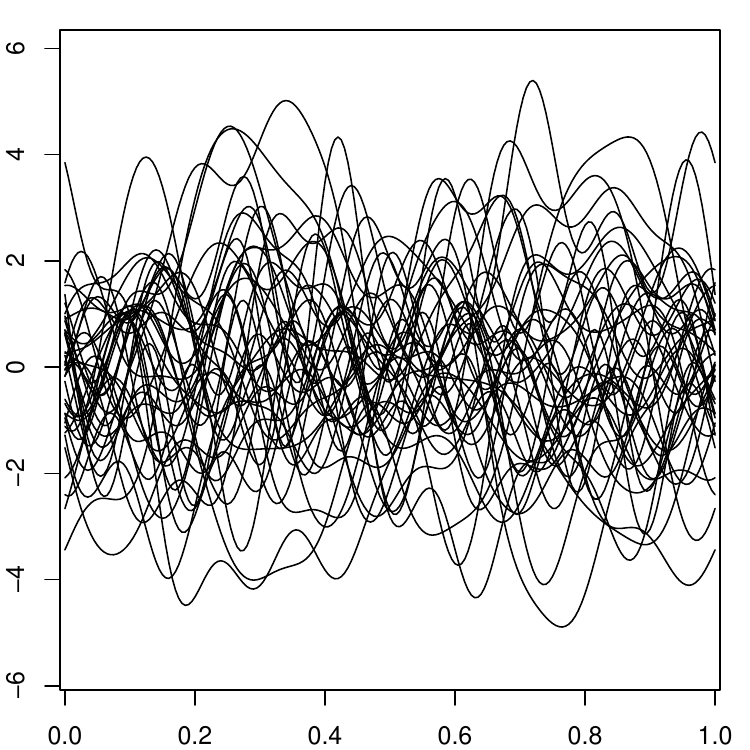}
\includegraphics[width=4cm]{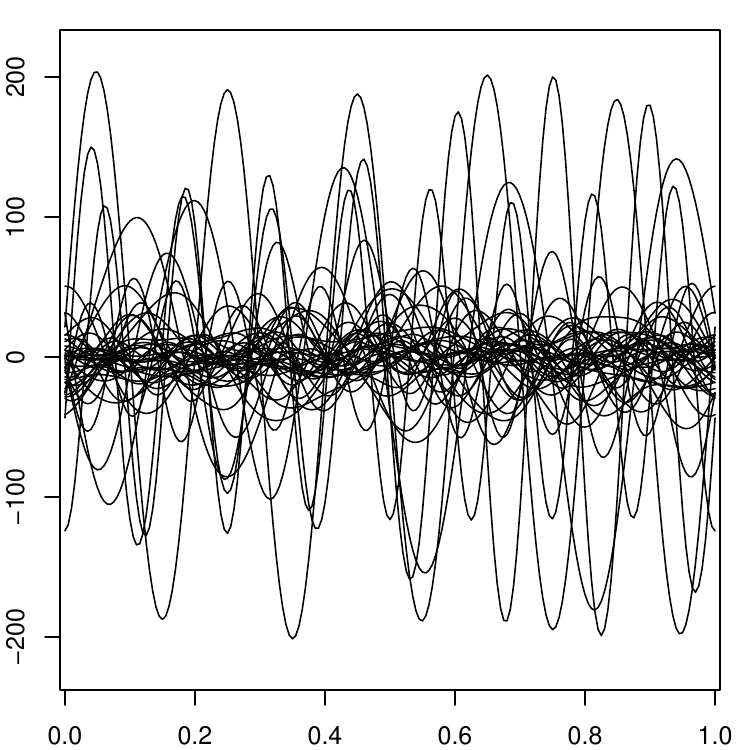}\includegraphics[width=4cm]{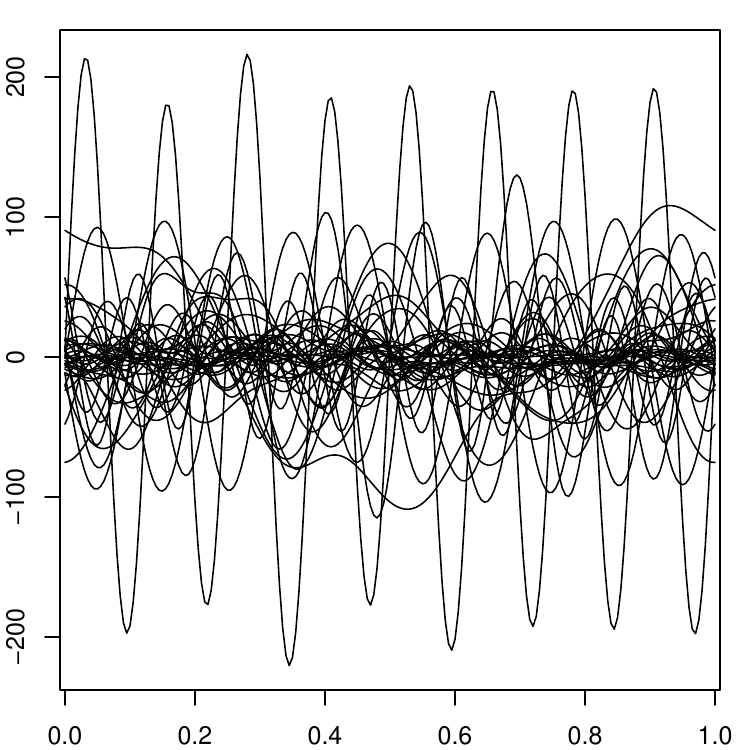}
\caption{\em Sample curves of $X$ and $Y$ in Settings (1) (first and second panel) and (2) (third and fourth panel), respectively. The sample curves in Setting (2) are heavy-tailed.\label{fig:0}}
\end{figure}

\begin{figure}
\centering
\includegraphics[width=4cm]{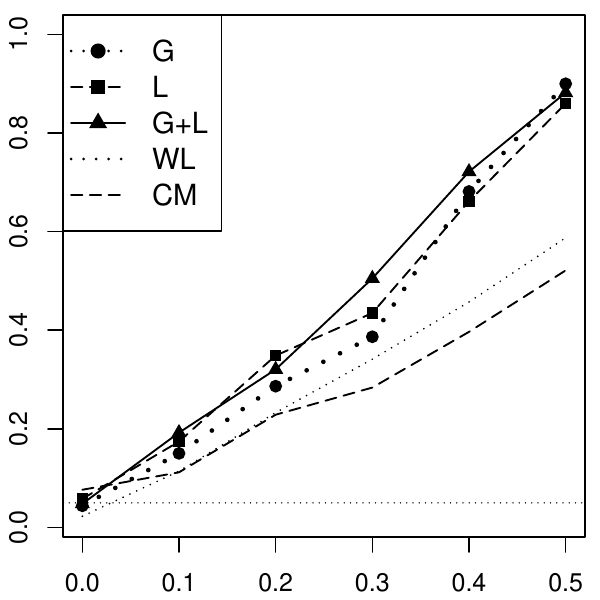}
\includegraphics[width=4cm]{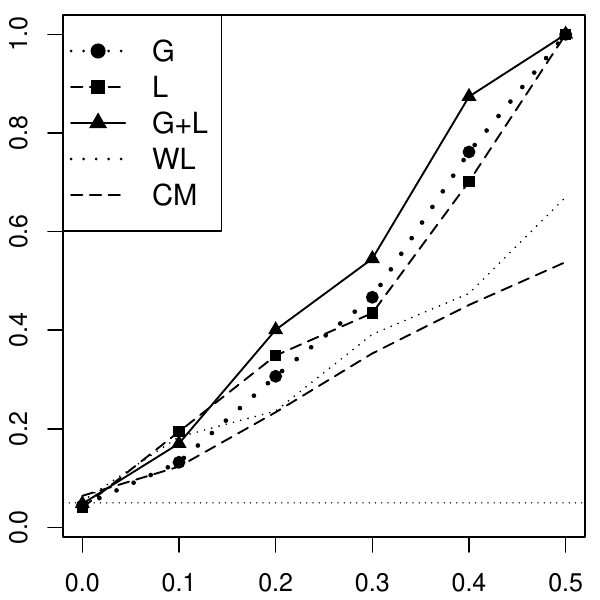}
\includegraphics[width=4cm]{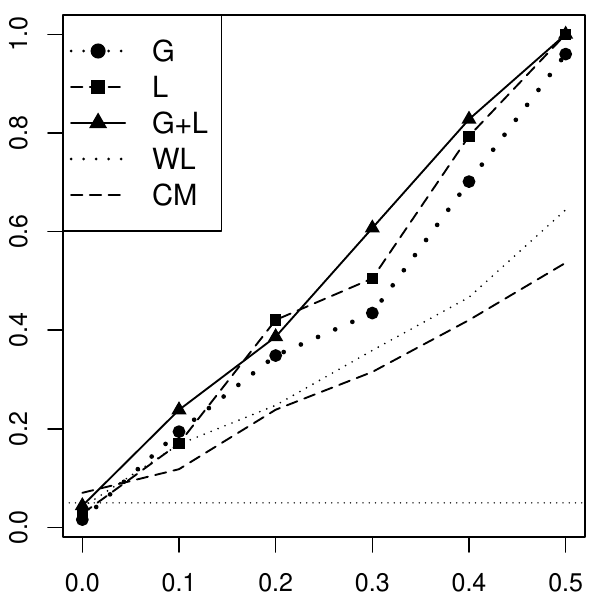}
\includegraphics[width=4cm]{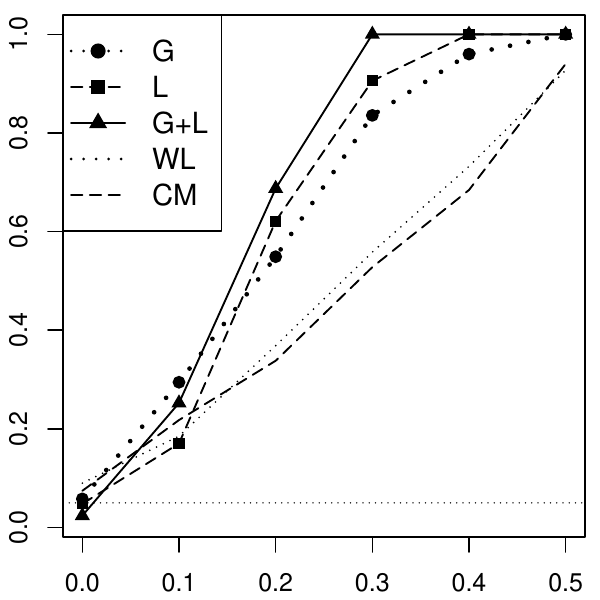}
\includegraphics[width=4cm]{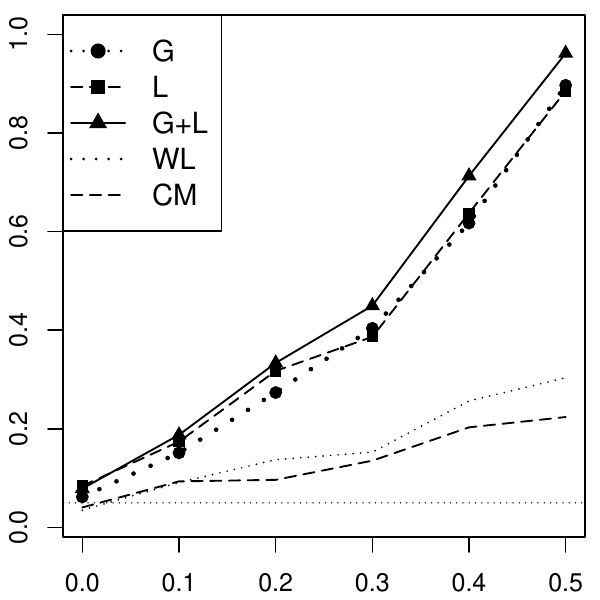}
\includegraphics[width=4cm]{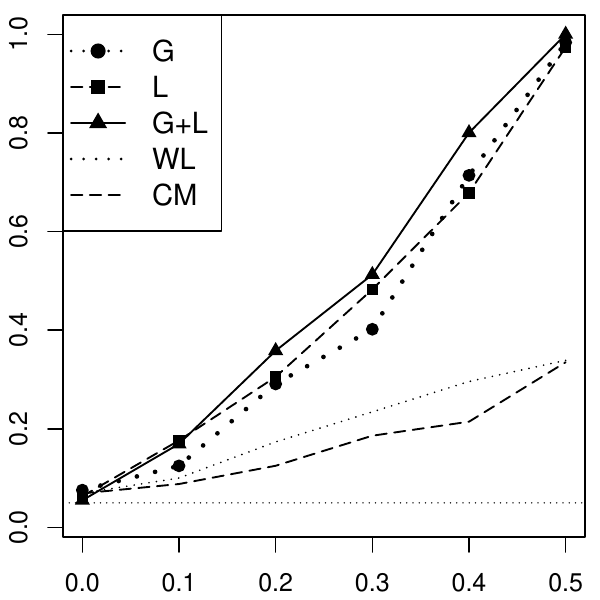}
\includegraphics[width=4cm]{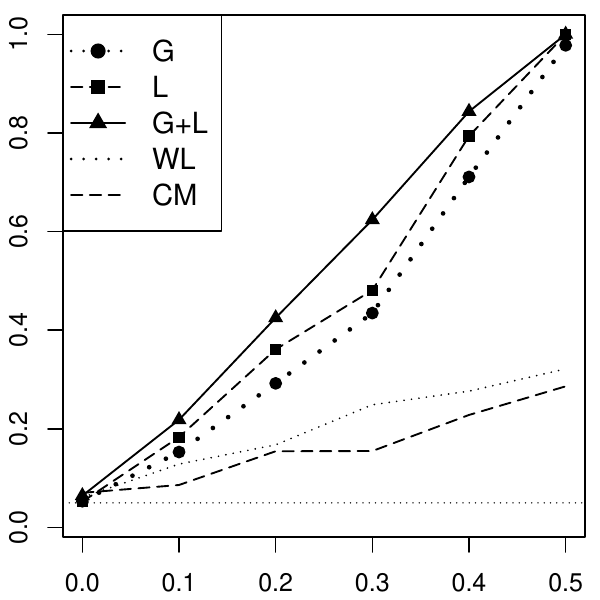}
\includegraphics[width=4cm]{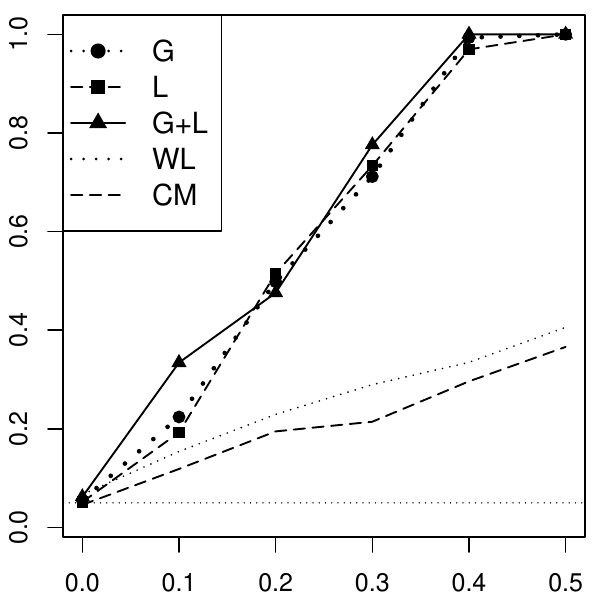}
\vspace{-1em}\caption{\em Empirical rejection probabilities (y-axis) under Settings (1) (first row) and (2) (second row), for various values of $\theta$ (x-axis); left to right: $(n,N)=(40,40), (40,200),(200,40),(200,200)$. 
Horizontal dotted line: nominal level $\alpha=0.05$.\label{fig:1}}
\end{figure}

We consider the  samples size $n=40$ and  $n=200$ of i.i.d.~bivariate $L^2([0,1])$-random variables $\{(X_j,Y_j)\}_{1\leq j\leq n}$ from the following two settings:

\begin{enumerate}[nolistsep,label={(\arabic*)}]
\item Let $\{ \eta_k \} _{k \geq }$ denote the common Fourier basis, that is $\eta_0(t)=1$
$\eta_{2 \ell-1}(t)=\sqrt{2} \cos (2 \pi \ell t)\,,\eta_{2 \ell}(t)=\sqrt{2} \sin (2 \pi \ell t)\,  $
$(\ell \geq 1)$
and define for $t\in[0,1]$
\begin{equation}
  \label{det53}  
      X_j(t)  =\sum_{k=1}^{16}U_{jk}\eta_k(t) ~,~~
      Y_j(t)=\sum_{k=1}^{16} V_{jk} \eta_k(t+0.2),
\end{equation}
where, $(U_{jk},V_{jk})^\top$ are i.i.d.~jointly mean-zero Gaussian with covariance matrix 
$$
\bigg 
(\begin{matrix}k^{-1.05} & \rho_k k^{-1.125} \\
\rho_k k^{-1.125} & k^{-1.2}
\end{matrix}
\bigg)
,
$$
and $\rho_k=0$ for $1\leq k\leq 8$, and $\rho_k=\theta=0,0.1,\ldots,0.5$ for $9\leq k\leq 16$.

\item  Let $U_{jk}$ and $V_{jk}$ be independent 
Cauchy distributed random variables 
and define for $t\in[0,1]$
\begin{equation}
  \label{det54} 
X_{j}(t) =\sum_{k=1}^{16}U_{jk}\eta_k(t) ~,~~
Y_{j}(t) =\sum_{k=1}^{16} \big\{\theta_{} U_{jk}\eta_k(t+0.2)+(1-\theta)V_{jk}\eta_k(t)\big\}
,
\end{equation}
where  $\theta=0, 0.1,\ldots,0.5$.

\end{enumerate}
Setting (1) is similar to the second setting in \cite{miao2022} with different covariance structures, and Setting (2) is a non-Gaussian case of heavy tailed distributions. Taking $\theta=0$ in Settings~(1) and (2) corresponds to  the null hypothesis of independence. The sample curves in Settings~(1) and (2), with $\theta=0.5,n=40,N=200$ are displayed in Figure~\ref{fig:0}.
For the different settings we compare the aggregated independence test 
$\hat\psi_{n}$ in \eqref{hatpsi:ind} with the wavelet method in \cite{miao2022} (denoted by WL) and the conditional mean approach in \cite{lee2020} (denoted by CM), which requires the existence of the conditional mean.  In Figure~\ref{fig:1} we show the empirical rejection probabilities of the different methods, where for the aggregated test $\hat\psi_{n}$, the  reference probability measures are chosen as Gaussian (G), Laplace (L), and both (G+L), respectively, with 
the different covariance kernels in \eqref{det54}.
The horizontal dotted line displays  the nominal level $\alpha=0.05$. The results show  a reasonable approximation of the nominal level under the null hypothesis for all methods under consideration. 
Moreover, we observe  advantages of the aggregation tests over the WL and CM approaches  in terms of power for  both scenarios.  While  the power of the new  test  \eqref{hatpsi:ind} is only slightly decreasing  in the presence of heavy-tailed data, we observe a stronger decrease in power for the WL and CM test in this case. In all cases the aggregation test with the Gaussian and Laplace reference measures yields the best results.

\begin{figure}
\centering
\includegraphics[width=4cm]{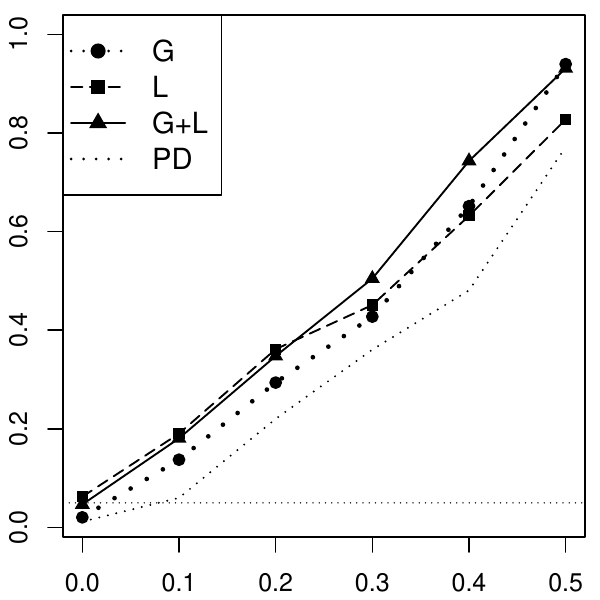}
\includegraphics[width=4cm]{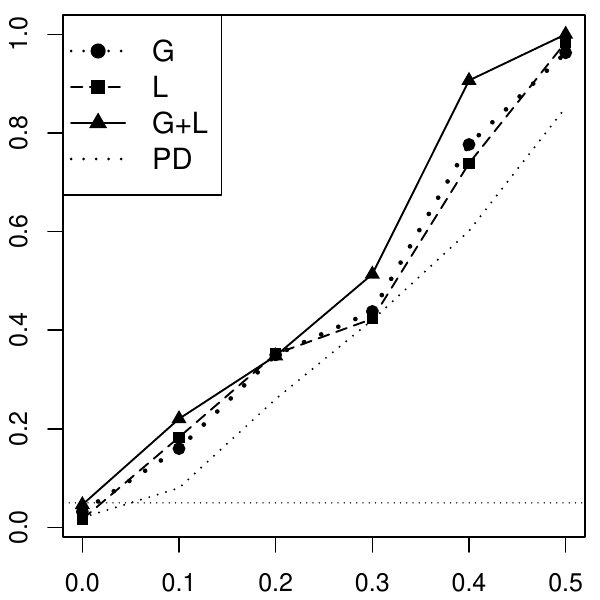}
\includegraphics[width=4cm]{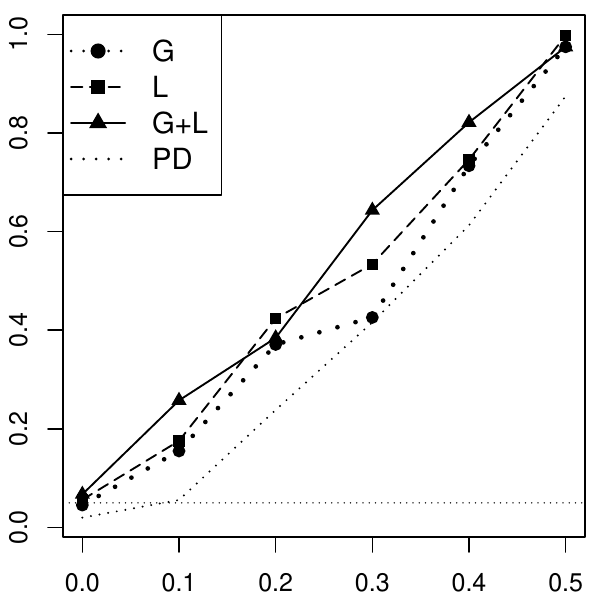}
\includegraphics[width=4cm]{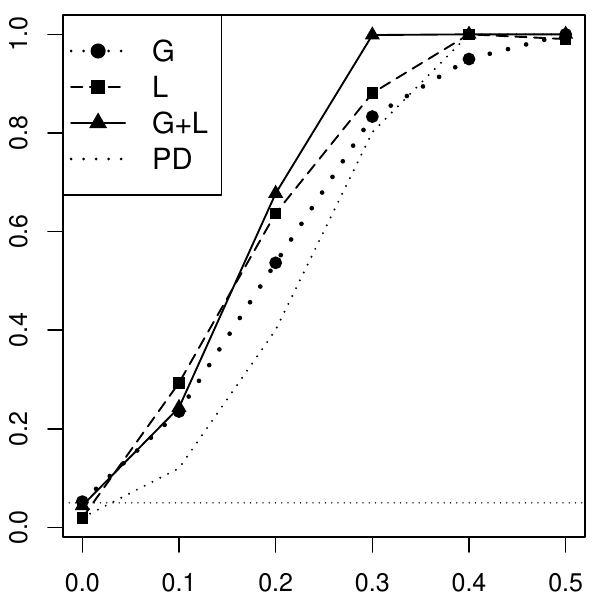}
\includegraphics[width=4cm]{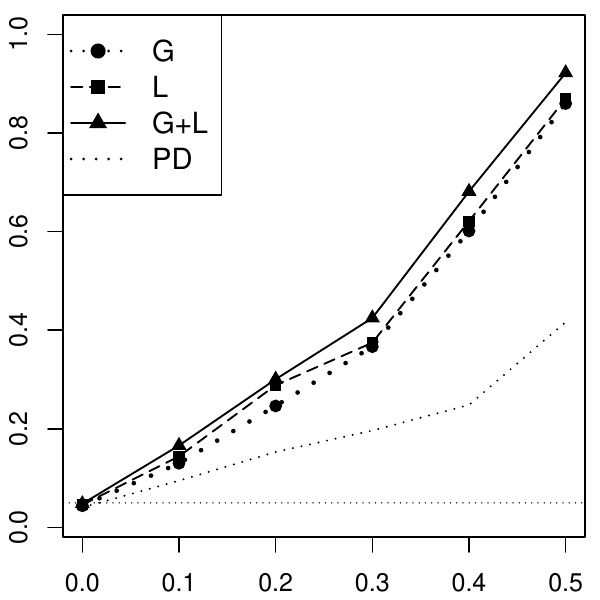}
\includegraphics[width=4cm]{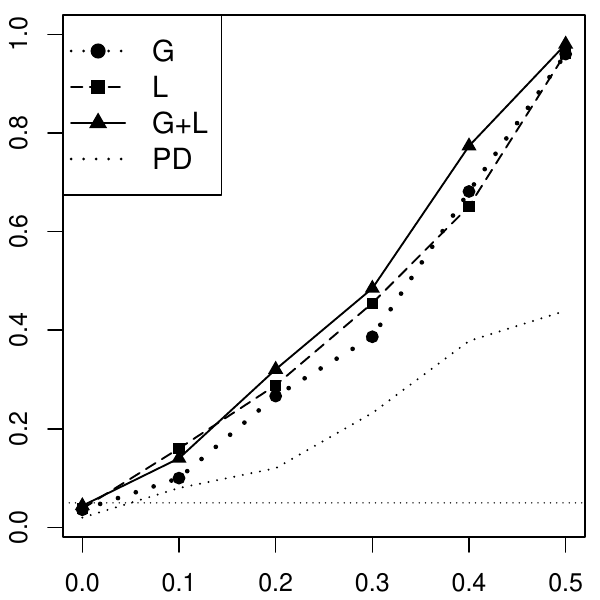}
\includegraphics[width=4cm]{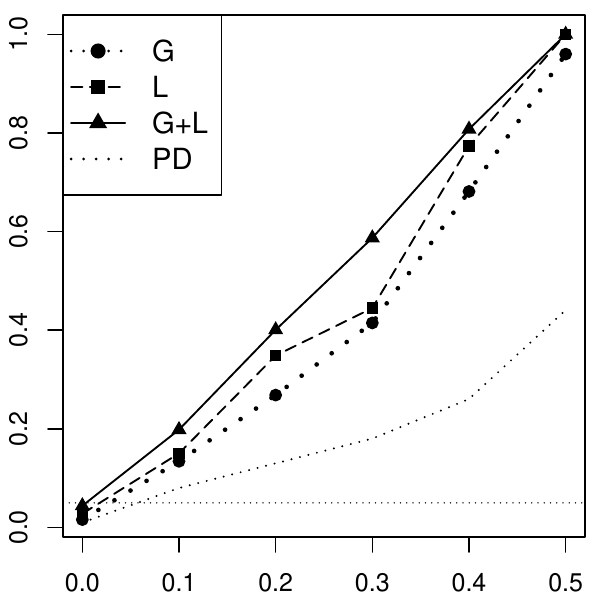}
\includegraphics[width=4cm]{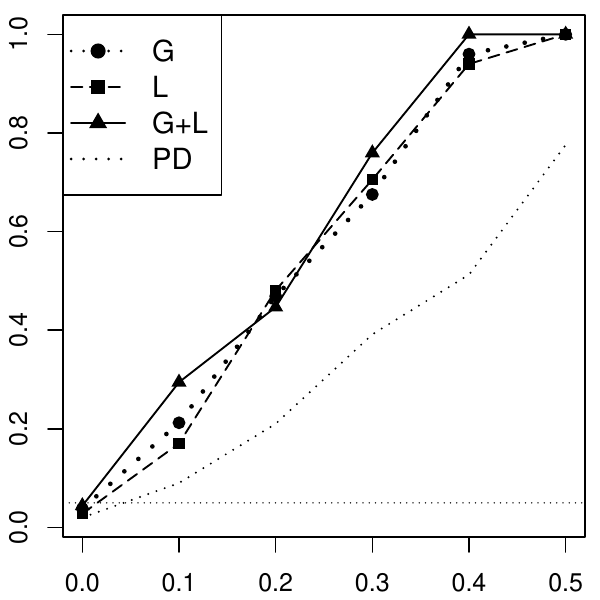}
\vspace{-1em}\caption{\em Empirical rejection probabilities (y-axis) under Settings (3) (first row) and (4) (second row), for two-sample inference, for various values of $\theta$ (x-axis); left to right: $(n,N)=(40,40), (40,200),(200,40),(200,200)$. Horizontal dotted line: nominal level $\alpha=0.05$. \label{fig:3}}
\end{figure}

\subsection{Two-sample tests}
For the two-sample test for the hypotheses \eqref{twosample} we consider independent  samples $X_1,\ldots,X_{n}$ and $Y_1,\ldots,Y_{n}$  of   i.i.d.~random variables, each of size $n=40$ or $n=200$, from the following two settings:

\begin{enumerate}[nolistsep]
\item[(3)] $X_j(t)=\sum_{k=1}^{16}U_{jk}\eta_k(t)$ and $Y_j(t)=\sum_{k=1}^{16} V_{j k} \eta_k(t)$,  for $t\in[0,1]$, where $\{\eta_k\}_{k=1}^{16}$ are Fourier basis defined in Setting (1), $U_{jk}$ and $V_{jk}$ are independent mean-zero Gaussian random variables with variances $k^{-1.05}$ and $\rho_kk^{-1.2}+k^{-1.05}$, respectively, where  $\rho_k=0$ for $1\leq k\leq 8$ and $\rho_k=\theta=0,0.1,\ldots,0.5$, for $9\leq k\leq 16$.

\item[(4)] same as Setting (2).
\end{enumerate}
Taking $\theta=0$ in Settings (3) and (4) corresponds to  null hypothesis of equal distributions in the two samples. Under different settings we compare the aggregation two-sample test (Remark~\ref{rem:agg:two}) with the homogeneity test based on the point-wise distribution function proposed in \cite{bugni2021permutation} who  propose a two-sample test based on point-wise distribution functions
(denoted by PD). In Figure~\ref{fig:3}  we display the empirical rejection probabilities of PD and of  the aggregated two-sample test described in Remark~\ref{rem:agg:two}, where the  reference probability measures are chosen as Gaussian (G), Laplace (L), and both (G+L), respectively, with the different covariance kernels in \eqref{det52}. From the results we observe that the nominal level $\alpha$ is well-approximated under the null hypothesis for both methods. However,  the new  aggregation tests entail higher empirical power and the power is only slightly decreasing under heavy-tailed distributions in Setting~(4). In all cases  the aggregation tests with Gaussian and Laplace reference measures yield the best results.

\medskip 
{\bf Acknowledgements}
This work was partially supported by the DFG Research unit 5381 \textit{Mathematical Statistics in the Information Age}, project number 460867398.  The authors would like to thank
Anatoly Zhigljavsky  for some very  helpful discussions on energy distances.


\baselineskip=11pt

\bibliographystyle{apalike}

\clearpage

\appendix

\setcounter{page}{1}

\appendix
\renewcommand{\thesection}{\Alph{section}}

\rhead{\textbf{\thepage}}
\lhead{\textbf{Online supplement}}

\section*{Supplementary Material}

\vspace{.5cm}

In this supplementary material we provide technical details of our theoretical results. 
In Section~\ref{app:proof} we provide the proofs of our main results in our main article. In Section~\ref{app:aux} we provide supporting lemmas that are used in the proofs in Section~\ref{app:proof}. In the sequel, we use ``$c$" to denote a generic positive constant (independent of $n$ and the reference probability measure $\nu$) which might differ from line to line. Let $I_n(n)$ denote  the set consisting of all permutations
of $ (1, \ldots , n)$.

\section{Proof of main results}\label{app:proof}
\subsection{Useful lemmas and definitions}

We first recall some useful definitions and propositions.
%
Recall the definitions of mean and covariance operators of a probability measure on $\H$; see, for example Definitions~2.1 and 2.3 in \cite{kuo1975}.
\begin{definition}[Mean and covariance operator]\label{def:meanvar}
Let $\nu$ be a Borel measure on $\H $. The mean $\mu_\nu$ and covariance operator $C_\nu$ of $\nu$ are defined by
\begin{align*}
\l \mu_\nu,x\r=\int_\H  \l z,x\r\,\nu(\d z)\,,\qquad \l C_\nu x,y\r=\int_\H \l x,z-\mu_\nu\r\l y,z-\mu_\nu\r\,\nu(\d z)\,,\qquad  x,y\in \H \,.
\end{align*}
\end{definition}

The following Minlos–Sazonov theorem 
states the necessary and sufficient conditions for a functional on $\H$ to be a well-defined characteristic functional of a probability measure on $\H$, which is useful for proving Theorems~\ref{thm:iff:ed}, \ref{thm:laplace} and \ref{thm:iff:dcov}; see for example Theorem~2.2 of \cite{kuo1975} or Theorem~4.5 in \cite{kukush2020} for a proof.
\begin{theorem}[The Minlos–Sazonov theorem]\label{thm:ms}
A functional $\phi:\H\to\mathbb C$ is the characteristic functional of a probability measure on $\H $ if and only if the following two statements are true.

\begin{enumerate}[nolistsep]
\item[{\rm (i)}] $\phi(0)=1$ and $\phi$ is a positive definite functional;

\item[{\rm (ii)}] For any $\epsilon>0$, there exists $S(\epsilon)\in \mathcal S$ (depending on $\epsilon$) such that $1-\Re\{\phi(x)\}\leq\l S(\epsilon) x,x\r+\epsilon$, for all $x\in \H$.
\end{enumerate}

\end{theorem}

The following proposition states the connection of a kernel and semimetric of negative-type; for its proof see for example Lemma~3.2.1 and Proposition~3.3.2 of \cite{berg1984}.

\begin{proposition}[Kernel and semimetric of negative-type]~\label{prop:kernel}
\begin{enumerate}[nolistsep,label={{\rm(\arabic*)}}]
\item Let $z\in\H$ be arbitrary. The kernel centered at $z$ defined in \eqref{knu2}
is a positive definite kernel if and only if $d_\nu$ satisfies \eqref{nt}.

\item Let $K_\nu$ denote any nondegenerate kernel on $\H$, that is, the Aronszajn map $z\mapsto K_\nu(\cdot,z)$ is injective. Then,
\begin{align*}
d_\nu(x,y)=\|K_\nu(x,\cdot)-K_\nu(y,\cdot)\|_{\H_{K_\nu}}^2=K_\nu(x,x)+K_\nu(y,y)-2K_\nu(x,y)
\end{align*}
defines a valid semimetric $d_\nu$ of negative type on $\H$.
\end{enumerate}

\end{proposition}


\subsection{Proof of Proposition~\ref{prop:ed}}\label{app:prop:ed}

In order to prove \eqref{ed:dnu}, observe that
\begin{align*}
&|\E\exp(\i\l w,X\r)-\E\exp(\i\l w,Y\r)|^2\\
&=\{\E\exp(\i\l w,X\r)-\E\exp(\i\l w,Y\r)\}\{\E\exp(-\i\l w,X\r)-\E\exp(-\i\l w,Y\r)\}\\
&=\E\exp(\i\l w,X\r)\E\exp(-\i\l w,X\r)-\E\exp(-\i\l w,X\r)\E\exp(\i\l w,Y\r)\\
&\qquad-\E\exp(\i\l w,X\r)\E\exp(-\i\l w,Y\r)+\E\exp(\i\l w,Y\r)\E\exp(-\i\l w,Y\r)\\
&=\E\exp(\i\l w,X-X'\r)-\E\exp(-\i\l w,X-Y'\r)-\E\exp(\i\l w,X-Y'\r)+\E\exp(\i\l w,Y-Y'\r)\,.
\end{align*}
Therefore, in view of \eqref{phiz}, by Fubini's theorem we obtain
\begin{align*}
\ed_\nu^2(X,Y)=\E\phi_{\nu}(X-X')-2\E\phi_{\nu}(X-Y')+\E\phi_{\nu}(Y-Y')\,.
\end{align*}

\subsection{Proof of Theorem~\ref{thm:iff:ed}}\label{app:thm:iff}

For sufficiency, suppose $\supp(\nu)=\H$. It is obvious that $P_1=P_2$ implies that $\ed_\nu^2(P_1,P_2)=0$. To prove the converse, note that $\phi_1-\phi_2$ is an uniformly continuous functional on $\H$; see for example p.~19 of \cite{kuo1975} for a proof. Then, if there exists $x_0\in\H$ such that $|\phi_1(x_0)-\phi_2(x_0)|^2>0$, there exists $r\in\mathbb R$ such that $|\phi_1(x_0)-\phi_2(x_0)|^2>0$ for any $x\in\mathcal B(x_0,r)$. Since $\supp(\nu)=\H$, it holds that $\nu\{\mathcal B(x_0,r)\}>0$. This implies that $\ed_\nu^2(P_1,P_2)\geq\int_{\mathcal B(x_0,r)}|\phi_1(x)-\phi_2(x)|^2\nu(\d x)>0$, which contradicts the fact that $\ed_\nu^2(P_1,P_2)=0$. The proves the sufficiency.

To prove necessity, suppose $\H\setminus\supp(\nu)\neq\varnothing$. Then, there exist an open non-empty set $S\subset\H$ such that $\nu(S)=0$.
Furthermore, there exists a ball $\mathcal B_r(w)=\{x\in\H:\|x-w\|^2\leq  r^2\}\subset S$ centered at $w$ with radius $r$ such that $0\notin \mathcal B_r(w)$.
Recalling that $\nu$ is symmetric, we deduce that $\nu[\mathcal B_r(w)\cup\{-\mathcal B_r(w)\}]=0$.
We start by defining a useful operator on $\H$.
Suppose $\{e_k\}_{k=1}^\infty$ is a complete orthornormal basis of $\H$.
Let $(c)_+=\max\{c,0\}$ for $c\in\mathbb R$ and define, for $x\in\H$,
let
\begin{align}\label{theta}
\theta(x)=\prod_{k=1}^\infty\Bigg[\bigg\{1-\frac{1}{3}\bigg(1-\frac{\pi k}{\sqrt{6}r}|\l x-w,e_k\r|\bigg)_+-\frac{1}{3}\bigg(1-\frac{\pi k}{\sqrt{6}r}|\l x+w,e_k\r|\bigg)_+\bigg\}\Bigg]^{-k^{-2}}-1\,.
\end{align}
Observe that $\theta(x)=\theta(-x)$ for any $x\in\H$. Furthermore, $\theta$ is supported on the set $\mathcal B_r(w)\cup\{-\mathcal B_r(w)\}$.
Next, we show that $\theta$ defined in \eqref{theta} is a positive definite functional on $\H$, that is, for any $m\in\mathbb N_+$, $a_1,\ldots,a_m\in\mathbb R$ and $x_1,\ldots,x_m\in\H$, it is true that
\begin{align}\label{goal3}
\sum_{j,q=1}^m a_ja_q\theta(x_j-x_q)\geq0\,.
\end{align}
In order to prove \eqref{goal3}, for any $t\in\mathbb R$, define
\begin{align}\label{gkt}
G_k(t)&:=\bigg(1-\frac{\pi k}{\sqrt{6}r}|t-\l w,e_k\r|\bigg)_++\bigg(1-\frac{\pi k}{\sqrt{6}r}|t+\l w,e_k\r|\bigg)_+\,,
\end{align}
so that in view of the definition of $\theta(x)$ in \eqref{theta}, we have
\begin{align*}
\theta(x)=\prod_{k=1}^\infty\bigg\{1-\frac{G_k(\l x,e_k\r)}{3}\bigg\}^{-k^{-2}}-1\,.
\end{align*}
Observing that $0\leq|G_k(t)|\leq 1$ for any $t\in\mathbb R$ and $k\in\mathbb N_+$, we deduce that
\begin{align}\label{thetaxx}
\prod_{k=1}^\infty\bigg\{1-\frac{G_k(\l x,e_k\r)}{3}\bigg\}^{-k^{-2}}\leq \prod_{k=1}^\infty\bigg(\frac{2}{3}\bigg)^{-k^2}=\bigg(\frac{3}{2}\bigg)^{\pi^2/6}\,.
\end{align}
We first show that $G_k$ is a positive definite function on $\mathbb R$. To achieve this, for any $k\in\mathbb N_+$, define
\begin{align}\label{psik}
\psi_k(t)=\bigg(1-\frac{\pi k}{\sqrt{6}r}|t|\bigg)_+\,,\qquad\text{for }t\in\mathbb R\,.
\end{align}
Note that the inverse Fourier transform of $\psi$ is given by
\begin{align}\label{t3}
\psi_k^{\rm invft}(s)&:=\frac{1}{2\pi}\int_{-\infty}^{+\infty}\exp(-\i st)\psi(t)\,\d t=\frac{1}{\pi}\int_{0}^{\sqrt{6}r(\pi k)^{-1}}\cos(st)\bigg(1-\frac{\pi k}{\sqrt{6}r}t\bigg)\,\d t\notag\\
&=\frac{\pi k\big[1-\cos \{\sqrt{6}r(\pi k)^{-1}s\}\big]}{\sqrt{6}r\pi s^2}\geq0\,,
\end{align} 
which implies that $\psi_k$ defined in \eqref{psik} is a positive definite function on $\mathbb R$. Now, we deduce that for each $k\in\mathbb N_+$, $G_k(t)$ defined in \eqref{gkt} satisfies
\begin{align*}
G_k(t)&=\int_{-\infty}^{+\infty}\Big[\exp\{\i t (s+\l w,e_k\r)\}+\exp\{\i t (s-\l w,e_k\r)\}\Big]\,\psi^{\rm invft}(s)\,\d s\,,
\end{align*}
so that $G_k$ defines a positive definite function on $\mathbb R$ due to Fourier shift (see for example Lemma~A.1(f) in \cite{meister2009}). To be more specific, in view of \eqref{t3}, we obtain that, for any $r\in\mathbb N_+$ and $b_1,\ldots,b_r,t_1,\ldots,t_r\in\mathbb R$, it holds that
\begin{align*}
&\sum_{j,\ell=1}^rb_jb_\ell\,G_k(t_j-t_\ell)\\
&=\sum_{j,\ell=1}^rb_jb_\ell\int_{-\infty}^{+\infty}\Big[\exp\{\i (t_j-t_\ell) (s+\l w,e_k\r)\}+\exp\{\i (t_j-t_\ell) (s-\l w,e_k\r)\}\Big]\,\psi^{\rm invft}(s)\,\d s\\
&=\int_{-\infty}^{+\infty}\Bigg[\bigg|\sum_{j=1}^rb_j\exp\{\i t_j(s+\l w,e_k\r)\}\bigg|^2+\bigg|\sum_{j=1}^rb_j\exp\{\i t_j(s-\l w,e_k\r)\}\bigg|^2\Bigg]\,\psi^{\rm invft}(s)\,\d s\geq0\,.
\end{align*}
The above equation shows that for each $k\in\mathbb N_+$, $G_k$ defined in \eqref{gkt} is a positive function on $\mathbb R$.

Next, for every $L\in\mathbb N_+$, we define
\begin{align*}
\theta_L(u)=\prod_{k=1}^L\bigg\{1-\frac{G_k(u_k)}{3}\bigg\}^{-k^{-2}}-1\,,\qquad\text{for }u=(u_1,\ldots,u_L)\trans\in\mathbb R^L\,,
\end{align*}
and we claim that $\theta_L$ is a positive definite function on $\mathbb R^L$. This claim is verified by applying Proposition~C.1.6 and the proof of Lemma~C.1.8 of \cite{bekka2008} to the positive definite kernels $K_{G_k}$ defined by $K_{G_k}(t,t')=G_k(t-t')$, for $t,t'\in\mathbb R$ and $k\in\mathbb N_+$.
For any $x\in\H$, let $x^{(L)}=(\l x,e_1\r,\ldots,\l x,e_L\r)\trans\in\mathbb R^L$.
Now, we have obtained that, for any $L\in\mathbb N_+$, we have, for any $m\in\mathbb N_+$, $a_1,\ldots,a_m\in\mathbb R$ and $x_1,\ldots,x_m\in\H$,
\begin{align}\label{goal2}
\sum_{j,q=1}^ma_ja_q\theta_L(x_j^{(L)}-x_q^{(L)})&\geq0\,.
\end{align}
Applying \eqref{thetaxx} and the fact that $1-e^{-x}\leq x$ for all $x\in\mathbb R$, we further deduce that
\begin{align*}
&\sum_{j,q=1}^m a_ja_q\{\theta_L(x_j^{(L)}-x_q^{(L)})-\theta(x_j-x_q)\}\\
&=\sum_{j,q=1}^m a_ja_q\Bigg[\prod_{k=1}^L\bigg\{1-\frac{G_k(\l x_j-x_q,e_k\r)}{3}\bigg\}^{-k^{-2}}\Bigg]\Bigg[1-\prod_{k=L+1}^\infty\bigg\{1-\frac{G_k(\l x_j-x_q,e_k\r)}{3}\bigg\}^{-k^{-2}}\Bigg]\\
&\leq(3/2)^{\pi^2/6}\sum_{j,q=1}^m a_ja_q\Bigg[1-\prod_{k=L+1}^\infty\bigg\{1-\frac{G_k(\l x_j-x_q,e_k\r)}{3}\bigg\}^{-k^{-2}}\Bigg]\\
&=(3/2)^{\pi^2/6}\sum_{j,q=1}^m a_ja_q\Bigg(1-\exp\bigg[-\sum_{k=L+1}^\infty\frac{1}{k^2}\log\bigg\{1-\frac{G_k(\l x_j-x_q,e_k\r)}{3}\bigg\}\bigg]\Bigg)\\
&\leq (3/2)^{\pi^2/6}\log(3/2)\bigg(\sum_{j=1}^m a_j\bigg)^2\sum_{k=L+1}^\infty\frac{1}{k^2}=o(1)
\end{align*}
as $L\to\infty$. Combining the above equation with \eqref{goal2} proves \eqref{goal3}, which implies that $\theta$ defined in \eqref{theta} is a positive definite functional.

Next, let
\begin{align}\label{phi1}
\phi_1(x)=\exp\bigg(-\frac{1}{2}\sum_{k=1}^\infty\frac{1}{k^2}\l x,e_k\r^2\bigg)
\end{align}
denote the characteristic functional of a mean-zero Gaussian measure on $\H$.
Since $0\notin\mathcal B_r(w)\cup\{-\mathcal B_r(w)\}$, there exists a constant $c_0>0$ such that for any $x\in\mathcal B_r(w)\cup\{-\mathcal B_r(w)\}$, it holds that $\sum_{k=1}^\infty k^{-2}\l x,e_k\r^2\geq c_0$, so that $\phi_1(x)\geq e^{-c_0}$. Observe from \eqref{thetaxx} that $\theta(x)\leq(3/2)^{\pi^2/6}-1$, for any $x\in\mathcal B_r(w)\cup\{-\mathcal B_r(w)\}$.
By taking 
\begin{align*}
a=2^{-1}(1-e^{-c_0})\{(3/2)^{\pi^2/6}-1\}^{-1}>0\,,
\end{align*}
and, for $\theta$ defined in \eqref{theta}, taking
\begin{align}\label{phi2}
\phi_2(x)=\phi_1(x)+a\theta(x)\,,\qquad ~~~x\in\H\,,
\end{align}
we have that $\phi_2(x)\leq1$ for any $x\in\H$. Since $\theta\geq0$ and $\phi_1$ is a characteristic functional of a probability measure on $\H$, then $\phi_2$ satisfies condition~(b) in Theorem~\ref{thm:ms}, which justifies that $\phi_2$ is also a well-defined characteristic functional of a probability measure on $\H$. Now, we have $\ed_\nu^2(P_1,P_2)=\int_\H|\phi_1(w)-\phi_2(w)|^2\nu(\d w)=0$ but $P_1\neq P_2$. This contradiction implies that $\supp(\nu)=\H$.

\subsection{Proof of Corollary~\ref{cor:4.11}}\label{app:cor:4.11}

In view of Theorem~\ref{thm:iff:ed}, to prove Corollary~\ref{cor:4.11}, it suffices to verify that $\ed_\nu$ defined in \eqref{ed} satisfies the triangle inequality. Let $F_1,F_2,F_3\in\mathcal P$ and let $\phi_1,\phi_2,\phi_3$ denote their corresponding characteristic functionals. Direct calculations yield
\begin{align*}
&\ed_\nu^2(F_1,F_2)=\int_\H|\phi_1(w)-\phi_2(w)|^2\,\nu(\d w)=\int_\H|\{\phi_1(w)-\phi_3(w)\}-\{\phi_1(w)-\phi_2(w)\}|^2\,\nu(\d w)\\
&=\int_\H\Big(|\phi_1(w)-\phi_3(w)|^2+|\phi_2(w)-\phi_3(w)|^2+\Re\big[\{\phi_1(w)-\phi_3(w)\}\{\overline{\phi_1(w)-\phi_2(w)}\}\big]\Big)\,\nu(\d w)\\
&\leq\int_\H|\phi_1(w)-\phi_3(w)|^2\nu(\d w)+\int_\H|\phi_2(w)-\phi_3(w)|^2\nu(\d w)+2\int_\H|\phi_1(w)-\phi_3(w)||\phi_1(w)-\phi_2(w)|\nu(\d w)\,.
\end{align*}
It follows from the above equation and the Cauchy-Schwarz inequality that
\begin{align*}
\ed_\nu^2(F_1,F_2)&\leq \int_\H|\phi_1(w)-\phi_3(w)|^2\nu(\d w)+\int_\H|\phi_2(w)-\phi_3(w)|^2\nu(\d w)\\
&\qquad+2\bigg\{\int_\H|\phi_1(w)-\phi_3(w)|^2\nu(\d w)\bigg\}^{1/2}\bigg\{\int_\H|\phi_2(w)-\phi_3(w)|^2\nu(\d w)\bigg\}^{1/2}\\
&=\{\ed_\nu(F_1,F_3)+\ed_\nu(F_2,F_3)\}^2\,,
\end{align*}
which implies $\ed_\nu(F_1,F_2)\leq\ed_\nu(F_1,F_3)+\ed_\nu(F_2,F_3)$ and concludes the proof.

\subsection{Proof of Theorem~\ref{thm:laplace}}\label{app:thm:laplace}


For sufficiency, suppose $\nu$ is the Laplace measure on $\H$ with mean $\mu$ and covariance operator $C$. For $x\in\H$, let $F_x$ denote the distribution of $\l x,\cdot\r$, which follows the Laplace distribution on $\mathbb R$ with mean $\mu_x$ and variance $2b_x^2$. We have
\begin{align*}
\phi(x)&=\int_{\H }\exp(\i\l x,y\r)\,\nu(\d y)=\int_{\mathbb R}\exp(\i t)\,F_x(\d t)=\frac{1}{2b_x}\int_{\mathbb R}\exp\Big(\i t-\frac{|t-\mu_x|}{b_x}\Big)\,\d t\,,
\end{align*}
Then, an application of change of variables and integration by parts yields
\begin{align}\label{o1}
\phi(x)&=\frac{1}{2b_x}\exp(\i\mu_x)\int_{\mathbb R}\exp(\i t-b_x^{-1}|t|)\,\d t=\exp(\i\mu_x)(1+b_x^2)^{-1}\,.
\end{align}
On the other hand, both $\mu_x,b_x$ are finite, and recall  that $\mu_x$ and $2b_x^2$ are the mean and variance of $F_x$, respectively, such that applying Definition~\ref{def:meanvar} yields
\begin{align}\label{mub}
&\mu_x=\int_{\mathbb R}t\,F_x(dt)=\int_\H \l x,y\r\,\nu(\d y)=\l \mu,x\r\,,\notag\\
&2b_x^2=\int_{\mathbb R}(t-\mu_x)^2\,F_x(dt)=\int_\H \l x,y-\mu_x\r^2\,\nu(\d y)=\int_\H \l x,y\r^2\,\nu(\d y)=\l Cx,x\r\,.
\end{align}

It follows from the Minlos–Sazonov theorem (Theorem~\ref{thm:ms}) that there exists an $\mathcal S$-perator $S_\e$ on $\H$ such that for any $\e\in(0,1/2)$, $1-\Re\{\phi(x)\}\leq\l S_\epsilon x,x\r+\epsilon$ for all $x\in\H$. 
We fix this $\epsilon$.
Let $\{\tau_j\}$ denote the non-zero eigenvalues of $S_\e$, and let $\{\eta_j\}$ denote their corresponding eigenvectors. In addition, let $\{f_j\}$ be an orthonormal basis of the $ker(S_\e)$. For any $x\in\H$, define
\begin{align}\label{C}
T(x)=\sum_{j}\tau_j\l x,\eta_j\r \eta_j+\sum_{j}j^{-2}\l x,f_j\r f_j\,.
\end{align}
Then, we have
\begin{align}\label{o4}
\sum_{j}\l T(\eta_j),\eta_j\r+\sum_{j}\l T(f_j),f_j\r=\sum_{j}\tau_j+\sum_{j}j^{-2}<\infty\,.
\end{align}
In addition, $ker(T)=\{0\}$, and, for any $x\in\H$,
\begin{align}\label{o3}
\l S_\epsilon(x),x\r=\sum_{j}\tau_j\l x,\eta_j\r^2\leq\sum_{j}\tau_j\l x,\eta_j\r^2+\sum_{j}j^{-2}\l x,f_j\r^2=\l T(x),x\r\,.
\end{align}
From \eqref{o1} we obtain that $\phi(x)=\exp(\i\mu_x)(1+b_x^2)^{-1}$, so that $|\phi(x)|=(1+b_x^2)^{-1}$.
Combining the above equation with \eqref{o3} yields
\begin{align*}
1-(1+b_x^2)^{-1}=1-|\phi(x)|\leq  1-\text{Re}\{\phi(x)\}\leq \l S_\epsilon (x),x\r+\epsilon\leq\l T (x),x\r+\epsilon\,.
\end{align*}
The above equation implies that
\begin{align}\label{o2}
b_x^2\leq\frac{4\epsilon}{1-2\epsilon}\,,\quad\text{for any }x\in\H\text{ such that }\l S_\epsilon x,x\r\leq\e\,.
\end{align}
Fix an arbitrary $x\in\H $, $x\neq0$, let $y=2^{-1/2}\l T(x),x\r^{-1/2}\epsilon^{1/2}x$ (note that the case where $x=0$ and $b_x^2=0$ is trivial). Then it is straightforward to verify $\l T(y),y\r=\e/2<\e$. Therefore, it follows from \eqref{o2} that
\begin{align}\label{cy}
b_y^2\leq\frac{4\epsilon}{1-2\epsilon}\,.
\end{align}
On the other hand, in view of Definition~\ref{def:meanvar}, we have
\begin{align*}
b_y^2=\int_{\H }\l y,z-\mu\r^2\,\nu(\d z)=2^{-1}\l T( x),x\r^{-1}\epsilon\int_{\H }\l x,z-\mu\r^2\,\nu(\d z)=2^{-1}\l T( x),x\r^{-1}\epsilon b_x^2\,.
\end{align*}
Combining the above equation with  \eqref{mub} and \eqref{cy} yields that
\begin{align*}
b_x^2=\frac 1 2\l Cx,x\r\leq \frac{8\l T(x),x\r}{1-2\epsilon}<\infty\,,\quad\text{for any }x\in\H\,.
\end{align*}
The finiteness of $b_x^2$ for any $x\in\H$ implies the existence of covariance operator $C$. Furthermore, in view of \eqref{o4}, it follows that
\begin{align*}
\sum_{j}\l C \eta_j, \eta_j\r+\sum_{j}\l C f_j, f_j\r\leq \frac{16}{1-2\e}\bigg\{\sum_{j}\l T(\eta_j),\eta_j\r+\sum_{j}\l T(f_j), f_j\r\bigg\}<\infty\,.
\end{align*}
which implies that $C$ is an $\mathcal S$-class operator.

Next, we prove necessity. First, $\phi(0)=1$. We show that $\phi$ is a positive definite functional. To see this, let $\{\lambda_\ell\}$ and $\{e_\ell\}$ be the eigenvalues and eigenvectors of $C$. We aim to show that for any arbitrary but fixed $m\in\mathbb N_+$,  $a_1,\ldots,a_m\in\mathbb R$ and $x_1,\ldots,x_m\in\H$, it holds that
\begin{align}\label{te1}
\sum_{j,k=1}^ma_ja_k\phi(x_j-x_k)=\sum_{j,k=1}^ma_ja_k\Big(1+\frac{1}{2}\sum_{\ell=1}^\infty\lambda_\ell\l x_j-x_k,e_\ell\r^2\Big)^{-1}\geq0\,.
\end{align}
For any $r\in\mathbb N_+$, let $x_r=(\l x,e_1\r,\ldots,\l x,e_r\r)\trans\in\mathbb R^r$ and $\mu_r=(\l\mu,e_1\r,\ldots,\l\mu,e_r\r)\trans\in\mathbb R^r$. Define
\begin{align*}
\phi_r(x)=\exp(\i\mu_r\trans x_r)(1+2^{-1}x_r\trans C_r x_r)^{-1}=\int_{\mathbb R^r}\exp(\i x_r\trans w)\,F_r(\d w)\,,
\end{align*}
where $F_r$ denotes the $r$-dimensional Laplace distribution with mean vector $\mu_r$ and  covariance matrix $C_r=\diag(\lambda_1,\ldots,\lambda_r)$. 
Let $x_{j,r}=(x_{j}^{(1)},\ldots,x_j^{(r)})\trans\in\mathbb R^r$, where $x_{j}^{(\ell)}=\l x_j,e_\ell\r$. 
It holds that
\begin{align}\label{te0}
\sum_{j,k=1}^ma_ja_k\phi_r(x_{j,r}-x_{k,r})
&=\sum_{j,k=1}^ma_ja_k\exp\Big(\i\sum_{\ell=1}^r\l\mu,e_\ell\r\l x,e_\ell\r\Big)\Big(1+\frac{1}{2}\sum_{\ell=1}^r\lambda_\ell\l x_j-x_k,e_\ell\r^2\Big)^{-1}\notag\\
&=\int_{\mathbb R^r}\Big|\sum_{j=1}^ma_j\exp(\i x_{j,r}\trans w)\Big|^2\,F_r(\d w)
\geq0\,.
\end{align}
Since $\sum_{\ell=1}^\infty\lambda_\ell<\infty$, it holds that, as $r\to\infty$,
\begin{align*}
&\bigg|\sum_{j,k=1}^m|a_j||a_k||\phi_r(x_{j,r}-x_{k,r})|-\sum_{j,k=1}^m|a_j||a_k||\phi(x_j-x_k)|\bigg|\\
&=\bigg|\sum_{j,k=1}^m|a_j||a_k|\Big(1+\frac{1}{2}\sum_{\ell=1}^\infty\lambda_\ell\l x_j-x_k,e_\ell\r^2\Big)^{-1}-\sum_{j,k=1}^m|a_j||a_k|\Big(1+\frac{1}{2}\sum_{\ell=1}^r\lambda_\ell\l x_j-x_k,e_\ell\r^2\Big)^{-1}\bigg|\\
&\leq \sum_{j,k=1}^\infty\sum_{\ell=r}^\infty|a_j||a_k|\lambda_\ell\l x_j-x_k,e_\ell\r^2\leq\sum_{j,k=1}^\infty|a_j||a_k|\|x_j-x_k\|^2\times\sum_{\ell=r}^\infty\lambda_\ell=o(1)\,.
\end{align*}
The above equation together with \eqref{te0} implies \eqref{te1}.

Consider the case where $\mu=0$. Note that $1-(1+a)^{-1}\leq a$ for any $a\geq0$, such that
\begin{align*}
1-\text{Re}\{\phi(x)\}=1-\Big(1+\frac{1}{2}\l C x,x\r\Big)^{-1}\leq\frac{1}{2}\l C x,x\r\,,
\end{align*}
where we used the fact that $\l C x,x\r=\int_\H\l x,y\r^2\nu(\d y)\geq0$ in view of Definition~\ref{def:meanvar}.
Therefore, by Theorem~\ref{thm:ms}, there exists a Borel measure $\omega$ on $\H$ such that $\phi=\omega^\ft$, that is,
\begin{align*}
\int_{\H }\exp(\i\l x,y\r)\,\omega(\d y)=\phi(x)=\Big(1+\frac{1}{2}\l Cx,x\r\Big)^{-1}\,,\qquad \text{for any }x\in\H \,.
\end{align*}
This implies that $F_x$ is the distribution of $\l x,\cdot\r$, so that
\begin{align*}
\int_{-\infty}^{+\infty}e^{\i t}\,F_x(\d t)=\Big(1+\frac{1}{2}\l Cx,x\r\Big)^{-1}\,.
\end{align*}
The above equation further implies that $F_x$ is a Laplace distribution with mean $0$ and variance $2b_x^2=\l Cx,x\r$. This verifies that $\mu$ is a Lpalace measure on $\H $. Next, consider the case $\mu\neq0$. Let $$\psi(x)=\Big(1+\frac{1}{2}\l Sx,x\r\Big)^{-1}\,,$$
such that $\phi(x)=\exp(\i\l \mu,x\r)\psi(x)$. Following the argument similar to the ones to deal with the case $\mu=0$, there exists a Laplace measure $\tilde\omega$ in $\H $ such that $\psi=\tilde\omega^\ft$. Now, define a Borel measure $\nu$ in $\H$ as follows: $\mu(E)=\nu(E-\mu)$, for any $\E\in\mathcal B(\H )$.
Then, $\nu$ is a Laplace measure on $\H$ and $\tilde w^\ft(x)=\exp(\i\l x_0,x\r)\phi_{\nu}(x)=\exp(\i\l x_0,x\r)\psi(x)=\phi(x)$.

\subsection{Proof of Proposition~\ref{prop:support}}\label{app:prop:support}

In order to prove $\supp(\nu_L)=\H$, it suffices to show that for any closed ball
\begin{align*}
\mathcal B(w,\sqrt{2}r)=\{x\in\H:\|x-w\|^2\leq 2r^2\}=\bigg\{x\in\H:\sum_{\ell=1}^\infty\l x-w,e_\ell\r^2\leq 2r^2\bigg\}
\end{align*}
centered at $w\in\H$ with $\|\cdot\|$-radius $\sqrt{2}r$. For any $M\in\mathbb N_+$, define the subsets of $\H$
\begin{align*}
A_M=\bigg\{x\in\H:\sum_{\ell=1}^M\l x-w,e_\ell\r^2\leq r^2\bigg\}\,,\qquad B_M=\bigg\{x\in\H:\sum_{\ell=M+1}^\infty\l x-w,e_\ell\r^2\leq r^2\bigg\}\,.
\end{align*}

Since $C_{\nu_L}$ is non-degenerate, then its eigenvalues $\{\lambda_k\}_{k=1}^\infty$ are all positive. Let
\begin{align*}
    c_M=\frac{\P\big[\sum_{k=1}^M \{\sqrt{\lambda_k}(L_{k}+L'_{k})-w_k\}^2\leq r^2\big]}{\P\big\{\sum_{k=1}^M (\sqrt{\lambda_k}L_{k}-w_k)^2\leq r^2\big\}},
\end{align*}
so that by assumption, $\lim\sup_{M\to\infty}c_M<\infty$.
Observe that
\begin{align}\label{mi}
\nu\{\mathcal B(w,r)\}\geq\nu(A_M\cap B_M)=\nu(A_M)-\nu(A_M\cap B_M^c)
\end{align}
Let $\tilde f_M$ and $f_M$ denote the probability density of $M$-dimensional mean-zero Laplace distribution with covariance matrix $\diag(\lambda_1,\ldots,\lambda_M)$ and $I_M$, respectively. Then, by a change of variables we obtain
\begin{align}\label{nua}
\nu(A_M)
&=\int_\H\one\bigg\{\sum_{k=1}^M\l x-w,e_k\r^2\leq r^2\bigg\}\,\nu(\d x)\notag\\
&=\int_{\mathbb R^M}\tilde f_M(x_1,\ldots,x_M)\,\one\bigg\{\sum_{k=1}^M(x_k-\l w,e_k\r)^2\leq r^2\bigg\}\,\d x_1\ldots \d x_M\notag\\
&=\int_{\mathbb R^M} f_M(x_1,\ldots,x_M)\,\one\bigg\{\sum_{k=1}^M(\sqrt{\lambda_k} x_k-\l w,e_k\r)^2\leq r^2\bigg\}\,\d x_1\ldots \d x_M\notag\\
&=\P\bigg\{\sum_{k=1}^M (\sqrt{\lambda_k}\mathfrak L_{1k}-\l w,e_k\r)^2\leq r^2\bigg\}\,.
\end{align}
It is true that $\nu(A_M)>0$.

In addition, for the minuend in \eqref{mi}, by Chebychev's inequality, we have
\begin{align}\label{nua2}
\nu(A_M\cap B_M^c)&=\nu\bigg(\bigg\{x\in\H:\sum_{k=1}^M\l x-w,e_k\r^2\leq r^2\bigg\}\cap\bigg\{x\in\H:\sum_{\ell=M+1}^\infty\l x-w,e_\ell\r^2\geq r^2\bigg\}\bigg)\notag\\
%
%
&\leq\frac{1}{r^2}\sum_{\ell=M+1}^\infty\int_\H\l x-w,e_\ell\r^2\,\one\bigg\{\sum_{k=1}^M\l x-w,e_k\r^2\leq r^2\bigg\}\,\nu(\d x)\notag\\
&\leq\frac{2}{r^2}\sum_{\ell=M+1}^\infty\int_\H\big(\l x,e_\ell\r^2+\l w,e_\ell\r^2\big)\,\one\bigg\{\sum_{k=1}^M\l x-w,e_k\r^2\leq r^2\bigg\}\,\nu(\d x)\notag\\
&\leq\frac{2\nu(A_M)}{r^2}\bigg(\sum_{\ell=M+1}^\infty\l w,e_\ell\r^2+c_M\bigg)\,,
%
\end{align}
where, by the change of variables,
\begin{align}\label{cm}
&c_M=\frac{1}{\nu(A_M)}\sum_{\ell=M+1}^\infty\int_\H\l x,e_\ell\r^2\,\one\bigg\{\sum_{k=1}^M\l x-w,e_k\r^2\leq r^2\bigg\}\,\nu(\d x)\\
&=\frac{1}{\nu(A_M)}\sum_{\ell=M+1}^\infty\int_{\mathbb R^{M+1}}x_\ell^2 \tilde f_{\ell}(x_1,\ldots,x_M,x_\ell)\,\one\bigg\{\sum_{k=1}^M(x_k-\l w,e_k\r)^2\leq r^2\bigg\}\,\d x_1\ldots\d x_M\,\d x_\ell\notag\\
&=\frac{1}{\nu(A_M)}\sum_{\ell=M+1}^\infty\lambda_\ell\int_{\mathbb R^{M+1}}x_\ell^2f_{\ell}(x_1,\ldots,x_M,x_\ell)\,\one\bigg\{\sum_{k=1}^M(\sqrt{\lambda_k}x_k-\l w,e_k\r)^2\leq r^2\bigg\}\,\d x_1\ldots\d x_M\,\d x_\ell\notag\\
&=\frac{1}{\nu(A_M)}\sum_{\ell=M+1}^\infty\lambda_\ell\int_{\mathbb R^{M}}\bigg[\int_{\mathbb R} x_\ell^2f_{\ell}(x_1,\ldots,x_M,x_\ell)\,\d x_\ell\bigg]\,\one\bigg\{\sum_{k=1}^M(\sqrt{\lambda_k}x_k-\l w,e_k\r)^2\leq r^2\bigg\}\,\d x_1\ldots\d x_M\,,\notag
\end{align}
where $\tilde f_\ell$ and $f_\ell$ denote the probability density function of $(M+1)$-dimensional mean-zero Laplace distribution with covariance matrix $\diag(\lambda_1,\ldots,\lambda_M,\lambda_\ell)$ and $I_{M+1}$, respectively. Note that the characteristic function of $f_\ell$ is given by
\begin{align*}
f_\ell^\ft(t_1,\ldots,t_M,t_\ell)=\bigg\{1+\frac{1}{2}\bigg( t_\ell^2+\sum_{k=1}^Mt_k^2\bigg)\bigg\}^{-1}\,.
\end{align*}
Note that the Fourier transform
\begin{align*}
&\bigg\{\int_{\mathbb R} x_\ell^2f_{\ell}\,\d x_\ell\bigg\}^\ft(t_1,\ldots,t_M)=\int_{\mathbb R^M}\exp\bigg(\i\sum_{k=1}^Mt_kx_k\bigg)\bigg\{\int_{\mathbb R} x_\ell^2f_{\ell}(x_1,\ldots,x_M,x_\ell)\,\d x_\ell\bigg\}\,\d x_1\ldots\d x_M\\
&\hspace{4cm}=-\frac{\partial^2 f_\ell^\ft}{\partial t_\ell^2}(t_1,\ldots,t_M,0)=\bigg(1+\frac{1}{2}\sum_{k=1}^Mt_k^2\bigg)^{-2}=\{f_M^\ft(t_1,\ldots,t_M)\}^2\,.
\end{align*}
Therefore, we obtain
\begin{align}
\int_{\mathbb R} x_\ell^2f_{\ell}(x_1,\ldots,x_M,x_\ell)\,\d x_\ell=f_M*f_M(x_1,\ldots,x_M)\,,
\end{align}
where $*$ denotes the convolution operator. Combining the above equation with \eqref{cm} yields that, for $M$ large enough,
\begin{align*}
c_M&=\frac{1}{\nu(A_M)}\int_{\mathbb R^{M}}f_M*f_M(x_1,\ldots,x_M)\,\one\bigg\{\sum_{k=1}^M(\sqrt{\lambda_k}x_k-\l w,e_k\r)^2\leq r^2\bigg\}\,\d x_1\ldots\d x_M\sum_{\ell=M+1}^\infty\lambda_\ell\\
&=\frac{\P\big[\sum_{k=1}^M \{\sqrt{\lambda_k}(\mathfrak L_{1k}+\mathfrak L_{2k})-\l w,e_k\r\}^2\leq r^2\big]}{\P\big\{\sum_{k=1}^M (\sqrt{\lambda_k}\mathfrak L_{1k}-\l w,e_k\r)^2\leq r^2\big\}}\sum_{\ell=M+1}^\infty\lambda_\ell\leq c_0\sum_{\ell=M+1}^\infty\lambda_\ell\,,
\end{align*}
where by condition \eqref{pp2}, $c_0>0$ does not depends on $M$.

Combining the above equation with \eqref{mi}--\eqref{nua2} yields
\begin{align}\label{c}
\nu\{\mathcal B(w,r)\}\geq\nu(A_M)\bigg\{1-\frac{2}{r^2}\bigg(\sum_{\ell=M+1}^\infty\l w,e_\ell\r^2+c_0\sum_{\ell=M+1}^\infty\lambda_\ell\bigg)\bigg\}
\end{align}
In addition, we have $$\sum_{\ell=M+1}^\infty\l w,e_\ell\r^2+c_0\sum_{\ell=M+1}^\infty\lambda_\ell\leq \frac{r^2}{4}$$ for $M$ large enough, since $\sum_{\ell=1}^\infty\l w,e_\ell\r^2<\infty$ and $\sum_{\ell=1}^\infty\lambda_\ell<\infty$. Therefore, we conclude from \eqref{c} that $\nu\{\mathcal B(w,r)\}\geq 2^{-1}\nu(A_M)>0$.

\subsection{Proof of Proposition~\ref{prop:mmded}}\label{proof:prop:mmded}

Since $K_\nu$ defined in \eqref{knu2} is bounded, then it holds that $\E_{W\sim P}\sqrt{K_\nu(W,W)}<\infty$ for any probability measure $P$ on $\H$.
Observe that $\|K_\nu(x,\cdot)\|_{\H_{K_\nu}}=\sqrt{K_\nu(x,x)}$ for any $x\in\H_{K_\nu}$.
Let $\mathfrak S_P(f)=\E_{W\sim P} f(W)$ denote a linear functional that maps $\H_{K_\nu}$ to $\mathbb R$.
For any $f\in\H_{K_\nu}$ and probability measure on $\H$, by the Cauchy-Schwarz inequality, it follows that
\begin{align*}
\mathfrak S_P(f)=\E_{W\sim P}\l f,K_\nu(W,\cdot)\r_{\H_{K_\nu}}\leq\|f\|_{\H_{K_\nu}}\E_{W\sim P}\|K_\nu(W,\cdot)\|_{\H_{K_\nu}}=\|f\|_{\H_{K_\nu}}\E_{W\sim P}\sqrt{K_\nu(W,W)}\,.
\end{align*}
Since $\E_{W\sim P}\sqrt{K_\nu(W,W)}<\infty$, the above equation implies that $\mathfrak S_P$ is a bounded linear functional on $\H_{K_\nu}$.
By the Riesz representation theorem, there exists a unique element $\mu_P\in\H_{K_\nu}$ such that $\E_{W\sim P} f(W)=\mathfrak S_P(f)=\l \mu_P,f\r_{\H_{K_\nu}}$. 
Note that $\mu_P(z)=\l\mu_P,K_\nu(z,\cdot)\r_{\H_{K_\nu}}=\E_{W\sim P} K_\nu(W,z)$.
Now, by the Cauchy-Schwarz inequality,
\begin{align}\label{mmdp}
&\mmd_\nu^2(P_1,P_2;\H_{K_\nu,1})=\sup_{f\in\H_{K_\nu},\|f\|_{\H_{K_\nu}}\leq1}|\mathfrak S_{P_1}(f)-\mathfrak S_{P_2}(f)|^2=\sup_{f\in\H_{K_\nu},\|f\|_{\H_{K_\nu}}\leq1}|\l \mu_{P_1}-\mu_{P_2},f\r_{\H_{K_\nu}}|^2\notag\\
&\hspace{2cm}=\|\mu_{P_1}-\mu_{P_2}\|_{\H_{K_\nu}}^2=\|\E_{W_1\sim P_1} K_\nu(W_1,\cdot)-\E_{W_2\sim P_2} K_\nu(W_2,\cdot)\|_{\H_{K_\nu}}^2\notag\\
&\hspace{2cm}=\|\E K_\nu(W_1,\cdot)\|_{\H_{K_\nu}}^2+\|\E K_\nu(W_2,\cdot)\|_{\H_{K_\nu}}^2-2\big\l \E K_\nu(W_1,\cdot),\E K_\nu(W_2,\cdot)\big\r_{\H_{K_\nu}}\,.
\end{align}
Suppose $W_1,W_1'\iidsim P_1$ and $W_2,W_2'\iidsim P_2$, it follows from Fubini's theorem that,
\begin{align*}
&\|\E K_\nu(W_j,\cdot)\|_{\H_{K_\nu}}^2=\E\l K_\nu(W_j,\cdot), K_\nu(W_j',\cdot)\r_{\H_{K_\nu}}=\E K_\nu(W_j,W_j')=\E K_\nu(W_j,W_j')\,,\quad\text{for }j=1,2\,;\\
&\l \E K_\nu(W_1,\cdot),\E K_\nu(W_2,\cdot)\r_{\H_{K_\nu}}=\E\l K_\nu(W_1,\cdot), K_\nu(W_2,\cdot)\r_{\H_{K_\nu}}=\E K_\nu(W_1,W_2)=\E K_\nu(W_1,W_2)\,.
\end{align*}
The proof is therefore complete by combining the above equations with \eqref{mmdp}.

\subsection{Proof of Proposition~\ref{thm:dcov}}\label{app:thm:dcov}
We derive the four equivalent formulas for $\dcov_\nu^2(X,Y)$.
First, observe that $\overline{\phi_{X,Y}(w_1,w_2)}=\phi_{X,Y}(-w_1,-w_2)$. Since $\nu$ is symmetric, we obtain that $\dcov_{\nu}^2(X,Y)=I_1+I_2 -2I_3 $, where
{
\begin{align*}
&I_1=\int_{\H^2}\phi_{X,Y}(w_1,w_2)\phi_{X,Y}(-w_1,-w_2)\nu(\d w_1)\nu(\d w_2)\,;\\
&I_2=\int_{\H^2}\phi_{X}(w_1)\phi_{X}(-w_1)\phi_{Y}(w_2)\phi_{Y}(-w_2)\,\nu(\d w_1)\nu(\d w_2)\,;\\
&I_3=\int_{\H^2}\phi_{X,Y}(w_1,w_2)\phi_{X}(-w_1)\phi_{Y}(-w_2)\nu(\d w_1)\nu(\d w_2)\,.
\end{align*}
}For the first term $I_1$, since $(X,Y)$ and $(X',Y')$ are independent and identically distributed, we have
\begin{align*}
\phi_{X,Y}(w_1,w_2)\phi_{X,Y}(-w_1,-w_2)&=\E\exp(\i\l X,w_1\r+\i\l Y,w_2\r)\E\exp(-\i\l X',w_1\r-\i\l Y',w_2\r)\\
&=\E\exp(\i\l X-X',w_1\r+\i\l Y-Y',w_2\r)\,.
\end{align*}
Since $|\phi_{\nu}(x)|\leq1$ for any $x\in\H$, applying Fubini's theorem yields
\begin{align*}
I_1&=\E\Big\{\int_{\H^2}\exp(\i\l X-X',w_1\r+\i\l Y-Y',w_2\r)\nu(\d w_1)\nu(\d w_2)\Big\}\notag\\
&=\E\Big\{\int_{\H }\exp(\i\l X-X',w_1\r)\,\nu(\d w_1)\int_{\H }\exp(\i\l Y-Y',w_2\r)\,\nu(\d w_2)\Big\}=\E\{\phi_{\nu}(X-X')\phi_{\nu}(Y-Y')\}\,.
\end{align*}
Following similar derivation, since $\phi_{\nu}$ is a symmetric functional, we obtain that, for $w_1,w_2\in\H$,
\begin{align*}
&I_2=\int_{\H}\E\exp(\i\l X-X',w_1\r)\nu(\d w_1)\int_\H\E\exp(\i\l Y-Y',w_2\r)\nu(\d w_2)=\E\{\phi_{\nu}(X-X')\}\E\{\phi_{\nu}(Y-Y')\}\,,\\
&I_3=\E\bigg\{\int_{\H }\exp(\i\l X'-X,w_1\r)\nu(\d w_1)\int_{\H }\exp(\i\l Y''-Y,w_2\r)\nu(\d w_2)\bigg\}=\E\{\phi_{\nu}(X-X')\phi_{\nu}(Y-Y'')\}\,.
\end{align*}

Second, observe that
{
\begin{align*}
&|\phi_{X,Y}(w_1,w_2)-\phi_{X}(w_1)\phi_{Y}(w_2)|^2=\big|\E\exp(\i\l X,w_1\r+\i\l Y,w_2\r)-\E\exp(\i\l X,w_1\r)\E\exp(\i\l Y,w_2\r)\big|^2\notag\\
&=\big|\E\big[\{\exp(\i\l X,w_1\r)-\E\exp(\i\l X,w_1\r)\}\{\exp(\i\l Y,w_2\r)-\E\exp(\i\l Y,w_2\r)\}\big]\big|^2\,.
\end{align*}
}
Since $(X,Y)$ and $(X',Y')$ are i.i.d., we obtain from the above equation that
\begin{align*}
&|\phi_{X,Y}(w_1,w_2)-\phi_{X}(w_1)\phi_{Y}(w_2)|^2=\E\big[\{\exp(\i\l X,w_1\r)-\E\exp(\i\l X,w_1\r)\}\{\exp(\i\l Y,w_2\r)-\E\exp(\i\l Y,w_2\r)\}\big]\notag\\
&\times\E\big[\{\exp(-\i\l X,w_1\r)-\E\exp(-\i\l X,w_1\r)\}\{\exp(-\i\l Y,w_2\r)-\E\exp(-\i\l Y,w_2\r)\}\big]
%
=\E\{J_1(w_1)J_2(w_2)\}\,,
\end{align*}
where, for $w_1,w_2\in\H$,
\begin{align*}
J_1(w_1)&=\{\exp(\i\l X,w_1\r)-\E\exp(\i\l X,w_1\r)\}\{\exp(-\i\l X',w_1\r)-\E\exp(-\i\l X',w_1\r)\}\,,\\
J_2(w_2)&=\{\exp(\i\l Y,w_2\r)-\E\exp(\i\l Y,w_2\r)\}\{\exp(-\i\l Y',w_2\r)-\E\exp(-\i\l Y',w_2\r)\}\,.
\end{align*}
By Fubini's theorem, we obtain
\begin{align*}
\dcov_\nu^2(X,Y)&=\int_{\H^2} \E\{J_1(w_1)J_2(w_2)\}\,\nu(\d w_1)\,\nu(\d w_2)=\E\bigg\{\int_{\H }J_1(w_1)\,\nu(\d w_1)\int_\H J_2(w_2)\,\nu(\d w_2)\bigg\}\,.
\end{align*}
For the first term $J_1(w)$, observe that
\begin{align*}
&J_1(w_1)=\exp(\i\l X-X',w_1\r)+\E\exp(\i\l X-X',w_1\r)\notag\\
&\qquad-\exp(\i\l X,w_1\r)\E\exp(-\i\l X',w_1\r)-\exp(-\i\l X',w_1\r)\E\exp(\i\l X,w_1\r)\notag\\
&=\exp(\i\l X-X',w_1\r)+\E\exp(\i\l X-X',w_1\r)-\E_{X'}\exp(\i\l X-X',w_1\r)-\E_{X}\exp(\i\l X-X',w_1\r)\,.
\end{align*}
Therefore, since $\phi_{\nu}$ is symmetric, we conclude the proof by applying Fubini's theorem to obtain,
\begin{align}\label{dc1}
&\int_{\H }J_1(w_1)\,\nu(\d w_1)
%
=\phi_{\nu}(X-X')+\E\phi_{\nu}(X-X')-\E_{X'}\phi_{\nu}(X-X')-\E_{X}\phi_{\nu}(X-X')\notag\\
&\hspace{3cm}=\phi_{\nu}(X-X')+\E\phi_{\nu}(X-X')-2\E_{X}\phi_{\nu}(X-X')\,,\notag\\
&\text{and similarly }\int_{\H }J_2(w_2)\,\nu(\d w_2)=\phi_{\nu}(Y-Y')+\E\phi_{\nu}(Y-Y')-2\E_{Y}\phi_{\nu}(Y-Y')\,.
\end{align}

Third, observing that $\phi_{\nu}(x-y)=1- d_\nu(x,y)$, we deduce that
\begin{align*}
\dcov_\nu^2(X,Y)
&=1-\E d_\nu(X,X')-\E d_\nu(Y,Y')+\E\{ d_\nu(X,X') d_\nu(Y,Y')\}\\
&\quad-2\Big[1-\E d_\nu(X,X')-\E d_\nu(Y,Y')+\E\{ d_\nu(X,X') d_\nu(Y,Y'')\}\Big]\\
&\quad+1-\E d_\nu(X,X')-\E d_\nu(Y,Y')+\E\{ d_\nu(X,X') d_\nu(Y,Y')\}\\
&=\E\{ d_\nu(X,X') d_\nu(Y,Y')\}+\E d_\nu(X,X')\E d_\nu(Y,Y')-2\E\{ d_\nu(X,X') d_\nu(Y,Y'')\}\,.
\end{align*}
Finally, in view of the definition of $d_\nu$ in \eqref{d}, the proof is complete by observing that
\begin{align*}
d_\nu(X,X')+\E d_\nu(X,X')-2\E_Xd_\nu(X,X')&=-\{\phi_{\nu}(X,X')+\E \phi_{\nu}(X,X')-2\E_X\phi_{\nu}(X,X')\}\,,\\
d_\nu(Y,Y')+\E d_\nu(Y,Y')-2\E_Yd_\nu(Y,Y')&=-\{\phi_{\nu}(Y,Y')+\E \phi_{\nu}(Y,Y')-2\E_Y\phi_{\nu}(Y,Y')\}\,.
\end{align*}

\subsection{Proof of Theorem~\ref{thm:iff:dcov}}\label{app:thm:sufficient}

For sufficiency, it is clear that if $X$ and $Y$ are independent, then $\phi_{X,Y}-\phi_X\otimes\phi_Y\equiv0$, so that $\dcov_\nu^2(X,Y)=0$. 
On the other hand, assume $\dcov_\nu^2(X,Y)=0$. Let $\psi(w_1,w_2)=\phi_{X,Y}(w_{1},w_{2})-\phi_X(w_{1})\phi_Y(w_{2})$.
Suppose there exists $(w_{10},w_{20})\in\H^2$ such that $|\psi(w_{10},w_{20})|=a>0$. Observe that
\begin{align*}
&|\psi(w_1,w_2)-\psi(w_{10},w_{20})|\leq|\E\exp(\i\l X,w_1\r+\i\l Y,w_2\r)-\E\exp(\i\l X,w_{10}\r+\i\l Y,w_{20}\r)|\\
&\qquad+|\E\exp(\i\l X,w_1\r)-\E\exp(\i\l X,w_{10}\r)|+|\E\exp(\i\l Y,w_2\r)-\E\exp(\i\l Y,w_{20}\r)|\\
&\leq \E|\exp(\i\l X,w_1-w_{10}\r+\i\l Y,w_2-w_{20}\r)-1|\\
&\qquad+\E|\exp(\i\l X,w_1-w_{10}\r)-1|+\E|\exp(\i\l Y,w_2-w_{20}\r)-1|\\
&\leq 2\E|\exp(\i\l X,w_1-w_{10}\r)-1|+2\E|\exp(\i\l Y,w_2-w_{20}\r)-1|\,.
\end{align*}
By the dominated convergence theorem, it holds that $|\psi(w_1,w_2)-\psi(w_{10},w_{20})|\to0$ as long as $\|w_1-w_{10}\|^2+\|w_2-w_{20}\|^2\to0$. This implies that there exists $r>0$ such that for any $(w_1,w_2)\in\H^2$ such that $\|w_1-w_{10}\|^2+\|w_2-w_{20}\|^2\leq r^2$. It therefore follows that $|\psi(w_1,w_2)-\psi(w_{10},w_{20})|\leq a/2$ such that $|\psi(w_1,w_2)|=|\phi_{X,Y}(w_1,w_2)-\phi_X(w_1)\phi_Y(w_2)|\geq a/2>0$.
Since $\supp(\nu)=\H$, we deduce that (for the product measure $\nu^{\otimes 2}$ on $\H^2$)
\begin{align*}
&\nu^{\otimes 2}\{(w_1,w_2)\in\H^2:\|w_1-w_{10}\|^2+\|w_2-w_{20}\|^2\leq r^2\}\\
&\geq\nu\{w_1\in\H:\|w_1-w_{10}\|^2\leq r^2/2\}\times\nu\{w_2\in\H:\|w_2-w_{20}\|^2\leq r^2/2\}>0\,.
\end{align*}
Therefore, we obtain 
\begin{align*}
\dcov_\nu^2(X,Y)\geq\int_{\H^2}|\psi(w_1,w_2)|\,\one\{\|w_1-w_{10}\|^2+\|w_2-w_{20}\|^2\leq r^2\}\,\nu(\d w_1)\,\nu(\d w_2)>0\,,
\end{align*}
which contradicts the assumption that $\dcov_\nu^2(X,Y)=0$. This implies that $\phi_{X,Y}\equiv\phi_X\otimes\phi_Y$, i.e.~$X$ and $Y$ are independent.

To prove necessity, we follow a similar derivation as the proof of Theorem~\ref{thm:iff:ed}. Suppose $\supp(\nu)\neq\H$, so that there exists a ball $\B(w_0,r)\subset\H$ centered at $w_0\in\H$ with radius $r>0$ such that $\nu\{\B(w_0,r)\}=0$ and $0\notin \B(w_0,r)$. For the Gaussian characteristic functional $\phi_1$ defined in \eqref{phi1}, let
\begin{align*}
\phi_{X,Y}(w_1,w_2)=\phi_1(w_1)\phi_2(w_2)+a^2\theta(w_1)\theta(w_2)\,,
\end{align*}
where $\theta$ and the constant $a$ are as in \eqref{theta} and \eqref{phi2}, respectively. Following the proof of Theorem~\ref{thm:iff:ed}, $\phi_{X,Y}$ defined in the above equation defines a well-defined characteristic functional on $\H^2$. Furthermore, $\phi_{X,Y}(w_1,w_2)-\phi_{X,Y}(w_1,0)\phi_{X,Y}(0,w_2)=a^2\theta(w_1)\theta(w_2)$, which implies that $\phi_{X,Y}\neq\phi_X\otimes\phi_Y$. On the other hand, we have $\dcov_\nu^2(X,Y)=0$. This contradiction implies that $\supp(\nu)=\H$.

\subsection{Proof of Proposition~\ref{prop:dcovn}}\label{app:prop:dcovn}

Observing \eqref{dcovn} and the definition of $\hat\phi_{X,Y}$ in \eqref{hatphixy}, we obtain
Therefore, we obtain from the above equation that
\begin{align}\label{tem}
&|\hat\phi_{X,Y}(w_1,w_2)-\hat\phi_X(w_1)\hat\phi_Y(w_2)|^2=\{\hat\phi_{X,Y}(w_1,w_2)-\hat\phi_X(w_1)\hat\phi_Y(w_2)\}\notag\\
&\qquad\qquad\times\{\hat\phi_{X,Y}(-w_1,-w_2)-\hat\phi_X(-w_1)\hat\phi_Y(-w_2)\}=\frac{1}{n^2}\sum_{k,\ell=1}^nI_{k,\ell}^{(X)}(w_1)I_{k,\ell}^{(Y)}(w_2)\,,
\end{align}
where, for $1\leq k,\ell\leq n$,
\begin{align*}
I_{k,\ell}^{(X)}(w_1)&=\Big\{\exp(\i\l X_{k},w_1\r)-\frac{1}{n}\sum_{k=1}^n\exp(\i\l X_{k},w_1\r)\Big\}\Big\{\exp(-\i\l X_{\ell},w_1\r)-\frac{1}{n}\sum_{{k}=1}^n\exp(-\i\l X_{k},w_1\r)\Big\}\,,
\end{align*}
and $I_{k,\ell}^{(Y,\nu)}$ is defined similarly.
Observe that
\begin{align*}
&\int_\H I_{k,\ell}^{(X)}(w_1)\,\nu(\d w_1)=\int_\H\Big\{\exp(\i\l X_{k}-X_{\ell},w_1\r)-\frac{1}{n}\sum_{j_1=1}^n\exp(\i\l X_{j_1}-X_{\ell},w_1\r)\\
&\hspace{3.5cm}-\frac{1}{n}\sum_{{j_2}=1}^n\exp(\i\l X_{k}-X_{j_2},w_1\r)+\frac{1}{n^2}\sum_{j_1,j_2=1}^n\exp(\i\l X_{j_1}-X_{j_2},w_1\r)\Big\}\,\nu(\d w_1)\\
&=\phi_{\nu}(X_{k}-X_{\ell})-\frac{1}{n}\sum_{j_1=1}^n\phi_{\nu}(X_{j_1}-X_{\ell})-\frac{1}{n}\sum_{j_2=1}^n\phi_{\nu}(X_{k}-X_{j_2})+\frac{1}{n^2}\sum_{j_1,j_2=1}^n\phi_{\nu}(X_{j_1}-X_{j_2})=\hat V^{(X,\nu)}_{k,\ell}\,.
\end{align*}
Similar derivation yields $\int_\H I_{k,\ell}^{(Y)}(w_2)\,\nu(\d w_2)=\hat V^{(Y,\nu)}_{k,\ell}$, which completes the proof.

\subsection{Proof of the conclusion in Remark~\ref{prop:hsic}}\label{app:proof:hsic}

We show that $\hsic_\nu^2(X,Y)=\dcov_\nu^2(X,Y)$. Following the same derivation as in \eqref{mmdp}, we obtain
\begin{align*}
\hsic_\nu^2(X,Y)&=\| \E\{K_\nu(\cdot,X)\otimes K_\nu(\cdot,Y)\}-\E  K_\nu(\cdot,X)\otimes\E K_\nu(\cdot,Y)\|_{\H _{K_\nu}\otimes\H _{K_\nu}}^2\\
&=\E \{K_\nu(X,X')K_\nu(Y,Y')\}+\E K_\nu(X,X')\E K_\nu(Y,Y')-2\E \{\E_X K_\nu(X,X')\E_Y K_\nu(Y,Y')\}\,.\notag
\end{align*} 
In view of the definition of $K_\nu$ in \eqref{knu2}, we obtain from the above equation that
\begin{align*}
\hsic_\nu^2(X,Y)&=\E [\{d_\nu(X,0)+d_\nu(X',0)-d_\nu(X,X')\}\{d_\nu(Y,0)+d_\nu(Y',0)-d_\nu(Y,Y')\}]\\
&+\E\{d_\nu(X,0)+d_\nu(X',0)-d_\nu(X,X')\}\E \{d_\nu(Y,0)+d_\nu(Y',0)-d_\nu(Y,Y')\}\\
&-2\E [\E_X \{d_\nu(X,0)+d_\nu(X',0)-d_\nu(X,X')\}\E_Y\{d_\nu(Y,0)+d_\nu(Y',0)-d_\nu(Y,Y')\}]\,.
\end{align*}
Expanding the above equation and observing the fact that
\begin{align*}
&\E\{d_\nu(X,0)d_\nu(Y,Y')\}=\E\{d_\nu(X',0)d_\nu(Y,Y')\}=\E[d_\nu(X',0)\E_Y\{d_\nu(Y,Y')\}]\,;\\
&\E\{d_\nu(X,X')d_\nu(Y,0)\}=\E\{d_\nu(X,X')d_\nu(Y',0)\}=\E[\E_X\{d_\nu(X,X')\}d_\nu(Y',0)]\,,
\end{align*}
we obtain $\hsic_\nu^2(X,Y)=\E\{ d_\nu(X,X') d_\nu(Y,Y')\}+\E d_\nu(X,X')\,\E d_\nu(Y,Y')-2\E\{ d_\nu(X,X') d_\nu(Y,Y'')\}$. The proof is complete in view of Proposition~\ref{thm:dcov}.

\subsection{Proof of Proposition~\ref{thm:dvar}}\label{app:thm:dvar}

For (i), observe that
\begin{align}\label{dvar1}
&\dvar_\nu^2(X)=\int_{\H^2}\big|\E \exp(\i\l X,w_1+w_2\r)-\E \exp(\i\l X,w_1\r)\E\exp(\i\l X,w_2\r)\big|^2\,\nu(\d w_1)\,\nu(\d w_2)\notag\\
&\hspace{1.3cm}=\int_{\H^2}\{I_1(w_1,w_2)+I_2(w_1,w_2)-I_3(w_1,w_2)-I_4(w_1,w_2)\}\nu(\d w_1)\nu(\d w_2)\,,
\end{align}
where
\begin{align*}
&I_1(w_1,w_2)=|\E \exp(\i\l X,w_1+w_2\r)|^2\,,\\
&I_2(w_1,w_2)=|\E \exp(\i\l X,w_1\r)|^2|\E\exp(\i\l X,w_2\r)|^2\,,\\
&I_3(w_1,w_2)=\E \exp(\i\l X,w_1+w_2\r)\E \exp(-\i\l X,w_1\r)\E\exp(-\i\l X,w_2\r)\,,\\
&I_4(w_1,w_2)=\E \exp(-\i\l X,w_1+w_2\r)\E \exp(\i\l X,w_1\r)\E\exp(\i\l X,w_2\r)\,.
\end{align*}
Recall that $(X',Y')$ and $(X'',Y'')$ are i.i.d.~copies of $(X,Y)$.  Note that for the first term $I_1$, we have $I_1(w_1,w_2)=\E \exp(\i\l X-X',w_1+w_2\r)$.
Therefore, by Fubini's theorem we obtain
\begin{align}\label{dvar2}
&\hspace{-1em}\int_{\H^2}I_1(w_1,w_2)\nu(\d w_1)\nu(\d w_2)=\E\Big\{\int_{\H^2} \exp(\i\l X-X',w_1+w_2\r)\,\nu(\d w_1)\,\nu(\d w_2)\Big\}\notag\\
&\hspace{-1em}=\E\Big\{\int_{\H }\exp(\i\l X-X',w_1\r)\nu(\d w_1)\int_{\H } \exp(\i\l X-X',w_2\r)\nu(\d w_2)\Big\}=\E\{{\phi_{\nu}}(X-X')\}^2\,.
\end{align}
For the second term $I_2$, note that $I_2(w_1,w_2)=\E \exp(\i\l X-X',w_1\r)\E\exp(\i\l X-X',w_2\r)$.
Therefore, we obtain from the above equation that
\begin{align}\label{dvar3}
&\int_{\H^2}I_2(w_1,w_2)\nu(\d w_1)\nu(\d w_2)=\int_{\H^2}\E \exp(\i\l X-X',w_1\r)\E\exp(\i\l X-X',w_2\r)\nu(\d w_1)\nu(\d w_2)\notag\\
&=\Big\{\int_{\H }\E \exp(\i\l X-X',w_1\r)\,\nu(\d w_1)\Big\}^2=\{\E\phi_{\nu}(X-X')\}^2\,.
\end{align}
For the third term $I_3$, we have $I_3(w_1,w_2)=\E \exp(\i\l X-X',w_1\r+\i\l X-X'',w_2\r)$. Therefore, we have
\begin{align}\label{dvar4}
&\int_{\H^2}I_3(w_1,w_2)\nu(\d w_1)\nu(\d w_2)=\E\Big\{\int_{\H^2}\exp(\i\l X-X',w_1\r+\i\l X-X'',w_2\r)\nu(\d w_1)\nu(\d w_2)\Big\}\notag\\
&=\E\{\phi_{\nu}(X-X')\phi_{\nu}(X-X'')\}\,.
\end{align}
For the fourth term $I_4$, we have $I_4(w_1,w_2)=\E \exp(-\i\l X-X',w_1\r-\i\l X-X'',w_2\r)$.
Hence, since $\phi_{\nu}$ is symmetric due to Lemma~\eqref{lem:sym}, we have
\begin{align}\label{dvar5}
&\int_{\H^2}I_4(w_1,w_2)\nu(\d w_1)\nu(\d w_2)=\E\Big\{\int_{\H^2}\exp(-\i\l X-X',w_1\r-\i\l X-X'',w_2\r)\nu(\d w_1)\nu(\d w_2)\Big\}\notag\\
&=\E\{\phi_{\nu}(X'-X)\phi_{\nu}(X''-X)\}=\E\{\phi_{\nu}(X-X')\phi_{\nu}(X-X'')\}\,.
\end{align}
Therefore, (i) follows by combining \eqref{dvar1}--\eqref{dvar5} and the fact that $d_\nu(x,x')=1-\phi_{\nu}(x-x')$.

Next, for (ii), we have
\begin{align*}
&|\phi_{X,X}(w_1,w_2)-\phi_X(w_1)\phi_X(w_2)|^2=\big|\E \exp(\i\l X,w_1+w_2\r)-\E \exp(\i\l X,w_1\r)\E\exp(\i\l X,w_2\r)\big|^2\\
&=\big|\E\big[\{\exp(\i\l X,w_1\r)-\E\exp(\i\l X,w_1\r)\}\{\exp(\i\l X,w_2\r)-\E\exp(\i\l X,w_2\r)\}\big]\big|^2\\
%
%
&=\E\big[\{\exp(\i\l X,w_1\r)-\E\exp(\i\l X,w_1\r)\}\{\exp(-\i\l X',w_1\r)-\E\exp(-\i\l X',w_1\r)\}\\
&\quad\times\{\exp(\i\l X,w_2\r)-\E\exp(\i\l X,w_2\r)\}\{\exp(-\i\l X',w_2\r)-\E\exp(-\i\l X',w_2\r)\}\big]
\end{align*}
Therefore, by Fubini's theorem, we find
\begin{align*}
&\dvar_\nu^2(X)=\E\bigg[\int_{\H }\{\exp(\i\l X,w_1\r)-\E\exp(\i\l X,w_1\r)\}\{\exp(-\i\l X',w_1\r)-\E\exp(-\i\l X',w_1\r)\}\nu(\d w_1)\bigg]^2\,.
\end{align*}
The proof of (ii) is complete in view of \eqref{dc1} and the fact that $d_\nu(x,x')=1-\phi_{\nu}(x-x')$.

For (iii), note that $\dvar_\nu^2(X)=0$ if and only if $\phi_X(w_1+w_2)=\phi_{X,X}(w_1,w_2)=\phi_X(w_1)\phi_X(w_2)$ for any $w_1,w_2\in\H$, which is equivalent to $\log\{\phi_X(w)\}$ being a (constant) linear functional of $w$. This is further equivalent to $X$ being degenerate.


\subsection{Proof of Theorem~\ref{thm:prop_dcor}}\label{app:thm:prop_dcor}

For (i), by Propositions~\ref{thm:dcov}, \ref{thm:dvar} and Cauchy-Schwarz inequality,
\begin{align}\label{cs}
&\dcov_\nu^2 (X,Y)\notag\\
&=\E\big[\{\phi_{\nu}(X-X')+\E \phi_{\nu}(X-X')-2\E_X\phi_{\nu}(X-X')\}\{\phi_{\nu}(Y-Y')+\E \phi_{\nu}(Y-Y')-2\E_Y\phi_{\nu}(Y-Y')\}\big]\notag\\
&\leq\big(\E\big[\{\phi_{\nu}(X-X')+\E \phi_{\nu}(X-X')-2\E_X\phi_{\nu}(X-X')\}^2\big]\big)^{1/2}\notag\\
&\quad\times
\big(\E\big[\{\phi_{\nu}(Y-Y')+\E \phi_{\nu}(Y-Y')-2\E_Y\phi_{\nu}(Y-Y')\}^2\big]\big)^{1/2}=\dvar_\nu(X)\times\dvar_\nu(Y)\,.
\end{align}

For (ii), suppose $\dcor_\nu^2(X)=0$. If $\dvar_\nu(X)\dvar_\nu(Y)\neq0$, we have $\dcov_\nu(X,Y)=0$, which by Theorem~\ref{thm:iff:dcov} is equivalent to $X$ and $Y$ being independent. If $\dvar_\nu(X)\dvar_\nu(Y)=0$, then either $X$ or $Y$ degenerates, which is equivalent to $X$ and $Y$ being independent.

For (iii), note that \eqref{cs} holds if and only if there exists a constant $a\in\mathbb R$ such that
\begin{align*}
\phi_{\nu}(X-X')+\E \phi_{\nu}(X-X')-2\E_X\phi_{\nu}(X-X')=
a\{\phi_{\nu}(Y-Y')+\E \phi_{\nu}(Y-Y')-2\E_Y\phi_{\nu}(Y-Y')\}\,.
\end{align*}
This is further equivalent to, for $(X,Y)$-almost every $(x,y),(x',y')\in\H^2$, and $c\in\mathbb R$,
\begin{align*}
&\phi_{\nu}(x-x')+\E \phi_{\nu}(X-X')-2\E\{\phi_{\nu}(X-X')|X'=x'\}\\
&=
c\big[\phi_{\nu}(y-y')+\E \phi_{\nu}(Y-Y')-a\E\{\phi_{\nu}(Y-Y')|Y'=y'\}\big]\,.
\end{align*}
By taking $x=x'$ and $y=y'$ in the above equation, we obtain
\begin{align*}
\E \phi_{\nu}(X-X')-2\E\{\phi_{\nu}(X-X')|X'=x'\}=
a\big[\E \phi_{\nu}(Y-Y')-2\E\{\phi_{\nu}(Y-Y')|Y'=y'\}\big]\,,
\end{align*}
which further implies that $\phi_{\nu}(x-x')=a \phi_{\nu}(y-y')$, for $(X,Y)$-almost every $(x,y)$ and $(x',y')$.


\subsection{Proof of Proposition~\ref{thm:converge}}\label{app:thm:converge}

We start with the following lemma, which is proved in Section~\ref{proof:lem:a.5}.

\begin{lemma}\label{lem:a.5}
Suppose Assumption~\ref{a:sym} is satisfied. Then, we have
\begin{align}\label{v}
\hat\dcov_\nu^2(X,Y)=\frac{1}{n^4}\sum_{j_1,j_2,j_3,j_4=1}^nh\{(X_{j_1},Y_{j_1}),(X_{j_2},Y_{j_2}),(X_{j_3},Y_{j_3}),(X_{j_4},Y_{j_4})\}\,,
\end{align}
where, for $(x_1,y_1),(x_2,y_2),(x_3,y_3),(x_4,y_4)\in\H^2$,
\begin{align}\label{h}
&h\{(x_1,y_1),(x_2,y_2),(x_3,y_3),(x_4,y_4)\}=\frac{1}{4!}\sum_{(k_1,k_2,k_3,k_4)\in I_4(4)}\big\{\phi_{\nu}(x_{k_1}-x_{k_2})\phi_{\nu}(y_{k_3}-y_{k_4})\notag\\
&\hspace{3cm}+\phi_{\nu}(x_{k_1}-x_{k_2})\phi_{\nu}(y_{k_1}-y_{k_2})-2\phi_{\nu}(x_{k_1}-x_{k_2})\phi_{\nu}(y_{k_1}-y_{k_3})\big\}\,.
\end{align}
\end{lemma}

For $(x_1,y_1),(x_2,y_2)\in\H^2$, define
\begin{align}\label{h2}
&h_1(x,y)=\E\big[h\big\{(x,y),(X^{(1)},Y^{(1)}),(X^{(2)},Y^{(2)}),(X^{(3)},Y^{(3)})\big\}\big]\,;\\
&h_2\{(x_1,y_1),(x_2,y_2)\}=\E\big[h\big\{(x_1,y_1),(x_2,y_2),(X^{(1)},Y^{(1)}),(X^{(2)},Y^{(2)})\big\}\big]-h_1(x_1,y_1)-h_1(x_2,y_2)\,,\notag
\end{align}
where $(X^{(1)},Y^{(1)}),(X^{(2)},Y^{(2)}),(X^{(3)},Y^{(3)})$ that are i.i.d.~copies of $(X,Y)$.
We will use the following lemma which shows that $\hat\dcov_\nu^2(X,Y)$ in \eqref{v} is a V-statistic of rank-1 in the case of independence, which is proved in Section~\ref{app:lem:dcovv} of the supplementary material.

\begin{lemma}\label{lem:dcovv}
Suppose Assumptions~\ref{a:sym} and \ref{a:supp} are satisfied. When $X\indep Y$, it holds that $h_1\equiv0$ and, for $(x_1,y_1),(x_2,y_2)\in\H^2$,  $h_2\{(x_1,y_1),(x_2,y_2)\}=6^{-1}g_X(x_1,x_2)\,g_Y(y_1,y_2)$, where
\begin{align}\label{gxy}
g_X(x_1,x_2)&=\phi_{\nu}(x_1-x_2)+\E \phi_{\nu}(X-X')-\E \phi_{\nu}(x_1-X)-\E \phi_{\nu}(X-x_2)\,,\notag\\
g_Y(y_1,y_2)&=\phi_{\nu}(y_1-y_2)+\E \phi_{\nu}(Y-Y')-\E \phi_{\nu}(y_1-Y)-\E \phi_{\nu}(Y-y_2)\,.
\end{align}
Furthermore, $\var\{h_1(X,Y)\}=0<\var[h_2\{(X^{(1)},Y^{(1)}),(X^{(2)},Y^{(2)})\}]$.

\end{lemma}

Now, for (i) in Proposition~\ref{thm:converge}, observe that the kernel $h$ in \eqref{h} is such that $\|h\|_\infty\leq 4$. The strong consistency in (i) therefore follows  from Lemma~\ref{lem:a.5} and the strong law of large numbers of V-statistics (see, for example, Theorem~5.2.9 in \cite{pena2012}).

For (ii), in view of \eqref{h2} and the fact that $h_1\equiv0$ in Lemma~\ref{lem:dcovv}, let
\begin{align*}
&h_3\{(x_1,y_1),(x_2,y_2),(x_3,y_3)\}=\E\big[h\big\{(x_1,y_1),(x_2,y_2),(x_3,y_3),(X,Y)\big\}\big]-\sum_{1\leq j<k\leq3}h_2\{(x_j,y_j),(x_k,y_k)\}\,,\\
&h_4\{(x_1,y_1),(x_2,y_2),(x_3,y_3),(x_4,y_4)\}=h\big\{(x_1,y_1),(x_2,y_2),(x_3,y_3),(x_4,y_4)\big\}\\
&\hspace{4cm}-\sum_{1\leq j<k\leq4}h_2\{(x_j,y_j),(x_k,y_k)\}-\sum_{1\leq j<k<\ell\leq4}h_3\{(x_j,y_j),(x_k,y_k),(x_\ell,y_\ell)\}\,.
\end{align*}
Recall from Lemma~\ref{lem:dcovv} that $h_1\equiv0$. Then, by Hoeffding's decomposition (see, for example, Section~1.3 in \cite{boroskikh1996}), we may decompose the V-statistic $\hat\dcov_\nu^2(X,Y)$ in \eqref{v} as the sum of canonical V-statistics as $\hat\dcov_\nu^2(X,Y)=V_{n,2}+V_{n,3}+V_{n,4}$, where
\begin{align}\label{v234}
&V_{n,2}=\frac{6}{n^2}\sum_{j_1,j_2=1}^n h_2\{(X_{j_1},Y_{j_1}),(X_{j_2},Y_{j_2})\}\,,\notag\\
&V_{n,3}=\frac{6}{n^3}\sum_{j_1,j_2,j_3=1}^n h_3\{(X_{j_1},Y_{j_1}),(X_{j_2},Y_{j_2}),(X_{j_3},Y_{j_3})\}\,,\notag\\
&V_{n,4}=\frac{1}{n^4}\sum_{j_1,j_2,j_3,j_4=1}^n h_4\{(X_{j_1},Y_{j_1}),(X_{j_2},Y_{j_2}),(X_{j_3},Y_{j_3}),(X_{j_4},Y_{j_4})\}\,.
\end{align}
We first show that, as $n\to\infty$,
\begin{align}\label{t2}
nV_{n,2}\converged \E \phi_{\nu}(X-X')\,\E d_{\nu}(Y-Y')+\sum_{j=1}^\infty\lambda_j(\zeta_j^2-1)\,.
\end{align}  
To see this, in view of Lemma~\ref{lem:dcovv}, we write $V_{n,2}=n^{-1}U_{n,2,1}+n^{-1}(n-1)U_{n,2,2}$, where
\begin{align*}
U_{n,2,1}&=\frac{1}{n}\sum_{j=1}^n g_X(X_j,X_j)\,g_Y(Y_j,Y_j)\,,\qquad U_{n,2,2}=\frac{1}{n(n-1)}\sum_{(j_1,j_2)\in I_2(n)} g_X(X_{j_1},X_{j_2})\,g_Y(Y_{j_1},Y_{j_2})\,.
\end{align*}
For $U_{n,2,1}$, note that in view of the definition of $g_X$ in \eqref{gxy},
\begin{align*}
\E g_X(X,X)=\E \phi_{\nu}(X-X')+\E \phi_{\nu}(X-X)-2\E \{\E \phi_{\nu}(X-X')|X'\}=-\E \phi_{\nu}(X-X')\,.
\end{align*}
By the law of large numbers, we obtain $U_{n,2,1}\converged \E \phi_{\nu}(X-X')\,\E \phi_{\nu}(X-X')$.
For $U_{n,2,2}$, Note that $g_X(x,X)=\phi_{\nu}(x-X)-\E \phi_{\nu}(x-X)$. It follows from Theorem~5.5.2 in \cite{serfling1980} that $nU_{n,2,2}\converged \sum_{j=1}^\infty\lambda_j(\zeta_j^2-1)$.

Next, for $V_{n,3}$ and $V_{n,4}$, we aim to show that $n(V_{n,2}+V_{n,3})=o_p(1)$. To achieve this, we further decompose these V-statistics as the linear combination of U-statistics. Following Section~1.3 in \cite{boroskikh1996}, we have
\begin{align*}
nV_{n,3}&=nU_{n,3,1}+nU_{n,3,2}+nU_{n,3,3}\,,\hspace{1cm} nV_{n,4}=nU_{n,4,1}+nU_{n,4,2}+nU_{n,4,3}+nU_{n,4,4}\,,
\end{align*}
where the rescaled U-statistics are given by
\begin{align*}
nU_{n,3,1}&=\frac{1}{n^2}\sum_{j=1}^n h_3\{(X_j,Y_j),(X_j,Y_j),(X_j,Y_j)\}\\
nU_{n,3,2}&=\frac{3}{n^2}\sum_{(j_1,j_2)\in I_2(n)}h_3\{(X_{j_1},Y_{j_1}),(X_{j_1},Y_{j_1}),(X_{j_2},Y_{j_2})\}\\
nU_{n,3,3}&=\frac{1}{n^2}\sum_{(j_1,j_2,j_3)\in I_3(n)}h_3\{(X_{j_1},Y_{j_1}),(X_{j_2},Y_{j_2}),(X_{j_3},Y_{j_3})\}\\
nU_{n,4,1}&=\frac{1}{n^3}\sum_{j=1}^n h_4\{(X_j,Y_j),(X_j,Y_j),(X_j,Y_j),(X_j,Y_j)\}\\
nU_{n,4,2}&=\frac{1}{n^3}\sum_{(j_1,j_2)\in I_2(n)}\big[4h_4\{(X_{j_1},Y_{j_1}),(X_{j_1},Y_{j_1}),(X_{j_1},Y_{j_1}),(X_{j_2},Y_{j_2})\}\\
&\hspace{3cm}+3h_4\{(X_{j_1},Y_{j_1}),(X_{j_1},Y_{j_1}),(X_{j_2},Y_{j_2}),(X_{j_2},Y_{j_2})\}\big]\\
nU_{n,4,3}&=\frac{6}{n^3}\sum_{(j_1,j_2,j_3)\in I_3(n)}h_4\{(X_{j_1},Y_{j_1}),(X_{j_1},Y_{j_1}),(X_{j_2},Y_{j_2}),(X_{j_3},Y_{j_3})\}\\
nU_{n,4,4}&=\frac{1}{n^3}\sum_{(j_1,j_2,j_3,j_4)\in I_4(n)}h_4\{(X_{j_1},Y_{j_1}),(X_{j_2},Y_{j_2}),(X_{j_3},Y_{j_3}),(X_{j_4},Y_{j_4})\}\,.
\end{align*}
For $U_{n,4,4}$ and $U_{n,3,3}$, it follows from Theorem~1.6.2 in \cite{lee1990} p.~28 that the symmetric kernels $h_3$ and $h_4$ are completely degenerate.
By the strong law of large numbers for degenerate U-statistics (see, for example, Theorem~2 in \cite{gine1992}), as $n\to\infty$,
\begin{align*}
&U_{n,3,3}=\frac{1}{n^2}\sum_{(j_1,j_2,j_3)\in I_3(n)}\big[h_3\{(X_{j_1},Y_{j_1}),\ldots,(X_{j_3},Y_{j_3})\}-\E h_3\{(X_{j_1},Y_{j_1}),\ldots,(X_{j_3},Y_{j_3})\}\big]=o_p(1)\,,\\
&U_{n,4,4}=\frac{1}{n^3}\sum_{(j_1,j_2,j_3,j_4)\in I_4(n)}\big[h_4\{(X_{j_1},Y_{j_1}),\ldots,(X_{j_4},Y_{j_4})\}-\E h_4\{(X_{j_1},Y_{j_1}),\ldots,(X_{j_4},Y_{j_4})\}\big]=o_p(1)\,.
\end{align*}
For $U_{n,3,2}$ and $U_{n,4,3}$, note that
\begin{align*}
&\E(U_{n,3,2})=6\E[h_3\{(X_{1},Y_{1}),(X_{1},Y_{1}),(X_{2},Y_{2})\}]=6\E[\tilde h_{3}\{(X_1,Y_1),(X_1,Y_1)\}]\,,\\
&\E(U_{n,4,3})=36\E[h_4\{(X_{1},Y_{1}),(X_{1},Y_{1}),(X_{2},Y_{2}),(X_{3},Y_{3})\}]=36\E[\tilde h_{4}\{(X_1,Y_1),(X_1,Y_1)\}]\,,
\end{align*}
where 
\begin{align*}
&\tilde h_{3}\{(x_1,y_1),(x_2,y_2)\}=\E[h_3\{(x_1,y_1),(x_2,y_2),(X,Y)\}]\,,\\
&\tilde h_{4}\{(x_1,y_1),(x_2,y_2)\}=\E[h_4\{(x_{1},y_{1}),(x_{2},y_{2}),(X_{2},Y_{2}),(X_{3},Y_{3})\}]\,.
\end{align*}
Since $h_3$ and $h_4$ are completely degenerate, it follows that $\tilde h_{3}\equiv0$ and $\tilde h_{4}\equiv0$, so that by the strong law of large numbers for U-statistics, we deduce that $nU_{n,3,2}=o_p(1)$ and $nU_{n,4,3}=o_p(1)$ as $n\to\infty$.
For $U_{n,3,1}$, $U_{n,4,1}$ and $U_{n,4,2}$, since $\|h_3\|_\infty+\|h_4\|_\infty<\infty$, direct bounds yields that $U_{n,3,1}+U_{n,4,1}+U_{n,4,2}=o_p(1)$ as $n\to\infty$. In conclusion, we have shown $n(V_{n,2}+V_{n,3})=o_p(1)$, combining which with \eqref{t2} and Slutsky's lemma concludes the proof.

\subsection{Proof of Theorem~\ref{thm:upper:ind}}\label{app:thm:upper:ind}

We start by introducing some useful notations.
Let $\Pi$ denote a generic random permutation of $(1,2,\ldots,n )$ 
and denote  by
\begin{align*}
\hat F_{\nu,B}(q)=\frac{1}{B}\sum_{b=1}^B\one\{\hat\dcov_\nu^2(X,Y,\Pi_b)\leq q\}\,
\end{align*}
the empirical cumulative distribution function of the permuted samples.
Define
\begin{align}\label{hatqnu}
\hat q_{\nu,1-\alpha}=\inf\big\{q\in\mathbb R:\hat F_{\nu,B}(q)\geq 1-\alpha\big\}\,.
\end{align}
Define  
\begin{align}\label{tildeqnu}
\tilde F_{\nu}(q)&=|I_n(n)|^{-1}\sum_{\Pi_0\in I_n(n)}\one\{\hat D_\nu(\T_n^{(\Pi_0)})\leq q\}\,,\notag
\\
\tilde q_{\nu,1-\alpha}&=\inf\big\{q\in\mathbb R:\tilde F_{\nu}(q)\geq 1-\alpha\big\}\,.
\end{align}
Write $\|\phi_{X,Y}-\phi_X\otimes\phi_Y\|_\nu^2=\int_{\H^2}|\phi_{X,Y}-\phi_X\otimes\phi_Y|^2\d(\nu\times\nu)$. In the sequel, suppose the joint distribution $F$ of $(X,Y)$ satisfies $F\in\mathcal F_\nu(\rho)$.

By the Dvoretzky–Kiefer–Wolfowitz inequality (\cite{massart1990}), for any $\e>0$, $\nu\in\V$,
\begin{align*}
\P_F\bigg\{\sup_{q\in\mathbb R}|\hat F_{\nu,B}(q)-\tilde F_{\nu}(q)|>\e\bigg\}\leq 2\exp(-2B\e^2)\,.
\end{align*}
Define the event $\mathcal E_\nu=\big\{\sup_{q\in\mathbb R}|\hat F_{\nu,B}(q)-\tilde F_{\nu}(q)|\leq\sqrt{(2B)^{-1}\log(6/\beta)}\big\}$.
By taking $\beta/3$ on the right-hand side of the above equation, we obtain that
\begin{align}\label{p5.5}
\sup_{\nu\in\V}\sup_{F\in\mathcal F_\nu(\rho)}\P_F(\mathcal E_\nu^c)=\sup_{\nu\in\V}\sup_{F\in\mathcal F_\nu(\rho)}\P_F\bigg\{\sup_{q\in\mathbb R}|\hat F_{\nu,B}(q)-\tilde F_{\nu}(q)|>\sqrt{(2B)^{-1}\log(6/\beta)}\bigg\}\leq\beta/3\,.
\end{align}
Recalling the definition of $\hat q_{\nu,1-\alpha}$ in \eqref{hatqnu}, we obtain from the above equation that, on event $\mathcal E_\nu$,
\begin{align*}
\hat q_{\nu,1-\alpha}&\leq \inf\big\{q\in\mathbb R:\tilde F_{\nu}(q)\geq 1-\alpha+\sqrt{(2B)^{-1}\log(6/\beta)}\big\}\leq \inf\big\{q\in\mathbb R:\tilde F_{\nu}(q)\geq1-\alpha/2\big\}=\tilde q_{\nu,1-\alpha/2}\,,
\end{align*}
where we used the assumption that $B\geq2\alpha^{-2}\log(6/\beta)$ and the definition of $\tilde q_{\nu,1-\alpha}$ in \eqref{tildeqnu}.
Combining the above equation with \eqref{p5.5} enables us to deduce that, the type-II error of the permutation test $\hat\psi_{\nu,n}(\alpha)$ in \eqref{hatpsi} satisfies
\begin{align}\label{p5.6}
\sup_{\nu\in\V}\sup_{F\in\mathcal F_\nu(\rho)}\P_F\{\hat\psi_{\nu,n}(\alpha)=0\}&=\sup_{\nu\in\V}\sup_{F\in\mathcal F_\nu(\rho)}\P_F\{\hat\dcov_\nu^2(X,Y)\leq \hat q_{\nu,1-\alpha}\}\notag\\
&\leq\sup_{\nu\in\V}\sup_{F\in\mathcal F_\nu(\rho)}\big[\P_F\{\hat\dcov_\nu^2(X,Y)\leq \hat q_{\nu,1-\alpha};\mathcal E_\nu\}+\P_F(\mathcal E_\nu^c)\big]\notag\\
&\leq\sup_{\nu\in\V}\sup_{F\in\mathcal F_\nu(\rho)}\P_F\{\hat\dcov_\nu^2(X,Y)\leq \tilde q_{\nu,1-\alpha/2}\}+\beta/3\,.
\end{align}

Next, we compute the probability that $\hat\dcov_\nu^2(X,Y)\leq \tilde q_{\nu,1-\alpha/2}$. Write $\|\phi_{X,Y}-\phi_X\otimes\phi_Y\|_\nu^2=\int_{\H^2}|\phi_{X,Y}-\phi_X\otimes\phi_Y|^2\d \nu$.
Let $\tilde w_{\nu,\gamma}$ denote the $\gamma$-quantile of $\tilde q_{\nu,1-\alpha/2}$; here for simplicity we notationally remove the dependency of $\tilde w_{\nu,\gamma}$ on $\alpha$. Then, by the definition of quantile, we deduce that
\begin{align}\label{p0}
\P_F\{\hat\dcov_\nu^2(X,Y)\leq \tilde q_{\nu,1-\alpha/2}\} &
\leq \P_F(\tilde q_{\nu,1-\alpha/2}> \tilde w_{\nu,1-\beta/3})+\P_F\{\hat\dcov_\nu^2(X,Y)\leq \tilde w_{\nu,1-\beta/3}\}\notag\\
&=\beta/3+\P_F\{\hat\dcov_\nu^2(X,Y)\leq \tilde w_{\nu,1-\beta/3}\}\,.
\end{align}
Furthermore, let $\kappa_{\nu,1-\beta/3}$ denote the $(1-\beta/3)$ quantile of $|\hat\dcov_\nu^2(X,Y)-\E_F\{\hat\dcov_\nu^2(X,Y)\}|$. Then, we obtain
\begin{align}\label{p1}
&\sup_{\nu\in\V}\sup_{F\in\mathcal F_\nu(\rho)}\P_F\big[\hat\dcov_\nu^2(X,Y)\leq \E_F\{\hat\dcov_\nu^2(X,Y)\}-\kappa_{\nu,1-\beta/3}\big]\notag\\
&\leq\sup_{\nu\in\V}\sup_{F\in\mathcal F_\nu(\rho)}\P_F\big[|\hat\dcov_\nu^2(X,Y)-\E_F\{\hat\dcov_\nu^2(X,Y)\}|\geq\kappa_{\nu,1-\beta/3}\big]\leq\beta/3\,.
\end{align}
Suppose
\begin{align*}
\|\phi_{X,Y}-\phi_X\otimes\phi_Y\|_\nu^2&\geq \sup_{\nu\in\V}\sup_{F\in\mathcal F_\nu(\rho)}\tilde w_{\nu,1-\beta/3}+\sup_{\nu\in\V}\sup_{F\in\mathcal F_\nu(\rho)}\kappa_{\nu,1-\beta/3}\\
&+\sup_{\nu\in\V}\sup_{F\in\mathcal F_\nu(\rho)}|\|\phi_{X,Y}-\phi_X\otimes\phi_Y\|_\nu^2-\E_F\{\hat\dcov_\nu^2(X,Y)\}|\,.
\end{align*}
Then, we deduce from \eqref{p5.6}--\eqref{p1} that 
\begin{align}\label{p01}
&\sup_{\nu\in\V}\sup_{F\in\mathcal F_\nu}\P_F\{\hat\psi_{\nu,n}(\alpha)=0\}\leq2\beta/3+\sup_{\nu\in\V}\sup_{F\in\mathcal F_\nu}\P_F\{\hat\dcov_\nu^2(X,Y)\leq \tilde w_{\nu,1-\beta/3}\}\notag\\
&\leq2\beta/3+\sup_{\nu\in\V}\sup_{F\in\mathcal F_\nu}\P_F\Big[\hat\dcov_\nu^2(X,Y)\leq \|\phi_{X,Y}-\phi_X\otimes\phi_Y\|_\nu^2-\kappa_{\nu,1-\beta/3}\notag\\
&\qquad-|\|\phi_{X,Y}-\phi_X\otimes\phi_Y\|_\nu^2-\E_F\{\hat\dcov_\nu^2(X,Y)\}|\Big]\notag\\
&\leq2\beta/3+\sup_{\nu\in\V} \sup_{F\in\mathcal F_\nu}\P_F\big[\hat\dcov_\nu^2(X,Y)\leq \E_F\{\hat\dcov_\nu^2(X,Y)\}-\kappa_{\nu,1-\beta/3}\big]\leq\beta\,.
\end{align}
Therefore, recalling the definition of the uniform separation rate $\rho$ defined in \eqref{seprate}, we deduce from the above equation that,
\begin{align}\label{qbeta}
\rho\{\hat\psi_{\nu,n}(\alpha),\beta\}&\leq \sup_{\nu\in\V}\sup_{F\in\mathcal F_\nu(\rho)}\tilde w_{\nu,1-\beta/3}+\sup_{\nu\in\V}\sup_{F\in\mathcal F_\nu(\rho)}\kappa_{\nu,1-\beta/3}\notag\\
&\quad+\sup_{\nu\in\V}\sup_{F\in\mathcal F_\nu(\rho)}|\|\phi_{X,Y}-\phi_X\otimes\phi_Y\|_\nu^2-\E_F\{\hat\dcov_\nu^2(X,Y)\}|\,.
\end{align}
Therefore, in order to show the upper bound for $\rho\{\hat\psi_{\nu,n}(\alpha),\beta\}$, it suffices obtain union bounds for the three terms $\tilde w_{\nu,1-\beta/3}$, $\kappa_{\nu,1-\beta/3}$ and $|\|\phi_{X,Y}-\phi_X\otimes\phi_Y\|_\nu^2-\E_F\{\hat\dcov_\nu^2(X,Y)\}|$ in \eqref{qbeta}.

\textbf{Step 1.} Bounding the quantile $\tilde w_{\nu,1-\beta/3}$.

To obtain an upper bound for $\tilde w_{\nu,1-\beta/3}$, we follow the decoupling idea of \cite{kim2022} and derive a concentration inequality for $\tilde q_{\nu,1-\alpha/2}$.
Observe from \eqref{hatv} that
\begin{align*}
\hat V_{k,\ell}^{(X,\nu)}=\frac{1}{n^2}\sum_{i,j=1}^n\Big\{\phi_{\nu}(X_k-X_\ell)-\phi_{\nu}(X_{i}-X_\ell)-\phi_{\nu}(X_k-X_{j})+\phi_{\nu}(X_{i}-X_{j})\Big\}\,,
\end{align*}
and similar formula is valid for $\hat V_{k,\ell}^{(Y,\nu)}$. Hence, in view of the definition of $\hat\dcov_\nu^2(X,Y,\Pi)$ in \eqref{hatdb}, we obtain
\begin{align}\label{p10}
&\hat\dcov_\nu^2(X,Y,\Pi)=\frac{1}{n^2}\sum_{k,\ell=1}^n\bigg[\frac{1}{n^2}\sum_{i,j=1}^n\big\{\phi_{\nu}(X_k-X_\ell)-\phi_{\nu}(X_{i}-X_\ell)-\phi_{\nu}(X_k-X_{j})+\phi_{\nu}(X_{i}-X_{j})\big\}\notag\\
&\times\frac{1}{n^2}\sum_{i',j'=1}^n\big\{\phi_{\nu}(Y_{\Pi(k)}-Y_{\Pi(\ell)})-\phi_{\nu}(Y_{\Pi(i')}-Y_{\Pi(\ell)})-\phi_{\nu}(Y_{\Pi(k)}-Y_{\Pi(j')})+\phi_{\nu}(Y_{\Pi(i')}-Y_{\Pi(j')})\big\}\bigg]\notag\\
&=\frac{1}{n^6}\sum_{j_1,j_2,j_3,j_4,j_5,j_6=1}^n\big\{\phi_{\nu}(X_{j_1}-X_{j_2})-\phi_{\nu}(X_{j_3}-X_{j_2})-\phi_{\nu}(X_{j_1}-X_{j_4})+\phi_{\nu}(X_{j_3}-X_{j_4})\big\}\notag\\
&\times\big\{\phi_{\nu}(Y_{\Pi(j_3)}-Y_{\Pi(j_4)})-\phi_{\nu}(Y_{\Pi(j_5)}-Y_{\Pi(j_4)})-\phi_{\nu}(Y_{\Pi(j_3)}-Y_{\Pi(j_6)})+\phi_{\nu}(Y_{\Pi(j_5)}-Y_{\Pi(j_6)})\big\}.
\end{align}
For $y_1,y_2,y_3,y_4\in\H$, define
\begin{align}\label{omeganu}
\omega_\nu(y_1,y_2,y_3,y_4)=\phi_{\nu}(y_1-y_2)-\phi_{\nu}(y_3-y_2)-\phi_{\nu}(y_1-y_4)+\phi_{\nu}(y_3-y_4)\,.
\end{align}
Then, in view of \eqref{p10} and \eqref{omeganu}, we obtain
\begin{align*}
\hat\dcov_\nu^2(X,Y,\Pi)&=\frac{1}{n^6}\sum_{j_1,j_2,j_3,j_4,j_5,j_6=1}^n\omega_\nu(X_{j_1},X_{j_2},X_{j_3},X_{j_4})\,
\omega_\nu(Y_{\Pi(j_3)},Y_{\Pi(j_4)},Y_{\Pi(j_5)},Y_{\Pi(j_6)})\,.
\end{align*}

Suppose $\K_1,\ldots,\K_{2\lfloor n/2\rfloor}
$ are i.i.d.~random variables and uniformly distributed over $\{1,2,\ldots,n\}$. Define
\begin{align}\label{hua}
&\hat\dcov_\nu^2(X,Y,\Pi,\K)=\frac{1}{\lfloor n/2\rfloor^3}\sum_{\ell_1,\ell_2,\ell_3=1}^{\lfloor n/2\rfloor}\big\{\omega_\nu(X_{\K_{\ell_1}},X_{\K_{\lfloor n/2\rfloor+\ell_1}},X_{\K_{\ell_2}},X_{\K_{\lfloor n/2\rfloor+\ell_2}})\notag\\
&\hspace{6cm}\times\omega_\nu(Y_{\Pi(\K_{\ell_2})},Y_{\Pi(\K_{\lfloor n/2\rfloor+\ell_2})},Y_{\Pi(\K_{\ell_3})},Y_{\Pi(\K_{\lfloor n/2\rfloor+\ell_3})})\big\}\,.
\end{align}
It holds that
\begin{align}\label{p3}
\hat\dcov_\nu^2(X,Y,\Pi,\K)=\E_{\K}\{\hat\dcov_\nu^2(X,Y,\Pi,\K)|\Pi,\T_n\}\,.
\end{align}
Since $\Pi$ follows a uniform distribution on $I_n(n)$, then, for $1\leq\ell\leq\lfloor n/2\rfloor$, it is true that $(Y_{\Pi(\K_{\ell})},Y_{\Pi(\K_{\lfloor n/2\rfloor+\ell})})$ is an exchangeable pair, that is,  
exchanging $Y_{\Pi(\K_{\ell})}$ and $Y_{\Pi(\K_{\lfloor n/2\rfloor+\ell})}$ does not alter the distribution of $\hat\dcov_\nu^2(X,Y,\Pi,\K,\delta)$. To be more specific, let $\delta_1,\ldots,\delta_{\lfloor n/2\rfloor}$ be i.i.d.~Bernoulli random variables such that $\P(\delta_j=1)=\P(\delta_j=0)=1/2$. Define
\begin{align*}
&\hat\dcov_\nu^2(X,Y,\Pi,\K,\delta)=\frac{1}{\lfloor n/2\rfloor^3}\sum_{\ell_1,\ell_2,\ell_3=1}^{\lfloor n/2\rfloor}\big\{\omega_\nu(X_{\K_{\ell_1}},X_{\K_{\lfloor n/2\rfloor+\ell_1}},X_{\K_{\ell_2}},X_{\K_{\lfloor n/2\rfloor+\ell_2}})\notag\\
&\hspace{5.5cm}\times\omega_\nu( W^{(\Pi,\K,\delta)}_{\ell_2},W^{(\Pi,\K,\delta)}_{\ell_3},W^{(\Pi,\K,\delta)}_{\lfloor n/2\rfloor+\ell_2},W^{(\Pi,\K,\delta)}_{\lfloor n/2\rfloor+\ell_3})\big\}\,,
\end{align*}
where, for $1\leq\ell\leq\lfloor n/2\rfloor$,
\begin{align*}
W^{(\Pi,\K,\delta)}_{\ell}=\delta_\ell Y_{\Pi(\K_{\ell})}+(1-\delta_\ell)Y_{\Pi(\K_{\lfloor n/2\rfloor+\ell})}\,,\qquad W^{(\Pi,\K,\delta)}_{\lfloor n/2\rfloor+\ell}=(1-\delta_\ell) Y_{\Pi(\K_{\ell})}+\delta_\ell Y_{\Pi(\K_{\lfloor n/2\rfloor+\ell})}\,.
\end{align*}
Then, $\hat\dcov_\nu^2(X,Y,\Pi,\K)$ and $\hat\dcov_\nu^2(X,Y,\Pi,\K,\delta)$ are identically distributed.

Furthermore, observe that $\omega_\nu$ in \eqref{omeganu} satisfies that, for $y_1,y_2,y_3,y_4\in\H$,
\begin{align}
&\omega_\nu(y_1,y_2,y_3,y_4)=-\omega_\nu(y_3,y_2,y_1,y_4)=-\omega_\nu(y_1,y_4,y_3,y_2)\,.\label{omega1}
\end{align}
We introduce i.i.d.~Rademacher random variables $\xi_1,\ldots,\xi_{\lfloor n/2\rfloor}$, that is, $\P(\xi_j=1)=\P(\xi_j=-1)=1/2$. Write $\xi=(\xi_1,\ldots,\xi_{\lfloor n/2\rfloor})$ and define
\begin{align*}
&\hat\dcov_\nu^2(X,Y,\Pi,\K,\xi)=\frac{1}{\lfloor n/2\rfloor^3}\sum_{\ell_1,\ell_2,\ell_3=1}^{\lfloor n/2\rfloor}\xi_{\ell_2}\xi_{\ell_3}\big\{\omega_\nu(X_{\K_{\ell_1}},X_{\K_{\lfloor n/2\rfloor+\ell_1}},X_{\K_{\ell_2}},X_{\K_{\lfloor n/2\rfloor+\ell_2}})\notag\\
&\hspace{6.5cm}\times\omega_\nu(Y_{\Pi(\K_{\ell_2})},Y_{\Pi(\K_{\lfloor n/2\rfloor+\ell_2})},Y_{\Pi(\K_{\ell_3})},Y_{\Pi(\K_{\lfloor n/2\rfloor+\ell_3})})\big\}\,.
\end{align*}
Then, in view of \eqref{omega1}, we deduce that $\hat\dcov_\nu^2(X,Y,\Pi,\K,\delta)$ and $\hat\dcov_\nu^2(X,Y,\Pi,\K,\xi)$ have the same distribution, so that
\begin{align}\label{p8}
\hat\dcov_\nu^2(X,Y,\Pi,\K)\ \text{and } \hat\dcov_\nu^2(X,Y,\Pi,\K,\xi) \text{ are identically distributed.}
\end{align}
Therefore, it follows from the Chernoff bound and \eqref{p3} that, for any $t>0$ and $\nu\in\V$,
\begin{align}\label{p9}
\P_\Pi\{\hat\dcov_\nu^2(X,Y,\Pi)\geq q|\T_n\}&\leq e^{-tq}\,\E_\Pi\big[\exp\{t\hat\dcov_\nu^2(X,Y,\Pi)\}|\T_n\big] \notag\\
&=e^{-tq}\,\E_\Pi\Big(\exp\big[t\E_{\K}\{\hat\dcov_\nu^2(X,Y,\Pi,\K)|\Pi,\T_n\}\big]\big|\T_n\Big) \notag\\
&\leq e^{-tq}\,\E_{\Pi,\K}\big[\exp\{t\hat\dcov_\nu^2(X,Y,\Pi,\K)\}|\T_n\big]\notag \\
&=e^{-tq}\,\E_{\Pi,\K,\xi}\big[\exp\{t\hat\dcov_\nu^2(X,Y,\Pi,\K,\xi)\}|\T_n\big]\,, 
\end{align}
where the first line follows from the Chernoff bound, the second line is due to \eqref{p3}, the third line is due to Jensen's inequality and the last line is due to \eqref{p8}.

For $1\leq\ell_1,\ell_2\leq\lfloor n/2\rfloor$, let
\begin{align*}
a_{\ell_2,\ell_3}(\T_n,\Pi,\K)&=\frac{1}{\lfloor n/2\rfloor^3}\,\omega_\nu(Y_{\Pi(\K_{\ell_2})},Y_{\Pi(\K_{\lfloor n/2\rfloor+\ell_2})},Y_{\Pi(\K_{\ell_3})},Y_{\Pi(\K_{\lfloor n/2\rfloor+\ell_3})})\\
&\hspace{2cm}\times\sum_{\ell_1=1}^{{\lfloor n/2\rfloor}}\omega_\nu(X_{\K_{\ell_1}},X_{\K_{\lfloor n/2\rfloor+\ell_1}},X_{\K_{\ell_2}},X_{\K_{\lfloor n/2\rfloor+\ell_2}})\,.
\end{align*}
Let $A_{\T_n,\Pi,\K}$ be the $\lfloor n/2\rfloor\times \lfloor n/2\rfloor$ matrix with entries $a_{\ell_1,\ell_2}(\T_n,\Pi,\K)$. Now, we have $$\hat\dcov_\nu^2(X,Y,\Pi,\K,\xi)=\xi'A_{\T_n,\Pi,\K}\xi.$$
By the Hanson--Wright inequality (e.g.~Theorem~1.1 of \cite{Rudelson2013}), there exist constants $c_1,c_2>0$ such that, for $0\leq t\leq c_2\|A_{\T_n,\Pi,\K}\|^{-1}$,
\begin{align}\label{p9.5}
&\E_{\Pi,\K,\xi}\big[\exp\{t\hat\dcov_\nu^2(X,Y,\Pi,\K,\xi)\}|\T_n\big]\leq \E_{\Pi,\K}\big\{\exp(c_1t\|A_{\T_n,\Pi,\K}\|_{\rm F}^2)\,|\T_n\big\}\,.
\end{align}
Observing the fact that $\|\phi_{\nu}\|_\infty\leq1$ for any $\nu\in\V$, by a line-by-line check of the proof of Theorem~1.1 of \cite{Rudelson2013}, it holds that the constants $c_1,c_2$ in equation \eqref{p9.5} is uniform in $\nu\in\V$.
Then, combining \eqref{p9} and \eqref{p9.5}, minimizing with respect to $t\in[0,c_2\|A_{\T_n,\Pi,\K}\|^{-1}]$ yields
\begin{align}\label{p11}
\P_\Pi\{\hat\dcov_\nu^2(X,Y,\Pi)\geq q|\T_n\}&\leq \E_{\Pi,\K}\bigg[\exp\bigg\{-c_3\min\bigg(\frac{q^2}{\|A_{\T_n,\Pi,\K}\|_{\rm F}^2},\frac{q}{\|A_{\T_n,\Pi,\K}\|}\bigg)\bigg\}\,\Big|\T_n\bigg]\notag\\
&\leq \E_{\Pi,\K}\bigg[\exp\bigg\{-c_3\min\bigg(\frac{q^2}{\|A_{\T_n,\Pi,\K}\|_{\rm F}^2},\frac{q}{\|A_{\T_n,\Pi,\K}\|_{\rm F}}\bigg)\bigg\}\,\Big|\T_n\bigg]
\end{align}
for some constant $c_3>0$ that is uniform in $\nu\in\V$, where in the last step we used the fact that $\|A_{\T_n,\Pi,\K}\|\leq\|A_{\T_n,\Pi,\K}\|_{\rm F}$.
Observing the definition of $\omega_\nu$ in \eqref{omega1}, by applying the Cauchy-Schwarz inequality and the fact that $\|\phi_{\nu}\|_\infty\leq 1$, we obtain, for any $\Pi$ and $\K$,
\begin{align}\label{p12}
\|A_{\T_n,\Pi,\K}\|_{\rm F}^2&=\frac{1}{\lfloor n/2\rfloor^6}\sum_{\ell_1,\ell_1',\ell_2,\ell_3=1}^{\lfloor n/2\rfloor}\omega_\nu^2(Y_{\Pi(\K_{\ell_2})},Y_{\Pi(\K_{\lfloor n/2\rfloor+\ell_2})},Y_{\Pi(\K_{\ell_3})},Y_{\Pi(\K_{\lfloor n/2\rfloor+\ell_3})})\notag\\
&\hspace{1cm}\times\omega_\nu(X_{\K_{\ell_1}},X_{\K_{\lfloor n/2\rfloor+\ell_1}},X_{\K_{\ell_2}},X_{\K_{\lfloor n/2\rfloor+\ell_2}})\,\omega_\nu(X_{\K_{\ell_1'}},X_{\K_{\lfloor n/2\rfloor+\ell_1'}},X_{\K_{\ell_2}},X_{\K_{\lfloor n/2\rfloor+\ell_2}})\notag\\
%
%
&\leq cn^{-4}\bigg[\sum_{r,\ell=1}^n\{\phi_{\nu}(X_r-X_\ell)\}^4\bigg]^{1/2}\bigg[\sum_{r,\ell=1}^n\{\phi_{\nu}(Y_r-Y_\ell)\}^4\bigg]^{1/2}\,.
\end{align}
Combining \eqref{p11} and \eqref{p12} yields that, for some absolute constant $c_3>0$,
\begin{align*}
\P_\Pi\{\hat\dcov_\nu^2(X,Y,\Pi)\geq q|\T_n\}&\leq \exp\bigg(-c_3qn^2\Big[\sum_{r,\ell=1}^n\{\phi_{\nu}(X_r-X_\ell)\}^4\Big]^{-1/4}\Big[\sum_{r,\ell=1}^n\{\phi_{\nu}(Y_r-Y_\ell)\}^4\Big]^{-1/4}\bigg)\,.
\end{align*}
Recalling the definition of $\tilde q_{\nu,1-\alpha/2}$ in \eqref{tildeqnu}, we therefore obtain from the above equation that
\begin{align*}
\tilde q_{\nu,1-\alpha/2}\leq c\log(2\alpha^{-1})n^{-2}\Big[\sum_{r,\ell=1}^n\{\phi_{\nu}(X_r-X_\ell)\}^4\Big]^{1/4}\Big[\sum_{r,\ell=1}^n\{\phi_{\nu}(Y_r-Y_\ell)\}^4\Big]^{1/4}\,.
\end{align*}
Therefore, by Markov's inequality and the fact that $\|\phi_{\nu}\|_\infty\leq1$, we obtain from the above equation that
\begin{align*}
\P_F(\tilde q_{\nu,1-\alpha/2}\geq q)&\leq ce^{-nq}\alpha^{-1}\E_F\bigg\{\exp\bigg(n^{-1}\Big[\sum_{r,\ell=1}^n\{\phi_{\nu}(X_r-X_\ell)\}^4\Big]^{1/4}\Big[\sum_{r,\ell=1}^n\{\phi_{\nu}(Y_r-Y_\ell)\}^4\Big]^{1/4}\bigg)\bigg\}\\
&\leq ce^{-nq}\alpha^{-1}\,.
\end{align*}
Recalling that $\tilde w_{\nu,1-\beta/3}$ is the $(1-\beta/3)$-quantile of $\tilde q_{\nu,1-\alpha/2}$, the above equation implies that
\begin{align}\label{p30}
\tilde w_{\nu,1-\beta/3}\leq cn^{-1}\{\log(\alpha^{-1})+\log(\beta^{-1})\}\,,
\end{align}
where the constant $c>0$ in the above equation is uniform in $\nu\in\V$.

\textbf{Step 2.} Bounding the quantile $\kappa_{\nu,1-\beta/3}$

We start by writing $\hat\dcov_\nu^2(X,Y)$ as a linear combination of U-statistics.
Recalling the definitions of $\hat\dcov_\nu^2(X,Y)$ and the empirical distance covariance in \eqref{v}, we have 
\begin{align}\label{d1}
\hat\dcov_\nu^2(X,Y)=\hat D_{\nu,1}(\T_n)+\hat D_{\nu,2}(\T_n)-2\hat D_{\nu,3}(\T_n)\,,
\end{align}
where
\begin{align*}
\hat D_{\nu,1}(\T_n)&=\frac{1}{n^4}\sum_{j_1,{j_2},{j_3},j_4=1}^n \phi_{\nu}(X_{j_1}-X_{j_2}) \,\phi_{\nu}(Y_{j_3}-Y_{j_4})\,,\\
\hat D_{\nu,2}(\T_n)&=\frac{1}{n^2}\sum_{{j_1},{j_2}=1}^n \phi_{\nu}(X_{j_1}-X_{j_2}) \,\phi_{\nu}(Y_{j_1}-Y_{j_2})\,,\\
\hat D_{\nu,3}(\T_n)&=\frac{1}{n^3}\sum_{{j_1},{j_2},{j_3}=1}^n \phi_{\nu}(X_{j_1}-X_{j_2}) \,\phi_{\nu}(Y_{j_1}-Y_{j_3})\,.
\end{align*}
For $k=2,3,4$, let $|I_{k}(n)|=k!{n\choose k}$ be the cardinality of $I_k(n)$.
Define 
\begin{align*}
\tilde D_{\nu,1}(\T_n)&=|I_4(n)|^{-1}\sum_{(j_1,{j_2},{j_3},j_4)\in I_4(n)} \phi_{\nu}(X_{j_1}-X_{j_2}) \,\phi_{\nu}(Y_{j_3}-Y_{j_4})\,,\\
%
\tilde D_{\nu,2}(\T_n)&=|I_2(n)|^{-1}\sum_{({j_1},{j_2})\in I_2(n)} \phi_{\nu}(X_{j_1}-X_{j_2}) \,\phi_{\nu}(Y_{j_1}-Y_{j_2})\,,\\
\tilde D_{\nu,3}(\T_n)&=|I_3(n)|^{-1}\sum_{({j_1},{j_2},{j_3})\in I_3(n)} \phi_{\nu}(X_{j_1}-X_{j_2})\, \phi_{\nu}(Y_{j_1}-Y_{j_3})\,,
\end{align*}
and let $\tilde D_\nu(\T_n)=\tilde D_{\nu,1}(\T_n)+\tilde D_{\nu,2}(\T_n)-2\tilde D_{\nu,3}(\T_n)$.

For $\hat D_{\nu,1}(\T_n)$, when $|\{j_1,j_2,j_3,j_4\}|\leq 3$, the cases of $\phi_{\nu}(X_{j_1}-X_{j_2})\phi_{\nu}(Y_{j_3}-Y_{j_4})\neq0$ is consist of: (i) $|\{j_1,j_2,j_3,j_4\}|=3$ and $|\{j_1,j_2\}\cap\{j_3,j_4\}|=1$; (ii) $|\{j_1,j_2,j_3,j_4\}|=2$ and $|\{j_1,j_2\}\cap\{j_3,j_4\}|=2$. 
Therefore, we obtain
\begin{align}\label{d2}
&n^4\hat D_{\nu,1}(\T_n)-|I_4(n)|\tilde D_{\nu,1}(\T_n)\notag\\
&=4\sum_{(j_1,j_2,j_3)\in I_3(n)}\phi_{\nu}(X_{j_1}-X_{j_2})\,\phi_{\nu}(Y_{j_1}-Y_{j_3})+2\sum_{(j_1,j_2)\in I_2(n)}\phi_{\nu}(X_{j_1}-X_{j_2})\,\phi_{\nu}(Y_{j_1}-Y_{j_2})\notag\\
&=4|I_3(n)|\tilde D_{\nu,3}(\T_n)+2|I_2(n)|\tilde D_{\nu,2}(\T_n)\,.
\end{align}
For $\hat D_{\nu,2}(\T_n)$, when $|\{j_1,j_2\}|=1$, it holds that $\phi_{\nu}(X_{j_1}-X_{j_2}) \phi_{\nu}(Y_{j_1}-Y_{j_2})=0$. Hence, we deduce that
\begin{align}\label{d3}
n^2\hat D_{\nu,2}(\T_n)-|I_2(n)|\tilde D_{\nu,2}(\T_n)=0\,.
\end{align}
For $\hat D_{\nu,3}(\T_n)$, when $|\{j_1,j_2,j_3\}|\leq2$, $\phi_{\nu}(X_{j_1}-X_{j_2})\phi_{\nu}(Y_{j_1}-Y_{j_3})\neq0$ if and only if $j_2=j_3\neq j_1$. Therefore,
\begin{align}\label{d4}
n^{3}\hat D_{\nu,3}(\T_n)-|I_3(n)|\tilde D_{\nu,3}(\T_n)=\sum_{(j_1,j_2)\in I_2(n)}\phi_{\nu}(X_{j_1}-X_{j_2})\phi_{\nu}(Y_{j_1}-Y_{j_2})=|I_2(n)|\tilde D_{\nu,2}(\T_n)\,.
\end{align}
Direct calculations by combining \eqref{d1}--\eqref{d4} yields
\begin{align}\label{p15}
\hat\dcov_\nu^2(X,Y)
&=\frac{n-1}{n}\bigg\{\frac{(n-2)(n-3)}{n^2}\tilde D_{\nu,1}(\T_n)+\frac{n^2-2n+2}{n^2}\tilde D_{\nu,2}(\T_n)-\frac{2(n-2)^2}{n^2}\tilde D_{\nu,3}(\T_n)\bigg\}\notag\\
&=\frac{3!{n-1\choose3}}{n^3}\tilde D_\nu(\T_n)+\frac{(n-1)(2n-4)}{n^3}\{\tilde D_{\nu,2}(\T_n)-\tilde D_{\nu,3}(\T_n)\}\,.
\end{align}

We therefore deduce from the above equation that $\hat\dcov_\nu^2(X,Y)-\E_F\hat\dcov_\nu^2(X,Y)=J_{\nu,n,1}+J_{\nu,n,2}$,
where
\begin{align}\label{jnun12}
J_{\nu,n,1}&=\frac{3!{n-1\choose3}}{n^3}\{\tilde D_\nu(\T_n)-\E\tilde D_\nu(\T_n)\}\,,\notag\\
J_{n,\nu,2}&=\frac{(n-1)(2n-4)}{n^3}\{\tilde D_{\nu,2}(\T_n)-\tilde D_{\nu,3}(\T_n)-\E\tilde D_{\nu,2}(\T_n)+\E\tilde D_{\nu,3}(\T_n)\}\,.
\end{align}
Therefore, we obtain that, for any $q>0$,
\begin{align*}
\P_F\big[|\hat\dcov_\nu^2(X,Y)-\E_F\{\hat\dcov_\nu^2(X,Y)\}|\geq q\big]\leq \P_F\big(|J_{\nu,n,1}|\geq q/2\big)+\P_F\big(|J_{\nu,n,2}|\geq q/2\big)\,.
\end{align*}
For $J_{\nu,n,1}$ in \eqref{jnun12}, observe the definition of the symmetric kernel $h$ in \eqref{h}, we have
\begin{align*}
\tilde D_\nu(\T_n)=\frac{1}{4!{n\choose 4}}\sum_{(j_1,j_2,j_3,j_4)\in I_4(n)}h\{(X_{j_1},Y_{j_1}),(X_{j_2},Y_{j_2}),(X_{j_3},Y_{j_3}),(X_{j_4},Y_{j_4})\}\,.
\end{align*}
By the Bernstein-type concentration inequality for general U-statistics in Theorem~2 in \cite{arcones1995}, it follows that for any $q\geq0$, $\nu\in\V$, and $F\in\mathcal F_\nu(\rho)$,
\begin{align}\label{q}
\P_F\big\{|\tilde D_\nu(\T_n)-\E_F\tilde D_\nu(\T_n)|\geq q\big\}\leq 4\exp\bigg(-\frac{nq^2}{2m^2\zeta_1+(2^{m+3}m^m+2m^{-1}/3)q}\bigg)\,,
\end{align}
where $\zeta_1=\var_F\big(\E_F\big[h\{(X_{1},Y_{1}),(X_{2},Y_{2}),(X_{3},Y_{3}),(X_{4},Y_{4})\}|(X_{1},Y_{1})\big]\big)$. In view of the definition of $h$ in \eqref{h} and the fact that $\|\phi_{\nu}\|_\infty\leq1$, we observe the fact that $\|h\|_\infty\leq4$ for any $F\in\mathcal F_\nu(\rho)$.
Hence, by the law of total variance, it follows that
\begin{align*}
\zeta_1&\leq\var_F[h\{(X_{1},Y_{1}),(X_{2},Y_{2}),(X_{3},Y_{3}),(X_{4},Y_{4})\}]\leq\|h\|_\infty^2\leq64\,.
\end{align*}
Since $|J_{\nu,n,1}|$ is bounded, so that we may assume $q\leq c$. In addition, when $n\geq4$, it holds that $3!{n-1\choose3}n^{-3}\geq 3/32$. We therefore deduce in view of \eqref{q} that, for any $q\geq0$,
\begin{align}\label{q2}
\P_F\big(|J_{\nu,n,1}|\geq q/2\big)&=\P_F\bigg\{|\tilde D_\nu(\T_n)-\E_F\tilde D_\nu(\T_n)|\geq \frac{n^3q}{2\times3!{n-1\choose3}}\bigg\}\notag\\
&\leq 4\exp\bigg(-\frac{9nq^2}{2^{25}+(2^{27}+2^{11}/3)q}\bigg)\leq 4\exp(-cnq^2)\,.
\end{align}
provided that $n\geq4$, where $c>0$ does not depend on $n$, or $\nu\in\V$ or $F\in\mathcal F_\nu(\rho)$.

Next, for  $J_{\nu,n,2}$ in \eqref{jnun12}, observe that $|\tilde D_{\nu,2}(\T_n)|\leq1$ and $|\tilde D_{\nu,3}(\T_n)|\leq1$ for any $\nu\in\V$ and $F\in\mathcal F_\nu(\rho)$, so that $2n|J_{\nu,n,2}|\leq 8$. An application of Markov's inequality directly yields that, for any $\nu\in\V$,  $F\in\mathcal F_\nu(\rho)$,  and $q>0$,
\begin{align*}
&\P_F(|J_{\nu,n,2}|\geq q/2)\leq e^{-nq}\E\exp(2n|J_{\nu,n,2}|)\leq \exp(16-nq)\,.
\end{align*}
Combining the above equation with \eqref{q2} yields that, for some constant $c>0$ that does not depend on $n$ and $\beta$,
\begin{align*}
\sup_{\nu\in\V}\sup_{F\in\mathcal F_\nu(\rho)}\P_F(|\hat\dcov_\nu^2(X,Y)-\E_F\hat\dcov_\nu^2(X,Y)|>q)\leq 4\exp(-cnq^2)+\exp(16-nq)\leq c_1\exp(-c_2nq^2)\,.
\end{align*}
By taking $q=c_2^{-1/2}\sqrt{n^{-1}\log(3c_1\beta^{-1})}$ in the above equation, we obtain that the above quantity is further bounded by $\beta/3$, which implies
\begin{align}\label{q3}
\sup_{\nu\in\V}\sup_{F\in\mathcal F_\nu(\rho)}\kappa_{\nu,1-\beta/3}\leq cn^{-1/2}\log(\beta^{-1})\,,
\end{align}
where $c>0$ does not depend on $n$ or $\beta$.

\textbf{Step 3.} Bounding the bias $|\|\phi_{X,Y}-\phi_X\otimes\phi_Y\|_\nu^2-\E_F\{\hat\dcov_\nu^2(X,Y)\}|$

Observe the fact that $\E_F\{\tilde D_\nu(\T_n)\}=\|\phi_{X,Y}-\phi_X\otimes\phi_Y\|_\nu^2\leq 4$, $|\E_F\{\tilde D_{\nu,2}(\T_n)\}|\leq 1$, and $|\E_F\{\tilde D_{\nu,3}(\T_n)\}|\leq 1$. 
From \eqref{p15} we obtain
\begin{align*}
&|\E_F\{\tilde D_\nu(\T_n)\}-\|\phi_{X,Y}-\phi_X\otimes\phi_Y\|_\nu^2|\\
&=\frac{-6n^2+11n-6}{n^3}\|\phi_{X,Y}-\phi_X\otimes\phi_Y\|_\nu^2+\frac{(n-1)(2n-4)}{n^3}\big[\E_F\{\tilde D_{\nu,2}(\T_n)\}-\E_F\{\tilde D_{\nu,3}(\T_n)\}\big]\leq cn^{-1}\,.
\end{align*}
We therefore deduce 
that, for some constant $c>0$ that does not depend on $\beta$, 
\begin{align*}
\sup_{\nu\in\V}\sup_{F\in\mathcal F_\nu(\rho)}|\E_F\{\tilde D_\nu(\T_n)\}-\|\phi_{X,Y}-\phi_X\otimes\phi_Y\|_\nu^2|\leq cn^{-1}\,.
\end{align*}
The proof is therefore complete by combining this bound with \eqref{qbeta}, \eqref{p30}, and \eqref{q3}.

\subsection{Proof of Theorem~\ref{thm:lower:ind}}\label{app:thm:lower:ind}

We follow Ingster's method (see \cite{ingster1987,ingster1993}) to derive the lower bound for $\rho^*$.
For any fixed test $\psi\in\Psi(\alpha)$, for the set of local alternative $\mathcal F_\nu(\rho)$ in \eqref{fnurhoa}, the maximum type-II error of $\psi$ is defined by
\begin{align*}
\beta_{n,\rho}(\psi)=\sup_{\nu\in\V}\sup_{F\in\mathcal F_\nu(\rho)}\P_F(\psi=0)\,.
\end{align*}
The minimax risk is then defined as the infimum value of $\beta_{n,\rho}(\psi)$ among all tests having controlled type-I error at level $\alpha$, that is,
\begin{align}\label{beta}
\beta_{n,\rho}^*=\inf_{\psi\in\Psi(\alpha)}\beta_{n,\rho}(\psi)=\inf_{\psi\in\Psi(\alpha)}\sup_{\nu\in\V}\sup_{F\in\mathcal F_\nu(\rho)}\P_F(\psi=0)\,.
\end{align}
Then, in view of the definition of $\rho^*(\alpha,\beta)$ in \eqref{rho*}, we obtain that, a sufficient condition for $\rho^*(\alpha,\beta)\geq \rho$ is given by $\beta_{n,\rho}^*\geq\beta$.

Let $\{\eta_k\}_{k\geq1}$ be an orthogonal basis on $\H$. Consider the following covariance operator on $\H$:
\begin{align}\label{cc}
C=\sum_{k=1}^\infty\lambda_k\l\eta_k,\cdot\r\eta_k\,,\qquad\text{where }\lambda_k=k^{-2}\,.
\end{align}
Suppose $\nu_0$ is a mean-zero Gaussian measure on $\H$ with covariance operator $C$ defined in \eqref{cc}.
For integer $M\geq2$ to be specified below, let $F_n$ denote the uniform measure over the set  $\{F_{n,2},\ldots,F_{n,M}\}$, which consists of $(M-1)$ joint distributions on $\H^2$, to be specified below as well.
Recalling the definition of $\beta_{n,\rho}^*$ in \eqref{beta} and the fact that $\P_{F_0}(\psi=1)\leq\alpha$ for any $\psi\in\Psi(\alpha)$, we deduce
\begin{align}\label{lb}
\beta_{n,\rho}^*&\geq\inf_{\psi\in\Psi(\alpha)}\sup_{F\in\mathcal F_{\nu_0}(\rho)}\P_{F}(\psi=0)\geq \inf_{\psi\in\Psi(\alpha)}\P_{F_n}(\psi=0)=1-\sup_{\psi\in\Psi(\alpha)}\P_{F_n}(\psi=1)\notag\\
&\geq1-\alpha-\sup_{\psi\in\Psi(\alpha)}|\P_{F_n}(\psi=1)-\P_{F_0}(\psi=1)|\geq1-\alpha-{\rm TV}(\P_{F_n},\P_{F_0})\,.
\end{align}
Here, the total variation distance between the probability measures $\P_{F_n}$ and $\P_{F_0}$ is defined as
\begin{align*}
{\rm TV}(\P_{F_n},\P_{F_0})=\sup_{A\in\mathcal A}|\P_{F_n}(A)-\P_{F_0}(A)|\,,
\end{align*}
where $\mathcal A$ denotes the space of measurable sets on $\H^{2n}$.
Next, observing the connection between the total variance distance and $L_1$-distance, direct calculations yield
\begin{align*}
{\rm TV}(\P_{F_n},\P_{F_0})=\frac{1}{2}\E_{\P_{F_0}}\bigg|\frac{\d\P_{F_n}(\T_n)}{\d\P_{F_0}(\T_n)}-1\bigg|
\leq\bigg[\E_{\P_{F_0}}\bigg\{\frac{\d\P_{F_n}(\T_n)}{\d\P_{F_0}(\T_n)}\bigg\}^2-1\bigg]^{1/2}\,.
\end{align*}
Here, $\frac{\d\P_{F_n}}{\d\P_{F_0}}$ denotes the Radon-Nikodym derivative of $F_n$ with respect to~the null distribution $F_0$.
Combining the above equation with \eqref{lb} yields
\begin{align*}
\beta_{n,\rho}^*\geq 1-\alpha-\bigg[\E_{\P_{F_0}}\bigg\{\frac{\d\P_{F_n}(\T_n)}{\d\P_{F_0}(\T_n)}\bigg\}^2-1\bigg]^{1/2}\,.
\end{align*}
Then, we deduce that $\beta_{n,\rho}^*\geq\beta$, provided that
\begin{align}\label{goal}
\E_{\P_{F_0}}\bigg\{\frac{\d\P_{F_n}(\T_n)}{\d\P_{F_0}(\T_n)}\bigg\}^2\leq1+4(1-\alpha-\beta)^2\,.
\end{align}

We then construct distributions $F_n,F_0$ on $\H^2$. Consider a set of $(M-1)$ joint distributions of random variables on $\H^2$, where $M\geq2$ is an integer to be specified later.
Let $\zeta_1^*,\zeta_2^*,\{\zeta_{\ell}\}_{\ell\geq1},\{\tilde\zeta_{\ell}\}_{\ell\geq1}$ be all independent and identically distributed real-valued standard normal random variables. For $2\leq k\leq M$, consider the following $\H$-valued random variables
\begin{align}\label{vw}
V_{(n,k)}&=\sqrt{\lambda_1}(\sqrt{\varrho_n}\zeta_1^*+\sqrt{1-\varrho_n}\zeta_{1})\eta_{1}+\sqrt{\lambda_k}(\sqrt{\varrho_n}\zeta_2^*+\sqrt{1-\varrho_n}\zeta_{k})\eta_{k}+\sum_{\ell\geq2,\ell\neq k}\sqrt{\lambda_\ell}\,\zeta_{\ell}\,\eta_\ell\,;\notag\\
W_{(n,k)}&=\sqrt{\lambda_1}(\sqrt{\varrho_n}\zeta_2^*+\sqrt{1-\varrho_n}\zeta_{1})\eta_{1}+\sqrt{\lambda_k}(\sqrt{\varrho_n}\zeta_1^*+\sqrt{1-\varrho_n}\tilde\zeta_{k})\eta_k+\sum_{\ell\geq2,\ell\neq k}\sqrt{\lambda_\ell}\,\tilde\zeta_{\ell}\,\eta_\ell\,,
\end{align}
where $\varrho_n\in(0,1)$ is some parameter to be specified below.
For $2\leq k\leq M$, let $F_{n,k}$ denote the joint distribution of $(V_{(n,k)},W_{(n,k)})$, and let  $\phi_{n,k}$ denote its characteristic functional. First, it follows from Lemma~\ref{lem:klexpan} that the marginal distribution of $V_{(n,k)}$ and $W_{(n,k)}$ are mean-zero Gaussian distribution with covariance operator $C$ defined in \eqref{cc}.

\begin{lemma}[Theorem~6.19 of \cite{stuart2010}]\label{lem:klexpan}
Let $\mathcal{C}$ be a self-adjoint, positive semi-definite, nuclear operator in a Hilbert space $\mathcal{H}$ and let $m \in \mathcal{H}$. Let $\{(\lambda_k,\eta_k)\}_{k=1}^{\infty}$ be an orthonormal set of eigenvalues and eigenvectors of $\mathcal{C}$ ordered so that
$
\lambda_1\geq \lambda_2 \geq \cdots.
$
Take $\{\zeta_k\}_{k=1}^{\infty}$ to be an i.i.d. sequence with $\zeta_1 \sim \mathcal{N}(0,1)$. Then the random variable $Z \in \mathcal{H}$ given by the Karhunen-Loève expansion
$
Z=m+\sum_{k=1}^{\infty} \sqrt{\lambda}_k \zeta_k \eta_k
$
is distributed according to $\mathcal{N}(m, \mathcal{C})$.
\end{lemma}

Next, for the joint distribution of $(V_{n,k},W_{n,k})$, in view of the definition of characteristic functionals and \eqref{vw}, direct calculations yields
\begin{align*}
&\phi_{n,k}(w_1,w_2)=\E\exp(\i\l V_{(n,k)},w_1\r+\i\l W_{(n,k)},w_2\r)\\
%
%
&=\prod_{\ell\geq2,\ell\neq k}\E\exp\big(\i\sqrt{\lambda_\ell}\l\eta_\ell,w_1\r\zeta_{\ell}\big)\times\prod_{\ell\geq2,\ell\neq k}\E\exp\big(\i\sqrt{\lambda_\ell}\l\eta_\ell,w_2\r\tilde\zeta_{\ell}\big)\notag\\
&\ \times\E\exp\{\i(\sqrt{\lambda_1}\l\eta_1,w_1\r+\sqrt{\lambda_k}\l\eta_k,w_2\r)\sqrt{\varrho_n}\zeta_1^*\}\times\E\exp\{\i(\sqrt{\lambda_1}\l\eta_1,w_2\r+\sqrt{\lambda_k}\l\eta_k,w_1\r)\sqrt{\varrho_n}\zeta_2^*\}\notag\\
&\ \times\E\exp(\i\sqrt{\lambda_k(1-\varrho_n)}\l\eta_k,w_1\r\zeta_k)\times\E\exp(\i\sqrt{\lambda_k(1-\varrho_n)}\l\eta_k,w_2\r\tilde\zeta_k)\,.
\end{align*}
Since $\zeta_1^*,\zeta_2^*,\{\zeta_{\ell}\}_{\ell\geq1},\{\tilde\zeta_{\ell}\}_{\ell\geq1}$ are i.i.d.~standard normal random variables, we deduce from the above equation that
{\small
\begin{align}\label{phink}
&\phi_{n,k}(w_1,w_2)=\exp\bigg[-\frac{1}{2}\bigg\{\sum_{\ell\geq2,\ell\neq k}\lambda_\ell\l\eta_\ell,w_1\r^2+\sum_{\ell\geq2,\ell\neq k}\lambda_\ell\l\eta_\ell,w_2\r^2+(\sqrt{\lambda_1}\l\eta_1,w_1\r+\sqrt{\lambda_k}\l\eta_k,w_2\r)^2\varrho_n\notag\\
&\hspace{2cm}+(\sqrt{\lambda_1}\l\eta_1,w_2\r+\sqrt{\lambda_k}\l\eta_k,w_1\r)^2\varrho_n+\lambda_k(1-\varrho_n)\l\eta_k,w_1\r^2+\lambda_k(1-\varrho_n)\l\eta_k,w_2\r^2\bigg\}\bigg]\notag\\
&=\exp\bigg\{-\frac{1}{2}\sum_{k=1}^\infty\lambda_k\l\eta_k,w_1\r^2-\frac{1}{2}\sum_{k=1}^\infty\lambda_k\l\eta_k,w_2\r^2-\varrho_n\sqrt{\lambda_1\lambda_k}(\l\eta_1,w_1\r\l\eta_k,w_2\r+\l\eta_k,w_1\r\l\eta_1,w_2\r)\bigg\}\,.
\end{align}
}
Recall the covariance operator $C$ and its eigenvalues $\lambda_k$ defined in \eqref{cc}, and let the covariance operator on $\H$ defined as
\begin{align}\label{ccnk}
C_{n,k}=\varrho_n\sqrt{\lambda_1\lambda_k}(\l\cdot,\eta_1\r\eta_k+\l\cdot,\eta_k\r\eta_1)\,.
\end{align}
By Theorem~2 in \cite{baker1973}, the joint Gaussian distribution $F_{n,k}$ of $(V_{(n,k)},W_{(n,k)})$ in \eqref{vw} is the Gaussian measure $\H^2$ with covariance operator
\begin{align}\label{cnuk}
C_{F_{n,k}}=\Bigg[\begin{matrix}
C & C_{n,k}\\
C_{n,k}& C
\end{matrix}\Bigg]\,;
\end{align}
that is, $C_{F_{n,k}}[w_1,w_2]\trans=[Cw_1+C_{n,k}w_2,C_{n,k}w_1+Cw_2]\trans$, for $(w_1,w_2)\in\H^2$.
In addition, observe
\begin{align}\label{cnk}
\l C_{n,k} w_1,w_2\r=\varrho_n\sqrt{\lambda_1\lambda_k}(\l w_1,\eta_1\r\l w_2,\eta_k\r+\l w_1,\eta_k\r\l w_2,\eta_1\r)\,,
\end{align}
so that in view of \eqref{phink} we obtain
\begin{align*}
\phi_{n,k}(w_1,w_2)=\exp\big(-2^{-1}\l Cw_1,w_1\r-2^{-1}\l Cw_2,w_2\r-\l C_{n,k}w_1,w_2\r\big)\,.
\end{align*}
In this case, the null distribution $F_0$ is the mean-zero Gaussian distribution with covariance operator
\begin{align}\label{cf0}
C_{F_{0}}=\Bigg[\begin{matrix}
C & 0\\
0& C
\end{matrix}\Bigg]\,.
\end{align}

Next, we compute the distance covariance $\dcov_{\nu_0}^2(V_{(n,k)},W_{(n,k)})$ between $V_{(n,k)}$ and $W_{(n,k)}$ defined in \eqref{vw}, with respect to~the reference probability measure $\nu_0$.
We obtain from the above equation that
\begin{align*}
&\phi_{n,k}(w_1,w_2)-\phi_{n,k}(w_1,0)\phi_{n,k}(0,w_2)\\
&=\exp\big(-2^{-1}\l Cw_1,w_1\r-2^{-1}\l Cw_2,w_2\r\big)\big\{\exp(-\l C_{n,k}w_1,w_2\r)-1\big\}\\
&\geq\frac{9}{16}|\l C_{n,k}w_1,w_2\r|^2\one\{|\l C_{n,k}w_1,w_2\r|\leq1/2\}+(1-e^{-1})\one\{|\l C_{n,k}w_1,w_2\r|>1/2\}\,.
\end{align*}
We use $I_{n,k}=\dcov_{\nu_0}^2(V_{(n,k)},W_{(n,k)})$ to denote the squared distance covariance for simplicity, and obtain that, for $2\leq k\leq M$,
\begin{align}\label{ink}
I_{n,k}&=\int_{\H^2}|\phi_{n,k}(w_1,w_2)-\phi_{n,k}(w_1,0)\phi_{n,k}(0,w_2)|^2\,{\nu_0}(\d w_1)\,{\nu_0}(\d w_2)\notag\\
&=\int_{\H^2}\exp(-\l Cw_1,w_1\r-\l Cw_2,w_2\r)|\exp(-\l C_{n,k}w_1,w_2\r)-1|^2\,{\nu_0}(\d w_1)\,{\nu_0}(\d w_2)\,.
\end{align}
Observe the fact that, for $a\in\mathbb R$, when $|a|\leq1/2$, it is true that $|1-e^{a}|\geq3|a|/4$; when $|a|>1/2$, it holds that $|1-e^{a}|\geq1-e^{-1/2}$. 
Therefore, we obtain
\begin{align*}
|\exp(-\l C_{n,k}w_1,w_2\r)-1|^2&\geq\frac{9}{16}|\l C_{n,k}w_1,w_2\r|^2\,\one\{|\l C_{n,k}w_1,w_2\r|\leq1/2\}\\
&\qquad ~~+(1-e^{-1/2})^2\,\one\{|\l C_{n,k}w_1,w_2\r|>1/2\}\,.
\end{align*}
Therefore, combining the above equation with \eqref{ink} yields $I_{n,k}\geq I_{n,k,1}+I_{n,k,2}$,
where
\begin{align}\label{ink12}
I_{n,k,1}&=\frac{9}{16}\int_{\H^2}\exp(-\l Cw_1,w_1\r-\l Cw_2,w_2\r)|\l C_{n,k}w_1,w_2\r|^2\one\{|\l C_{n,k}w_1,w_2\r|\leq1/2\}{\nu_0}(\d w_1){\nu_0}(\d w_2)\,,\notag\\
I_{n,k,2}&=(1-e^{-1/2})^2\int_{\H^2}\exp(-\l Cw_1,w_1\r-\l Cw_2,w_2\r)\one\{|\l C_{n,k}w_1,w_2\r|>1/2\}{\nu_0}(\d w_1){\nu_0}(\d w_2)\,.
\end{align}
For the first term $I_{n,k,1}$, observing the fact that $\exp(\l Cw_1,w_1\r+\l Cw_2,w_2\r)\geq1$, we obtain directly from \eqref{cnk} that
\begin{align}\label{ink1}
I_{n,k,1}&\geq\frac{9}{16}\varrho_n^2\lambda_1\lambda_k\int_{\H^2}\Big[\exp(-2\l Cw_1,w_1\r-2\l Cw_2,w_2\r)(\l w_1,\eta_1\r\l w_2,\eta_k\r+\l w_1,\eta_k\r\l w_2,\eta_1\r)^2\notag\\
&\qquad\hspace{3cm}\times\one\{|\l C_{n,k}w_1,w_2\r|\leq1/2\}\Big]\,{\nu_0}(\d w_1)\,{\nu_0}(\d w_2)\,.
\end{align}
For the second term $I_{n,k,2}$ in \eqref{ink12}, suppose that
\begin{align*}
\varrho_n^2\leq4\lambda_1^{-1}\min_{1\leq k\leq M}\lambda_k=4\lambda_1^{-1}M^{-2}\,.
\end{align*}
Then it follows that
\begin{align*}
&\lambda_1\lambda_k(\l w_1,\eta_1\r\l w_2,\eta_k\r+\l w_1,\eta_k\r\l w_2,\eta_1\r)^2\\
&\leq 4^{-1}\lambda_1\lambda_k(\l w_1,\eta_1\r^2+\l w_2,\eta_k\r^2+\l w_1,\eta_k\r^2+\l w_2,\eta_1\r^2)^2\\
&\leq4^{-1}\lambda_1\lambda_k^{-1}\big\{\lambda_1(\l w_1,\eta_1\r^2+\l w_2,\eta_1\r^2)+\lambda_k(\l w_1,\eta_k\r^2+\l w_2,\eta_k\r^2)\big\}^2\\
&\leq\varrho_n^{-2}\big\{\lambda_1(\l w_1,\eta_1\r^2+\l w_2,\eta_1\r^2)+\lambda_k(\l w_1,\eta_k\r^2+\l w_2,\eta_k\r^2)\big\}^2\,.
\end{align*}
By observing the fact that $\exp(a)\geq a^2$ for any $a>0$, we obtain from the above equation that
\begin{align*}
&\exp(\l Cw_1,w_1\r+\l Cw_2,w_2\r)=\exp\Big\{\sum_{j=1}^\infty\lambda_j(\l w_1,\eta_j\r^2+\l w_2,\eta_j\r^2)\Big\}\\
&\geq\exp\big\{\lambda_1(\l w_1,\eta_1\r^2+\l w_2,\eta_1\r^2)+\lambda_k(\l w_1,\eta_k\r^2+\l w_2,\eta_k\r^2)\big\}\\
&\geq\big\{\lambda_1(\l w_1,\eta_1\r^2+\l w_2,\eta_1\r^2)+\lambda_k(\l w_1,\eta_k\r^2+\l w_2,\eta_k\r^2)\big\}^2 \\
&\geq\varrho_n^2{\lambda_1\lambda_k}(\l w_1,\eta_1\r\l w_2,\eta_k\r+\l w_1,\eta_k\r\l w_2,\eta_1\r)^2\,.
\end{align*}
This equation implies that
\begin{align*}
I_{n,k,2}&\geq (1-e^{-1/2})^2\varrho_n^2\lambda_1\lambda_k\int_{\H^2}\Big[\exp(-2\l Cw_1,w_1\r-2\l Cw_2,w_2\r)\\
&\hspace{1cm}\times(\l w_1,\eta_1\r\l w_2,\eta_k\r+\l w_1,\eta_k\r\l w_2,\eta_1\r)^2\one\{|\l C_{n,k}w_1,w_2\r|>1/2\}\Big]\,{\nu_0}(\d w_1)\,{\nu_0}(\d w_2)\,.
\end{align*}
Therefore, combining the above equation with \eqref{ink1} yields
\begin{align}\label{ink0}
I_{n,k}&\geq (1-e^{-1/2})^2\varrho_n^2\lambda_1\lambda_k\int_{\H^2}\Big[\exp(-2\l Cw_1,w_1\r-2\l Cw_2,w_2\r)\notag\\
&\hspace{5cm}\times(\l w_1,\eta_1\r\l w_2,\eta_k\r+\l w_1,\eta_k\r\l w_2,\eta_1\r)^2\Big]\,{\nu_0}(\d w_1)\,{\nu_0}(\d w_2)\notag\\
&=2(1-e^{-1/2})^2\varrho_n^2\lambda_1\lambda_k\int_\H\exp(-2\l Cw,w\r)\l w,\eta_1\r^2{\nu_0}(\d w)\int_\H\exp(-2\l Cw,w\r)\l w,\eta_k\r^2{\nu_0}(\d w)\notag\\
&\quad+(1-e^{-1/2})^2\varrho_n^2\lambda_1\lambda_k\bigg\{\int_\H\exp(-2\l Cw,w\r)\l w,\eta_1\r\l w,\eta_k\r{\nu_0}(\d w)\bigg\}^2\notag\\
&\geq 2(1-e^{-1/2})^2\varrho_n^2G_{n,1}G_{n,k}\,.
\end{align}
where, for $2\leq k\leq M$,
\begin{align*}
G_{n,k}=\lambda_k\int_\H\exp(-2\l Cw,w\r)\l w,\eta_k\r^2\,{\nu_0}(\d w)\,.
\end{align*}
In order to compute $G_{n,k}$, in the sequel, for simplicity, we use $\N_{\lambda}$ to denote mean-zero Gaussian measure on $(\mathbb R,\mathscr B(\mathbb R))$ with variance $\lambda>0$, that is, $\N_\lambda(B)=(2\pi\lambda)^{-1/2}\int_B\exp(-2^{-1}\lambda^{-1}x^2)\d x$ for $B\in\mathscr B(\mathbb R)$.
Observe that direct calculations yields
\begin{align}\label{t1}
&\int_{\mathbb R}\exp(-2\lambda_kx^2)\N_{\lambda_k}(\d x)=\dfrac{1}{(4\lambda_k^2+1)^{1/2}},~~~\lambda_k\int_{\mathbb R}\exp(-2\lambda_kx^2)x^2\N_{\lambda_k}(\d x)=\frac{\lambda_k^2}{(4\lambda_k^2+1)^{3/2}}.
\end{align}
Let $P_m$ denote the projection mapping on $\H$ defined by $P_m(x)=\sum_{j=1}^m\l x,\eta_j\r\eta_j$, for $x\in\H$.
It is true that $x=\lim_{m\to\infty}P_m(x)$, so that by the dominated convergence theorem, 
\begin{align*}
G_{n,k}=\lambda_k\lim_{m\to\infty}\int_\H\exp\{-2\l CP_m(w),P_m(w)\r\}\l P_m(w),\eta_k\r^2\,{\nu_0}(\d w)\,.
\end{align*}
For any $2\leq k\leq M$ and $m>k$, in view of \eqref{t1}, we obtain
\begin{align*}
&\lambda_k\int_\H\exp\{-2\l C(P_mw),P_m(w)\r\}\l P_m(w),\eta_k\r^2\,{\nu_0}(\d w)\\
&=\lambda_k\int_\H\exp\bigg(-2\sum_{j=1}^m\lambda_j\l w,\eta_j\r^2\bigg)\l w,\eta_k\r^2{\nu_0}(\d w)\\
&=\lambda_k\int_{\mathbb R}\exp(-2\lambda_k w_k^2)w_k^2\,\N_{\lambda_k}(\d w_k)\times\prod_{1\leq j\leq m;\,j\neq k}\int_{\mathbb R}\exp(-2\lambda_j w_j^2)\,\N_{\lambda_j}(\d w_j)\\
&=\frac{\lambda_k\int_{\mathbb R}\exp(-2\lambda_k w_k^2)w_k^2\,\N_{\lambda_k}(\d w_k)}{\int_{\mathbb R}\exp(-2\lambda_k w_k^2)\,\N_{\lambda_k}(\d w_k)}\times\prod_{j=1}^m\int_{\mathbb R}\exp(-2\lambda_j w_j^2)\,\N_{\lambda_j}(\d w_j)\\
&=\frac{\lambda_k^2}{4\lambda_k^2+1}\times\prod_{j=1}^m\int_{\mathbb R}\exp(-2\lambda_j w_j^2)\,\N_{\lambda_j}(\d w_j)\,.
\end{align*}
The above equation implies that
\begin{align*}
G_{n,k}&=\frac{\lambda_k^2}{4\lambda_k^2+1}\times\lim_{m\to\infty}\int_\H\exp\{-2\l CP_m(w),P_m(w)\r\}\,{\nu_0}(\d w)\\
&=\frac{\lambda_k^2}{4\lambda_k^2+1}\times\int_\H\exp(-2\l Cw,w\r)\,{\nu_0}(\d w)=\frac{\lambda_k^2}{4\lambda_k^2+1}\times\E{\phi_{\nu_0}}(2\tilde X)\,,
\end{align*}
where $\tilde X$ is mean zero Gaussian with covariance operator $2C$. Observe that $\E{\phi_{\nu_0}}(2\tilde X)\geq c$ where $c>0$ do not depend on $n$. Combining the above equation with \eqref{ink0} yields that, for $2\leq k\leq M$
\begin{align}\label{p2}
\dcov_{\nu_0}^2(V_{(n,k)},W_{(n,k)})=I_{n,k}&\geq 2(1-e^{-1/2})^2\varrho_n^2\frac{\lambda_1^2\lambda_k^2}{(4\lambda_1^2+1)(4\lambda_k^2+1)}\{\E{\phi_{\nu_0}}(2\tilde X)\}^2\geq c\varrho_n^2\,.
\end{align}
This proves that $F_{n,2},\ldots,F_{n,M}\in\mathcal F_{\nu_0}(\rho)$, by taking $\rho=c\varrho_n^2$.

Finally, we verify \eqref{goal}.
We use the following lemma which gives an explicit expression for Radon-Nikodym derivatives between Gaussian measures; we refer to Corollary~6.4.11 of \cite{Bogachev1998} or Theorem~11 of \cite{minh2021} for a proof.
\begin{lemma}\label{lem:rd}
Suppose ${\nu}_1$ and ${\nu}_2$ are Gaussian measures on $\H$ with covariance operators $C_1$ and $C_2$, respectively, and ${\nu}_1$ is absolutely continuous with respect to~${\nu}_2$ (i.e., $\nu_1\ll\nu_2$). Then, the Radon-Nikodym derivative is given by
\begin{align*}
\frac{\d{\nu}_1(x)}{\d{\nu}_2(x)}
%
&=\exp\bigg[-\frac{1}{2}\sum_{k=1}^\infty\Big\{\log(1-\alpha_k)+\frac{\alpha_k}{1-\alpha_k}\l x,C_2^{-1/2}\xi_k\r^2\Big\}\bigg]\,,\qquad\text{for }x\in\H\,,
\end{align*}
where the $\alpha_k$'s are the eigenvalues of the operator $S=I-C_2^{-1/2}C_1C_2^{-1/2}$ and the $\xi_k$'s are the corresponding eigenvectors.

\end{lemma}

We apply Lemma~\ref{lem:rd} to compute the Radon-Nikodym derivative in \eqref{goal}.
Recalling the definitions of $C$ and $C_{n,k}$ in \eqref{ccnk}, we observe that $C^{-1/2}C_{n,k}C^{-1/2}=\varrho_n(\l\cdot,\eta_1\r\eta_k+\l\cdot,\eta_k\r\eta_1)$. Furthermore, $C^{-1/2}C_{n,k}C^{-1/2}$ is a rank 2 operator on $\H$, and
\begin{align*}
&C^{-1/2}C_{n,k}C^{-1/2}(\eta_1+\eta_k)=\varrho_n(\eta_1+\eta_k)\,,\\
&C^{-1/2}C_{n,k}C^{-1/2}(\eta_1-\eta_k)=-\varrho_n(\eta_1-\eta_k)\,.
\end{align*}
This implies that the nonzero eigenvalues of $C^{-1/2}C_{n,k}C^{-1/2}$ are $\varrho_n$ and $-\varrho_n$ with corresponding eigenvectors $2^{-1/2}(\eta_1+\eta_k)$ and $2^{-1/2}(\eta_1-\eta_k)$. Consider the following operator on $\H^2$:
\begin{align*}
S_{n,k}&:=I-\Bigg[\begin{matrix}
C^{-1/2} & 0\\
0 & C^{-1/2}
\end{matrix}\Bigg]\Bigg[\begin{matrix}
C & C_{n,k}\\
C_{n,k} & C
\end{matrix}\Bigg]
\Bigg[\begin{matrix}
C^{-1/2} & 0\\
0 & C^{-1/2}
\end{matrix}\Bigg]=\Bigg[\begin{matrix}
0 & -C^{-1/2}C_{n,k}C^{-1/2}\\
-C^{-1/2}C_{n,k}C^{-1/2} & 0
\end{matrix}\Bigg]\,.
%
\end{align*}
Note that $S_{n,k}$ is a rank-4 operator on $\H^2$ and direct calculations yield
\begin{align*}
&S_{n,k}[-(\eta_1+\eta_k),\eta_1+\eta_k]\trans=\varrho_n[-(\eta_1+\eta_k),\eta_1+\eta_k]\trans\,,\\
&S_{n,k}[-(\eta_1-\eta_k),\eta_1-\eta_k]\trans=\varrho_n[-(\eta_1-\eta_k),\eta_1-\eta_k]\trans\,,\\
&S_{n,k}[\eta_1+\eta_k,\eta_1+\eta_k]\trans=-\varrho_n[\eta_1+\eta_k,\eta_1+\eta_k]\trans\,,\\
&S_{n,k}[\eta_1-\eta_k,\eta_1-\eta_k]\trans=-\varrho_n[\eta_1-\eta_k,\eta_1-\eta_k]\trans\,.
\end{align*}
This implies that the non-zero eigenvalues of $S_{n,k}$ are $\alpha_{k,1}=\alpha_{k,2}=\varrho_n$, $\alpha_{k,3}=\alpha_{k,4}=-\varrho_n$ with corresponding eigenvectors
\begin{align*}
&\xi_{k,1}=[-(\eta_1+\eta_k),\eta_1+\eta_k]\trans/2\,,\qquad
\xi_{k,2}=[-(\eta_1-\eta_k),\eta_1-\eta_k]\trans/2\,,\\
&\xi_{k,3}=[\eta_1+\eta_k,\eta_1+\eta_k]\trans/2\,,\hspace{1.5cm}
\xi_{k,4}=[\eta_1-\eta_k,\eta_1-\eta_k]\trans/2\,.
\end{align*}
Therefore, recalling the null covariance operator $C_{F_0}=\diag(C,C)$ in \eqref{cf0} and observing $C_{F_0}^{-1/2}=\diag(C^{-1/2},C^{-1/2})$, we obtain that, for $(x,y)\in\H^2$,
\begin{align*}
\frac{\d\P_{F_{n,k}}(x,y)}{\d\P_{F_0}(x,y)}
&=(1-\varrho_n^2)^{-1}\exp\bigg[-\frac{1}{2}\Big\{\frac{\alpha_{k,1}}{1-\alpha_{k,1}}\l (x,y),C_{F_0}^{-1/2}\xi_{k,1}\r^2+\frac{\alpha_{k,2}}{1-\alpha_{k,2}}\l (x,y),C_{F_0}^{-1/2}\xi_{k,2}\r^2\\
&\hspace{3.5cm}+\frac{\alpha_{k,3}}{1-\alpha_{k,3}}\l (x,y),C_{F_0}^{-1/2}\xi_{k,3}\r^2+\frac{\alpha_{k,4}}{1-\alpha_{k,4}}\l (x,y),C_{F_0}^{-1/2}\xi_{k,4}\r^2\Big\}\bigg]\\
&=(1-\varrho_n^2)^{-1}\exp\bigg(-\frac{1}{8}\Big[\frac{\varrho_n}{1-\varrho_n}\{-\l x,C^{-1/2}(\eta_1+\eta_k)\r+\l y,C^{-1/2}(\eta_1+\eta_k)\r\}^2\\
&\hspace{4cm}+\frac{\varrho_n}{1-\varrho_n}\{-\l x,C^{-1/2}(\eta_1-\eta_k)\r+\l y,C^{-1/2}(\eta_1-\eta_k)\r\}^2\\
&\hspace{4cm}+\frac{-\varrho_n}{1+\varrho_n}(\l x,C^{-1/2}(\eta_1+\eta_k)\r+\l y,C^{-1/2}(\eta_1+\eta_k)\r)^2\\
&\hspace{4cm}+\frac{-\varrho_n}{1+\varrho_n}(\l x,C^{-1/2}(\eta_1-\eta_k)\r+\l y,C^{-1/2}(\eta_1-\eta_k)\r)^2\Big]\bigg)\\
%
%
&=(1-\varrho_n^2)^{-1}\exp\bigg\{-\frac{\varrho_n}{4(1-\varrho_n)}\big(\lambda_1^{-1}\l x-y,\eta_1\r^2+\lambda_k^{-1}\l x-y,\eta_k\r^2\big)\\
&\hspace{3cm}~~~+\frac{\varrho_n}{4(1+\varrho_n)}\big(\lambda_1^{-1}\l x-y,\eta_1\r^2+\lambda_k^{-1}\l x-y,\eta_k\r^2\big)\bigg\}\,.
\end{align*}
Recalling that $F_n$ is uniform measure over $\{F_{n,2},\ldots,F_{n,M}\}$ and that $(X_1,Y_1),\ldots,(X_n,Y_n)$ are i.i.d., we obtain from the above equation that
\begin{align*}
&\frac{\d\P_{F_{n}}(\T_n)}{\d\P_{F_0}(\T_n)}=(M-1)^{-1}(1-\varrho_n^2)^{-n}\sum_{k=2}^{M}\prod_{j=1}^n\exp\Big\{-\frac{\varrho_n}{4(1-\varrho_n)}\big(\lambda_1^{-1}\l X_j- Y_j,\eta_1\r^2+\lambda_k^{-1}\l X_j-Y_j,\eta_k\r^2\big)\\
&\hspace{6.5cm}~+\frac{\varrho_n}{4(1+\varrho_n)}\big(\lambda_1^{-1}\l X_j+Y_j,\eta_1\r^2+\lambda_k^{-1}\l X_j+Y_j,\eta_k\r^2\big)\Big\} \,.
\end{align*}
Therefore, in view of the above equation and the fact that by construction $C$ is the covariance operator of the mean-zero Gaussian measure ${\nu_0}\in\V$, we have
\begin{align*}
\E_{\P_{F_0}}\bigg\{\frac{\d\P_{F_n}(\T_n)}{\d\P_{F_0}(\T_n)}\bigg\}^2=J_{n,1}+J_{n,2}\,,
\end{align*}
where, by denoting $x_k=\l x,\eta_k\r$ and $y_k=\l y,\eta_k\r$ for $k\geq1$, the $J_{n,1}$ and $J_{n,2}$ is defined as
\begin{align}\label{jn1jn2}
J_{n,1}&=(M-1)^{-2}(1-\varrho_n^2)^{-2n}\sum_{2\leq k_1\neq k_2\leq M}\notag\\
&\quad\bigg(\int_{\H^2}\exp\Big[-\frac{\varrho_n}{4(1-\varrho_n)}\big\{2\lambda_1^{-1}(x_1-y_1)^2+\lambda_{k_1}^{-1}(x_{k_1}-y_{k_1})^2+\lambda_{k_2}^{-1}(x_{k_2}- y_{k_2})^2\big\}\notag\\
&\hspace{1cm}+\frac{\varrho_n}{4(1+\varrho_n)}\big\{2\lambda_1^{-1}(x_1+ y_1)^2+\lambda_{k_1}^{-1}(x_{k_1}+y_{k_1})^2+\lambda_{k_2}^{-1}( x_{k_2}+y_{k_2})^2\big\}\Big]\,{\nu_0}(\d x)\,{\nu_0}(\d y)\bigg)^n\notag\\
J_{n,2}&=(M-1)^{-2}(1-\varrho_n^2)^{-2n}\sum_{k=2}^{M}\bigg(\int_{\H^2}\exp\Big[-\frac{\varrho_n}{2(1-\varrho_n)}\big\{\lambda_1^{-1}(x_1-y_1)^2+\lambda_k^{-1}(x_k-y_k)^2\big\}\notag\\
&\hspace{1cm}+\frac{\varrho_n}{2(1+\varrho_n)}\big\{\lambda_1^{-1}(x_1+y_1)^2+\lambda_k^{-1}(x_k+y_k)^2\big\}\Big]\,{\nu_0}(\d x)\,{\nu_0}(\d y)\bigg)^n\,.
\end{align}

For the first term $J_{n,1}$, it follows from the definition of Gaussian measures and the definition of $\nu_0$ that $J_{n,1}$ in \eqref{jn1jn2} reduces to
\begin{align*}
&J_{n,1}=(M-1)^{-2}(1-\varrho_n^2)^{-2n}\sum_{2\leq k_1\neq k_2\leq M}\notag\\
&\hspace{1.5cm}\bigg(\int_{\mathbb R^6}\exp\Big[-\frac{\varrho_n}{4(1-\varrho_n)}\big\{2\lambda_1^{-1}(x_1-y_1)^2+\lambda_{k_1}^{-1}(x_{k_1}-y_{k_1})^2+\lambda_{k_2}^{-1}(x_{k_2}- y_{k_2})^2\big\}\notag\\
&\hspace{2.5cm}+\frac{\varrho_n}{4(1+\varrho_n)}\big\{2\lambda_1^{-1}(x_1+ y_1)^2+\lambda_{k_1}^{-1}(x_{k_1}+y_{k_1})^2+\lambda_{k_2}^{-1}( x_{k_2}+y_{k_2})^2\big\}\Big]\notag\\
&\hspace{3cm}\N_{\lambda_1}(\d x_1)\,\N_{\lambda_{k_1}}(\d x_{k_1})\,\N_{\lambda_{k_2}}(\d x_{k_2})\,\N_{\lambda_1}(\d y_1)\,\N_{\lambda_{k_1}}(\d y_{k_1})\,\N_{\lambda_{k_2}}(\d y_{k_2})\bigg)^n\,.
\end{align*}
Furthermore, an application of change of variables and the above equation implies
\begin{align}\label{jn1}
&J_{n,1}=(M-1)^{-2}(1-\varrho_n^2)^{-2n}\sum_{2\leq k_1\neq k_2\leq M}\notag\\
&\hspace{1.5cm}\bigg(\int_{\mathbb R^6}\exp\Big[-\frac{\varrho_n}{4(1-\varrho_n)}\big\{2 (s_1-t_1)^2+ (s_{k_1}-t_{k_1})^2+ (s_{k_2}- t_{k_2})^2\big\}\notag\\
&\hspace{2.5cm}+\frac{\varrho_n}{4(1+\varrho_n)}\big\{2 (s_1+ t_1)^2+ (s_{k_1}+ t_{k_1})^2+ ( s_{k_2}+t_{k_2})^2\big\}\Big]\notag\\
&\hspace{3cm}\N_{1}(\d s_1)\,\N_{1}(\d s_{k_1})\,\N_{1}(\d s_{k_2})\,\N_{1}(\d t_1)\,\N_{1}(\d t_{k_1})\,\N_{1}(\d t_{k_2})\bigg)^n\notag\\
&=(M-1)^{-2}(1-\varrho_n^2)^{-2n}\sum_{2\leq k_1\neq k_2\leq M}\bigg\{\int_{\mathbb R^6}\exp\big([s_1,t_1]\Sigma_1[s_1,t_1]\trans+2^{-1}[s_{k_1},t_{k_1}]\Sigma_1[s_{k_1},t_{k_1}]\trans\notag\\
&+2^{-1}[s_{k_2},t_{k_2}]\Sigma_1[s_{k_2},t_{k_2}]\trans\big)\,\N_{1}(\d s_1)\,\N_{1}(\d s_{k_1})\,\N_{1}(\d s_{k_2})\,\N_{1}(\d t_1)\,\N_{1}(\d t_{k_1})\,\N_{1}(\d t_{k_2})\bigg\}^n,
\end{align}
where
\begin{align}\label{sigma}
\Sigma_1=\Bigg[\begin{matrix}
-\frac{\varrho_n^2}{1-\varrho_n^2} & \frac{\varrho_n}{1-\varrho_n^2}\\
\frac{\varrho_n}{1-\varrho_n^2}&-\frac{\varrho_n^2}{1-\varrho_n^2} \\
\end{matrix}\Bigg]\,.
\end{align}
Observe the fact that $\det(I_2-\Sigma_1)=(1-\varrho_n^2)^{-1}$ and $\det(I_2-2\Sigma_1)=1$. In addition, by applying the formula for the moment generating functions of the Gaussian quadratic forms, we obtain that, for invertible matrix $\Sigma$ and $U\in\mathbb R^p$ that follows a standard multivariate normal distribution $N(0,I_p)$, it holds that
\begin{align}\label{formula}
\E\exp(U\trans\Sigma U)=\det(I_p-2\Sigma)^{-1/2}\,.
\end{align}
We therefore deduce from \eqref{jn1} and \eqref{formula} that
\begin{align}\label{jn12}
&J_{n,1}=(M-2)(M-1)^{-1}(1-\varrho_n^2)^{-2n}\det(I_2-2\Sigma_1)^{-n/2}\det(I_2-\Sigma_1)^{-n}\notag\\
&=(M-2)(M-1)^{-1}(1-\varrho_n^2)^{-n}\leq(1-\varrho_n^2)^{-n}= \exp(n\varrho_n^2)\{1+o(1)\}\leq\exp(2n\varrho_n^2)\,.
\end{align}

For the second term $J_{n,2}$ in \eqref{jn1jn2}, we have
\begin{align*}
J_{n,2}&=(M-1)^{-2}(1-\varrho_n^2)^{-2n}\sum_{k=2}^{M}\bigg(\int_{\mathbb R^4}\exp\Big[-\frac{\varrho_n}{2(1-\varrho_n)}\big\{\lambda_1^{-1}(x_1-y_1)^2+\lambda_k^{-1}(x_k-y_k)^2\big\}\\
&\quad+\frac{\varrho_n}{2(1+\varrho_n)}\big\{\lambda_1^{-1}(x_1+y_1)^2+\lambda_k^{-1}(x_k+y_k)^2\big\}\Big]\,\N_{ \lambda_1}(\d x_1)\,\N_{ \lambda_k}(\d x_k)\,\N_{ \lambda_1}(\d y_1)\,\N_{ \lambda_k}(\d y_k)\bigg)^n\\
&=(M-1)^{-1}(1-\varrho_n^2)^{-2n}\bigg[\int_{\mathbb R^2}\exp\Big\{-\frac{\varrho_n(s-t)^2}{2(1-\varrho_n)}+\frac{\varrho_n(s+t)^2}{2(1+\varrho_n)}\Big\}\,\N_{1}(\d s)\,\N_{1}(\d t)\bigg]^{2n}\\
&=(M-1)^{-1}(1-\varrho_n^2)^{-2n}\bigg\{\int_{\mathbb R^2}\exp([s,t]\Sigma_1[s,t]\trans)\,\N_{1}(\d s)\,\N_{1}(\d t)\bigg\}^{2n}\,,
\end{align*}
where $\Sigma_1$ is defined in \eqref{sigma}. By applying \eqref{formula}, we obtain
\begin{align}\label{jn22}
J_{n,2}&=(M-1)^{-1}(1-\varrho_n^2)^{-2n}\det(I_2-2\Sigma_1)^{-n/2}=(M-1)^{-1}(1-\varrho_n^2)^{-2n}\notag\\
&=(M-1)^{-1}\exp(2n\varrho_n^2)\{1+o(1)\}\leq(M-1)^{-1}\exp(4n\varrho_n^2)\,.
\end{align}
By taking $M=\lceil 4^{-1}(1-\alpha-\beta)^{-2}\rceil+2>4^{-1}(1-\alpha-\beta)^{-2}$ and $\varrho_n=(cn)^{-1/2}$, where 
\begin{align*}
c=2^{-1}\log\big\{2^{-1}\sqrt{M^2+4M+4(1-\alpha-\beta)^2M}-M/2\big\}>0\,,
\end{align*}
we obtain by combining \eqref{jn12} and \eqref{jn22} that, 
\begin{align*}
\E_{\P_{F_0}}\bigg\{\frac{\d\P_{F_n}(\T_n)}{\d\P_{F_0}(\T_n)}\bigg\}^2\leq \exp(2n\varrho_n^2)+(M-1)^{-1}\exp(4n\varrho_n^2)\leq 1+4(1-\alpha-\beta)^2\,.
\end{align*}
Now, we have shown \eqref{goal} and completes the proof.

\subsection{Proof of Proposition~\ref{prop:edhat}}\label{app:prop:edhat}

For \eqref{eed}, observe that
\begin{align*}
&|\hat\phi_X(w)-\hat\phi_Y(w)|^2\\
&=\bigg\{\frac{1}{n_1}\sum_{j=1}^{n_1}\exp(\i\l X_j,w\r)-\frac{1}{n_2}\sum_{j=1}^{n_2}\exp(\i\l Y_j,w\r)\bigg\}
\bigg\{\frac{1}{n_1}\sum_{j=1}^{n_1}\exp(-\i\l X_j,w\r)-\frac{1}{n_2}\sum_{j=1}^{n_2}\exp(-\i\l Y_j,w\r)\bigg\}\\
&=\frac{1}{n_1^2}\sum_{j,k=1}^{n_1}\exp(\i\l X_{j}-X_k,w\r)+\frac{1}{n_2^2}\sum_{j,k=1}^{n_2}\exp(\i\l Y_j-Y_k,w\r)\\
&\qquad-\frac{1}{n_1n_2}\sum_{j=1}^{n_1}\sum_{k=1}^{n_2}\{\exp(\i\l X_j-Y_k,w\r)+\exp(\i\l -X_j+Y_k,w\r)\}\,,
\end{align*}
so that \eqref{eed} follows.
For \eqref{new}, observe the fact that
\begin{align*}
\sum_{j,k=1}^{n_1+n_2}\phi_\nu(Z_j-Z_k)=\sum_{j,k=1}^{n_1}\phi_\nu(Z_j-Z_k)+\sum_{j,k=1}^{n_2}\phi_\nu(Z_{n_1+j}-Z_{n_1+k})+2\sum_{j=1}^{n_1}\sum_{k=1}^{n_2}\phi_\nu(Z_j-Z_k)\,,
\end{align*}
so that in view of \eqref{eed}, we conclude that
\begin{align*}
&\hat\ed_\nu^2(X,Y)=\frac{1}{n_1^2}\sum_{j,k=1}^{n_1}\phi_{\nu}( X_{j}-X_k)+\frac{1}{n_2^2}\sum_{j,k=1}^{n_2}\phi_{\nu}( Y_j-Y_k)\\
&\qquad\qquad\qquad-\frac{1}{n_1n_2}\bigg\{\sum_{j,k=1}^{n_1+n_2}\phi_\nu(Z_j-Z_k)-\sum_{j,k=1}^{n_1}\phi_\nu(Z_j-Z_k)-\sum_{j,k=1}^{n_2}\phi_\nu(Z_{n_1+j}-Z_{n_1+k})\bigg\}\\
&=\frac{n_1+n_2}{n_1n_2}\bigg\{\frac{1}{n_1}\sum_{j,k=1}^{n_1}\phi_{\nu}(Z_j-Z_k)+\frac{1}{n_2}\sum_{j,k=1}^{n_2}\phi_{\nu}(Z_{n_1+j}-Z_{n_1+k})-\frac{1}{n_1+n_2}\sum_{j,k=1}^{n_1+n_2}\phi_{\nu}(Z_j-Z_k)\bigg\}\,.
\end{align*}

\subsection{Proof of Proposition~\ref{thm:converge:ed}}\label{app:thm:converge:ed}

Note that (i) follows from the strong law of large numbers for V-statistics and the fact that $\phi_\nu$ is uniformly bounded. For (ii), recall the definition of $(Z_1,\ldots,Z_{n_1},Z_{n_1+1},\ldots,Z_{n_1+n_2})=(X_1,\ldots,X_{n_1},Y_1,\ldots,Y_{n_2})$ in Proposition~\ref{prop:edhat}.
When $F_X=F_Y=F_Z$,
for $w_1,w_2\in\H$, define
\begin{align*}
h(w_1,w_2)=\phi_{\nu}(w_1-w_2)-\E\phi_{\nu}(w_1-Z)-\E\phi_{\nu}(Z-w_2)+\E\phi_{\nu}(Z-Z')\,.
\end{align*}
We obtain from \eqref{eed} that
\begin{align*}
\hat\ed_\nu^2(X,Y)
&=\frac{1}{n_1^2}\sum_{j,k=1}^{n_1}\{h(X_{j},X_k)+\E_Z\phi_{\nu}(X_j-Z)+\E_Z\phi_{\nu}(Z-X_k)-\E\phi_{\nu}(Z-Z')\}\\
&\quad+\frac{1}{n_2^2}\sum_{j,k=1}^{n_2}\{h(Y_{j},Y_k)+\E_Z\phi_\nu(Y_j-Z)+\E_Z\phi_\nu(Z-Y_k)-\E\phi_{\nu}(Z-Z')\}\\
&\quad-\frac{2}{n_1n_2}\sum_{j=1}^{n_1}\sum_{k=1}^{n_2}\{h(X_j,Y_k)+\E_Z\phi_\nu(X_j-Z)+\E_Z\phi_\nu(Z-Y_k)-\E\phi_{\nu}(Z-Z')\}\\
&=\frac{1}{n_1^2}\sum_{j,k=1}^{n_1}h(X_{j},X_k)+\frac{1}{n_2^2}\sum_{j,k=1}^{n_2}h(Y_{j},Y_k)-\frac{2}{n_1n_2}\sum_{j=1}^{n_1}\sum_{k=1}^{n_2}h(X_j,Y_k)\,.
\end{align*}
Observe that $Z_j=X_j$ for $1\leq j\leq n_1$ and $Z_{k+n_1}=Y_k$ for $1\leq k\leq n_2$. Therefore, we deduce from the above equation that
\begin{align}\label{te}
n_1\hat\ed_\nu^2(X,Y)&=\frac{n_1+n_2}{n_2}(\tilde V_1+\tilde V_2-\tilde V_3)\,,
\end{align}
where
%
\begin{align*}
&\tilde V_1=\frac{1}{n_1}\sum_{(j,k)\in I_2(n_1)}h(Z_{j},Z_k)\,;
\\
&\tilde V_2=\frac{1}{n_2}\sum_{(j,k)\in I_2(n_2)}h(Z_{n_1+j},Z_{n_1+k})\,;
\\
&\tilde V_3=\frac{1}{n_1+n_2}\sum_{(j,k)\in I_2(n_1+n_2)}h(Z_{j},Z_k)\,;\\
&R=\frac{1}{n_1}\sum_{j=1}^{n_1}h(Z_j,Z_j)+\frac{1}{n_2}\sum_{j=1}^{n_2}h(Z_{n_1+j},Z_{n_1+j})-\frac{1}{n_1+n_2}\sum_{j=1}^{n_1+n_2}h(Z_{j},Z_j)\,.
\end{align*}
Note that $h(\cdot,Z)=\phi_{\nu}(\cdot-Z)-\E\phi_{\nu}(\cdot,Z)$, so that $\E h(\cdot,Z)\equiv0$.
In addition, in view of (ii) in Theorem~\ref{thm:dvar}, $\var\{h(Z,Z')\}=\E h^2(Z,Z')=\dvar_\nu^2(Z)>0$.
This implies that $\tilde V_1,\tilde V_2,\tilde V_3$ are re-scaled V-statistics of rank-1. Let $\{\psi_{\nu,\ell}\}_{\ell\geq1}$  denote the corresponding eigenvalues of $\{\tau_{\nu,\ell}\}_{\ell\geq1}$ of the operator $R_\nu$. Following the proof of Theorem 5.5.2 in \cite{serfling1980}, we obtain, for $z_1,z_2\in\H$,
\begin{align*}
h(z_1,z_2)=\sum_{\ell=1}^\infty\tau_{\nu,\ell}\,\psi_{\nu,\ell}(z_1)\psi_{\nu,\ell}(z_2)\,.
\end{align*}
For $\ell\geq1$, define
\begin{align*}
&W_{\ell,n_1}=\frac{1}{\sqrt{n_1}}\sum_{j=1}^{n_1}\psi_{\nu,\ell}(Z_j)\,;\qquad W_{\ell,n_2}=\frac{1}{\sqrt{n_2}}\sum_{j=n_1+1}^{n_1+n_2}\psi_{\nu,\ell}(Z_j)\,;\\
&T_{\ell,n_1}=\frac{1}{n_1}\sum_{j=1}^{n_1}\psi_{\nu,\ell}^2(Z_j)\,;\qquad~~~~ T_{\ell,n_2}=\frac{1}{n_2}\sum_{j=n_1+1}^{n_1+n_2}\psi_{\nu,\ell}^2(Z_j)\,,
\end{align*}
so that
\begin{align*} &\tilde V_1=\sum_{\ell=1}^\infty\tau_{\nu,\ell} (W_{\ell,n_1}^2-T_{\ell,n_1})\,,\qquad \tilde V_2=\sum_{\ell=1}^\infty\tau_{\nu,\ell} (W_{\ell,n_2}^2-T_{\ell,n_2})\,,\\
&\tilde V_3=\sum_{\ell=1}^\infty\tau_{\nu,\ell}\bigg\{\Big(\frac{\sqrt{n_1}W_{\ell,n_1}+\sqrt{n_2}W_{\ell,n_2}}{\sqrt{n_1+n_2}}\Big)^2-\frac{{n_1}T_{\ell,n_1}+{n_2}T_{\ell,n_2}}{{n_1+n_2}}\bigg\}\,.
\end{align*}
Examining the proof of Theorem 5.5.2 in \cite{serfling1980}, we obtain that
\begin{align*}
\tilde V_1\converged\sum_{\ell=1}^\infty\tau_{\nu,\ell}(\xi_\ell^2-1)\,;\quad \tilde V_2\converged\sum_{\ell=1}^\infty\tau_{\nu,\ell}(\zeta_\ell^2-1)\,;\quad \tilde V_3\converged\sum_{\ell=1}^\infty\tau_{\nu,\ell}\{(1+c)^{-1}(\xi_\ell+\sqrt{c}\zeta_\ell)^2-1\}\,.
\end{align*}
In addition, by the law of large numbers, since $h$ is bounded, we have 
\begin{align*}
R\to \E h(Z,Z')=-\E\phi_{\nu}(Z-Z')\qquad\as
\end{align*}
The result therefore follows in view of \eqref{te}.

\subsection{Proof of Theorem~\ref{thm:upper:two}}\label{app:thm:upper:two}

\subsubsection*{Upper bound for the uniform separation rate of $\hat\psi_{\nu,n,2}$}

The proof follows similar ideas as in the proof of Theorem~\ref{thm:upper:ind}.
Let $\Pi$ denote a generic random permutation of $(1,2,\ldots,n )$.
Define
\begin{align*}
\hat F_{\nu,B}(q)=\frac{1}{B}\sum_{b=1}^B\one\{\ed_\nu^2(\hat{\mathbb P}^X,\hat{\mathbb P}^Y,\Pi_b)\leq q\}\,,\qquad \tilde F_{\nu}(q)=|I_n(n)|^{-1}\sum_{\Pi\in I_n(n)}\one\{\ed_\nu^2(\hat{\mathbb P}^X,\hat{\mathbb P}^Y)\leq q\}\,.
\end{align*}
Observe that $\hat q_{\nu,1-\alpha}=\inf\big\{q\in\mathbb R:\hat F_{\nu,B}(q)\geq 1-\alpha\big\}$. Let $\tilde q_{\nu,1-\alpha}=\inf\big\{q\in\mathbb R:\tilde F_{\nu}(q)\geq 1-\alpha\big\}$, and let $\tilde w_{\nu,\gamma}$ denote the $\gamma$-quantile of $\tilde q_{\nu,1-\alpha/2}$.
Following similar derivation as in the proof of Theorem~\ref{thm:upper:ind}, we obtain that, for $\mathcal F_\nu(\rho)$ defined in \eqref{fnurhoa}, and for any $\nu\in\V$ and $F\in\mathcal F_\nu(\rho)$,
\begin{align}\label{q4}
\P_F\{\hat\psi_{\nu,n}(\alpha)=0\}&\leq\P_F\bigg\{\sup_{q\in\mathbb R}|\hat F_{\nu,B}(q)-\tilde F_{\nu}(q)|>\sqrt{(2B)^{-1}\log(6/\beta)}\bigg\}\notag\\
&\quad+\P_F(\tilde q_{\nu,1-\alpha/2}> \tilde w_{\nu,1-\beta/3})+\P_F\{\ed_\nu^2(\hat{\mathbb P}^X,\hat{\mathbb P}^Y)\leq \tilde w_{\nu,1-\beta/3}\}\notag\\
&\leq2\beta/3+\P_F\{\ed_\nu^2(\hat{\mathbb P}^X,\hat{\mathbb P}^Y)\leq \tilde w_{\nu,1-\beta/3}\}\,,
\end{align}
provided that $B\geq2\alpha^{-2}\log(6/\beta)$. Let $\kappa_{\nu,1-\beta/3}$ denote the $(1-\beta/3)$-quantile of $|\ed_\nu^2(\hat{\mathbb P}^X,\hat{\mathbb P}^Y)-\E_F\{\ed_\nu^2(\hat{\mathbb P}^X,\hat{\mathbb P}^Y)\}|$. 
Write for simplicity $\|\phi_1-\phi_2\|_\nu^2=\int_{\H}|\phi_1-\phi_2|^2\d \nu$. Suppose that
\begin{align*}
\|\phi_{1}-\phi_2\|_\nu^2&\geq \sup_{\nu\in\V}\sup_{F\in\mathcal P_1}\tilde w_{\nu,1-\beta/3}+\sup_{\nu\in\V}\sup_{F\in\mathcal P_1}\kappa_{\nu,1-\beta/3}+\sup_{\nu\in\V}\sup_{F\in\mathcal P_1}|\|\phi_{1}-\phi_2\|_\nu^2-\E_F\{\ed_\nu^2(\hat{\mathbb P}^X,\hat{\mathbb P}^Y)\}|\,,
\end{align*}
then, we obtain, for any $\nu\in\V$ and $F\in\mathcal F_\nu$,
\begin{align*}
&\P_F\{\ed_\nu^2(\hat{\mathbb P}^X,\hat{\mathbb P}^Y)\leq \tilde w_{\nu,1-\beta/3}\}\\
&\leq\P_F\big[\ed_\nu^2(\hat{\mathbb P}^X,\hat{\mathbb P}^Y)\leq \|\phi_{1}-\phi_2\|_\nu^2-\kappa_{\nu,1-\beta/3}-|\|\phi_{1}-\phi_2\|_\nu^2-\E_F\{\ed_\nu^2(\hat{\mathbb P}^X,\hat{\mathbb P}^Y)\}|\big]\notag\\
&\leq\P_F\big[\ed_\nu^2(\hat{\mathbb P}^X,\hat{\mathbb P}^Y)\leq \E_F\{\ed_\nu^2(\hat{\mathbb P}^X,\hat{\mathbb P}^Y)\}-\kappa_{\nu,1-\beta/3}\big]\leq\beta/3\,.
\end{align*}
Combining the above equation with \eqref{q3} yields
\begin{align*}
\rho\{\hat\psi_{\nu,n}(\alpha),\beta\}&\leq\sup_{\nu\in\V}\sup_{F\in\mathcal P_1}\tilde w_{\nu,1-\beta/3}+\sup_{\nu\in\V}\sup_{F\in\mathcal P_1}\kappa_{\nu,1-\beta/3}+\sup_{\nu\in\V}\sup_{F\in\mathcal P_1}|\|\phi_{1}-\phi_2\|_\nu^2-\E_F\{\ed_\nu^2(\hat{\mathbb P}^X,\hat{\mathbb P}^Y)\}|\,.
\end{align*}
It then suffices to bound three terms in the right-hand side of the above equation.

For $\tilde w_{\nu,1-\beta/3}$, we first derive an equivalent formula for $\ed_\nu^2(\hat{\mathbb P}^X,\hat{\mathbb P}^Y)$ in \eqref{new} that is easier to handle. From the proof of Proposition~\ref{prop:edhat} in Section~\ref{app:prop:edhat}, we obtain
\begin{align*}
\ed_\nu^2(\hat{\mathbb P}^X,\hat{\mathbb P}^Y,\Pi)&=\int_\H\bigg\{\frac{1}{n_1}\sum_{j_1=1}^{n_1}\exp(\i\l Z_{\Pi(j_1)},w\r)-\frac{1}{n_2}\sum_{k_1=1}^{n_2}\exp(\i\l Z_{\Pi(n_1+k_1)},w\r)\bigg\}\\
&\qquad\times\bigg\{\frac{1}{n_1}\sum_{j_2=1}^{n_1}\exp(-\i\l Z_{\Pi(j_2)},w\r)-\frac{1}{n_2}\sum_{k_2=1}^{n_2}\exp(-\i\l Z_{\Pi(n_1+k_2)},w\r)\bigg\}\,\nu(\d w)\\
%
%
&=\frac{1}{n_1^2n_2^2}\sum_{j_1,j_2=1}^{n_1}\sum_{k_1,k_2=1}^{n_2}\omega_\nu(Z_{\Pi(j_1)},Z_{\Pi(j_2)},Z_{\Pi(n_1+k_1)},Z_{\Pi(n_1+k_2)})\,,
\end{align*}
where $\omega_\nu$ is defined in \eqref{omeganu}. 
Without loss of generality assume that $n_1\leq n_2$. Suppose $\K_1,\ldots,\K_{n_1}$ are i.i.d.~and follow uniform distribution on $I_{n_2}(n_2)$. Write $\K=(\K_1,\ldots,\K_{n_1})$ and following similar arguments as in the proof of Theorem~\ref{thm:upper:ind},
we define
\begin{align*}
&\ed_\nu^2(\hat{\mathbb P}^X,\hat{\mathbb P}^Y,\Pi,\K)=\frac{1}{n_1^2}\sum_{r,\ell=1}^{n_1}\omega_\nu(Z_{\Pi(r)},Z_{\Pi(\ell)},Z_{\Pi(n_1+\K_{r})},Z_{\Pi(n_1+\K_{\ell})})\,.
\end{align*}
Observe that
\begin{align}\label{p32}
\ed_\nu^2(\hat{\mathbb P}^X,\hat{\mathbb P}^Y,\Pi)=\E_{\K}\{\ed_\nu^2(\hat{\mathbb P}^X,\hat{\mathbb P}^Y,\Pi,\K)|\T_{n,2},\Pi\}
\end{align}
Since $\Pi$ follows a uniform distribution that are uniformly distributed over the set of all permutations of $(1,2,\ldots,n_1+n_2)$, then, for any $r=1,2,\ldots,n_1$, exchanging $Z_{\Pi(r)}$ and $Z_{\Pi(n_1+\K_{r})}$  does not alter the distribution of $\ed_\nu^2(\hat{\mathbb P}^X,\hat{\mathbb P}^Y,\Pi,\K)$. To be specific, let $\delta_1,\ldots,\delta_{n_1}$ be i.i.d.~Bernoulli random variables such that $\P(\delta_j=1)=\P(\delta_j=0)=1/2$. Define
\begin{align*}
&\ed_\nu^2(\hat{\mathbb P}^X,\hat{\mathbb P}^Y,\Pi,\K,\delta)=\frac{1}{n_1^2}\sum_{r,\ell=1}^{n_1}\omega_\nu(W_r^{(\Pi,\K,\delta)},W_\ell^{(\Pi,\K,\delta)},W_{n_1+r}^{(\Pi,\K,\delta)},W_{n_1+\ell}^{(\Pi,\K,\delta)})\,,
\end{align*}
where, for $1\leq r\leq n_1$,
\begin{align*}
W_r^{(\Pi,\K,\delta)}=(1-\delta_r) Z_{\Pi(r)}+\delta_r Z_{\Pi(n_1+\K_r)}\,,\qquad W_{r+n_1}^{(\Pi,\K,\delta)}=\delta_r Z_{\Pi(r)}+(1-\delta_r) Z_{\Pi(n_1+\K_r)}\,.
\end{align*}
Then, $\ed_\nu^2(\hat{\mathbb P}^X,\hat{\mathbb P}^Y,\Pi,\K)$ and $\ed_\nu^2(\hat{\mathbb P}^X,\hat{\mathbb P}^Y,\Pi,\K,\delta)$ are identically distributed.
Furthermore, introducing i.i.d.~Rademacher random variables $\xi_1,\ldots,\xi_{n_1}$, that is, $\P(\xi_j=1)=\P(\xi_j=-1)=1/2$, let $\xi=(\xi_1,\ldots,\xi_{n_1})\in\mathbb R^{n_1}$ and define
\begin{align*}
&\ed_\nu^2(\hat{\mathbb P}^X,\hat{\mathbb P}^Y,\Pi,\K,\xi)=\frac{1}{n_1^2}\sum_{r,\ell=1}^{n_1}\xi_r\,\xi_\ell\,\omega_\nu(Z_{\Pi(r)},Z_{\Pi(\ell)},Z_{\Pi(n_1+\K_{r})},Z_{\Pi(n_1+\K_{\ell})})\,.
\end{align*}
Then, in view of the property of $\omega_\nu$ in \eqref{omega1}, we deduce that $\ed_\nu^2(\hat{\mathbb P}^X,\hat{\mathbb P}^Y,\Pi,\K,\delta)$ and $\ed_\nu^2(\hat{\mathbb P}^X,\hat{\mathbb P}^Y,\Pi,\K,\xi)$ have the same distribution, so that
\begin{align}\label{p82}
\ed_\nu^2(\hat{\mathbb P}^X,\hat{\mathbb P}^Y,\Pi,\K)\text{ and }\ed_\nu^2(\hat{\mathbb P}^X,\hat{\mathbb P}^Y,\Pi,\K,\xi) \text{ are identically distributed.}
\end{align}
Therefore, following the same derivation in \eqref{p9}, we obtain that, for any $t>0$,
\begin{align}\label{p92}
\P_\Pi\{\ed_\nu^2(\hat{\mathbb P}^X,\hat{\mathbb P}^Y,\Pi)\geq q|\T_{n,2}\}&\leq e^{-tq}\,\E_{\Pi,\K,\xi}\big[\exp\{t\ed_\nu^2(\hat{\mathbb P}^X,\hat{\mathbb P}^Y,\Pi,\K,\xi)\}|\T_{n,2}\big]\,.
\end{align}

For $1\leq r,\ell\leq n_1$, let
\begin{align*}
a_{r,\ell}(\T_{n,2},\Pi,\K)&=n_1^{-2}\omega_\nu(Z_{\Pi(r)},Z_{\Pi(\ell)},Z_{\Pi(n_1+\K_{r})},Z_{\Pi(n_1+\K_{\ell})})\,.
\end{align*}
Let $A_{\T_{n,2},\Pi,\K}$ be the $n\times n$ matrix with entries $a_{r,\ell}(\T_{n,2},\Pi,\K)$. Now, we have $\hat D_{\nu}(\T_{n,2}^{(\Pi)},\K,\xi)=\xi'A_{\T_{n,2},\Pi,\K}\xi$.
By the Hanson--Wright inequality (e.g.~Theorem~1.1 of \cite{Rudelson2013}), there exist constants $c_1,c_2>0$ such that, for $0\leq t\leq c_2\|A_{\T_{n,2},\Pi,\K}\|^{-1}$,
\begin{align}\label{p9.52}
&\E_{\Pi,\K,\xi}\big[\exp\{t\ed_\nu^2(\hat{\mathbb P}^X,\hat{\mathbb P}^Y,\Pi,\K,\xi)|\T_{n,2}\big]\leq \E_{\Pi,\K}\big\{\exp(c_1t\|A_{\T_{n,2},\Pi,\K}\|_{\rm F}^2)\,|\T_{n,2}\big\}\,.
\end{align}
Observing the fact that $\|\phi_{\nu}\|_\infty\leq1$ for any $\nu\in\V$, by a line-by-line check of the proof (Steps 2--4) of Theorem~1.1 of \cite{Rudelson2013}, it holds that the constants $c_1,c_2$ in equation \eqref{p9.52} is uniform in $\nu\in\V$.
Then, combining \eqref{p92} and \eqref{p9.52}, minimizing with respect to $0\leq t\leq c_2\|A_{\T_{n,2},\Pi,\K}\|^{-1}$ yields
\begin{align}\label{p112}
\P_\Pi\{\ed_\nu^2(\hat{\mathbb P}^X,\hat{\mathbb P}^Y,\Pi)\geq q|\T_{n,2}\}&\leq \E_{\Pi,\K}\bigg[\exp\bigg\{-c_3\min\bigg(\frac{q^2}{\|A_{\T_{n,2},\Pi,\K}\|_{\rm F}^2},\frac{q}{\|A_{\T_{n,2},\Pi,\K}\|}\bigg)\bigg\}\,\Big|\T_{n,2}\bigg]\notag\\
&\hspace{-1cm}\leq \E_{\Pi,\K}\bigg[\exp\bigg\{-c_3\min\bigg(\frac{q^2}{\|A_{\T_{n,2},\Pi,\K}\|_{\rm F}^2},\frac{q}{\|A_{\T_{n,2},\Pi,\K}\|_{\rm F}}\bigg)\bigg\}\,\Big|\T_{n,2}\bigg]\,.
\end{align}
for some constant $c_3>0$ that is uniform in $\nu\in\V$, where in the last step we used the fact that $\|A_{\T_{n,2},\Pi,\K}\|\leq\|A_{\T_{n,2},\Pi,\K}\|_{\rm F}$.
Therefore, we obtain that for any $\Pi$ and $\K$,
\begin{align}\label{p122}
&\|A_{\T_{n,2},\Pi,\K}\|_{\rm F}^2=\frac{1}{n^4}\sum_{r,\ell=1}^{n_1}\omega_\nu^2(Z_{\Pi(r)},Z_{\Pi(\ell)},Z_{\Pi(n_1+\K_{r})},Z_{\Pi(n_1+\K_{\ell})})\notag\\
&=\frac{1}{n^4}\sum_{r,\ell=1}^{n_1}\{\phi_{\nu}(Z_{\Pi(r)}-Z_{\Pi(\ell)})+\phi_{\nu}(Z_{\Pi(n_1+\K_{r})}-Z_{\Pi(n_1+\K_{\ell})})\notag\\
&\quad-\phi_{\nu}(Z_{\Pi(r)}-Z_{\Pi(n_1+\K_{r})})-\phi_{\nu}(Z_{\Pi(\ell)}-Z_{\Pi(n_1+\K_{\ell})})\}^2\leq cn^{-4}\sum_{r,\ell=1}^n\{\phi_{\nu}(Z_r-Z_\ell)\}^2\,.
\end{align}
Combining \eqref{p112} and \eqref{p122} yields
\begin{align*}
\P_\Pi\{\ed_\nu^2(\hat{\mathbb P}^X,\hat{\mathbb P}^Y,\Pi)\geq q|\T_{n,2}\}&\leq \exp\bigg(-c_3qn^4\sum_{r,\ell=1}^n\{\phi_{\nu}(Z_r-Z_\ell)\}^2\bigg)\,.
\end{align*}
Recalling that $\tilde q_{\nu,1-\alpha/2}$ is the $(1-\alpha/2)$-quantile of $\ed_\nu^2(\hat{\mathbb P}^X,\hat{\mathbb P}^Y,\Pi)$, we therefore obtain from the above equation that
\begin{align*}
\tilde q_{\nu,1-\alpha/2}\leq c\log(\alpha^{-1})n^{-2}\bigg[\sum_{r,\ell=1}^n\{\phi_{\nu}(Z_r-Z_\ell)\}^2\bigg]^{1/2}\,.
\end{align*}
Therefore, by Markov's inequality and the fact that $\|\phi_{\nu}\|_\infty\leq1$, we obtain from the above equation
\begin{align*}
\P(\tilde q_{\nu,1-\alpha/2}\geq q)&\leq ce^{-nq}\alpha^{-1}\E\bigg\{\exp\bigg(n^{-1}\bigg[\sum_{r,\ell=1}^n\{\phi_{\nu}(Z_r-Z_\ell)\}^2\bigg]^{1/2}\bigg)\bigg\}\leq ce^{-nq}\alpha^{-1}\,.
\end{align*}
Recalling that $\tilde w_{\nu,1-\beta/3}$ is the $(1-\beta/3)$-quantile of $\tilde q_{\nu,1-\alpha/2}$, the above equation implies that
\begin{align}\label{p302}
\tilde w_{\nu,1-\beta/3}\leq cn^{-1}\{\log(\alpha^{-1})+\log(\beta^{-1})\}\,,
\end{align}
where the constant $c>0$ in the above equation is uniform in $\nu\in\V$ and $F\in\mathcal P_1$.

\textbf{Step 2.} Bounding the quantile $\kappa_{\nu,1-\beta/3}$.

We start by re-writing $\ed_\nu^2(\hat{\mathbb P}^X,\hat{\mathbb P}^Y)$ to obtain a equivalent formula that is easier to handle. Direct calculations yields $\ed_\nu^2(\hat{\mathbb P}^X,\hat{\mathbb P}^Y)-\E_F\{\ed_\nu^2(\hat{\mathbb P}^X,\hat{\mathbb P}^Y)\}=J_{\nu,1}+J_{\nu,2}+J_{\nu,3}$, where
\begin{align*}
J_{\nu,1}&=\frac{(n_1+n_2)(n_1-1)}{n_1n_2}\big\{\tilde D_{\nu,1}(\T_{n,2})-\E_F\tilde D_{\nu,1}(\T_{n,2})\big\}\,,\\
J_{\nu,2}&=\frac{(n_1+n_2)(n_2-1)}{n_1n_2}\big\{\tilde D_{\nu,2}(\T_{n,2})-\E_F\tilde D_{\nu,2}(\T_{n,2})\big\}\,,\\
J_{\nu,3}&=-\frac{(n_1+n_2)(n_1+n_2-1)}{n_1n_2}\big\{\tilde D_{\nu,3}(\T_{n,2})-\E_F\tilde D_{\nu,3}(\T_{n,2})\big\}\,.
\end{align*}
Here, the U-statistics are defined by
\begin{align}\label{q6}
&\tilde D_{\nu,1}(\T_{n,2})=\frac{1}{n_1(n_1-1)}\sum_{(j,k)\in I_2(n_1)}\phi_{\nu}( Z_{j}-Z_{k})\,,\notag\\
&\tilde D_{\nu,2}(\T_{n,2})=\frac{1}{n_2(n_2-1)}\sum_{(j,k)\in I_2(n_2)}\phi_{\nu}(Z_{n_1+j}-Z_{n_1+k})\,,\notag\\
&\tilde D_{\nu,3}(\T_{n,2})=\frac{1}{(n_1+n_2)(n_1+n_2-1)}\sum_{(j,k)\in I_2(n_1+n_2)}\phi_{\nu}(Z_j-Z_k)\,.
\end{align}

For the first term $J_{\nu,1}$, by Theorem~2 in \cite{arcones1995}, it follows that for any $q\geq0$, $\nu\in\V$, and $F\in\mathcal P_1$,
\begin{align}\label{q1}
\P_F\big\{|\tilde D_{\nu,1}(\T_{n,2})-\E_F\tilde D_{\nu,1}(\T_{n,2})|\geq q\big\}\leq 4\exp\bigg(-\frac{nq^2}{2m^2\zeta_1+(2^{m+3}m^m+2m^{-1}/3)q}\bigg)\,,
\end{align}
where $\zeta_1=\var_F\big[\E_F\{\phi_{\nu}(X_2-X_1)|X_{1}\}\big]$. In view of the fact that $\|\phi_{\nu}\|_\infty\leq1$ and the law of total variance, it follows that $\zeta_1\leq\var_F\{\phi_{\nu}(X_2-X_1)\}\leq\|\phi_{\nu}\|_\infty^2\leq1$.
Since $|J_{\nu,1}|$ is bounded, so that we may assume further that $q$ in \eqref{q1} is such that $q\leq c$. In view of the fact that $n_1\asymp n_3$, we therefore deduce from \eqref{q1} that, for any $0\leq q\leq c$,
\begin{align}
\P_F\big(|J_{\nu,1}|\geq q/3\big)&=\P_F\bigg\{|\tilde D_{\nu,1}(\T_{n,2})-\E_F\tilde D_{\nu,1}(\T_{n,2})|\geq \frac{n_1n_2q}{3(n_1+n_2)(n_1-1)}\bigg\}\notag\\
&\leq 4\exp(-cnq^2)\,.
\end{align}
where $c>0$ does not depend on $n$, $\nu\in\V$ or $F\in\mathcal P_1$. Applying the same calculations on $J_{\nu,2}$ and $J_{\nu,3}$, we obtain that
\begin{align*}
\P_F(|\ed_\nu^2(\hat{\mathbb P}^X,\hat{\mathbb P}^Y)-\E\ed_\nu^2(\hat{\mathbb P}^X,\hat{\mathbb P}^Y)|>q)&\leq \P_F\big(|J_{\nu,1}|\geq q/3\big)+\P_F\big(|J_{\nu,2}|\geq q/3\big)+\P_F\big(|J_{\nu,3}|\geq q/3\big)\\
&\leq c_1\exp(-c_2nq^2)\,.
\end{align*}
where $c_1,c_2>0$ are absolute constants.
By taking $q=c_2^{-1/2}\sqrt{n^{-1}\log(3c_1\beta^{-1})}$, we obtain that the above quantity is further bounded by $\beta/3$, which implies
\begin{align}\label{p303}
\sup_{\nu\in\V}\sup_{F\in\mathcal P_1}\kappa_{\nu,1-\beta/3}\leq cn^{-1/2}\log(\beta^{-1})\,,
\end{align}
where $c>0$ is an absolute constant.

\textbf{Step 3.} Bounding the bias $|\|\phi_{1}-\phi_2\|_\nu^2-\E_F\{\ed_\nu^2(\hat{\mathbb P}^X,\hat{\mathbb P}^Y)\}|$.

Observing Proposition~\ref{prop:edhat}, direct calculations yield
\begin{align*}
\E_F\{\ed_\nu^2(\hat{\mathbb P}^X,\hat{\mathbb P}^Y)\}
%
&=(1-n_1^{-1})\E_F\phi_{\nu}(X_1-X_2)+(1-n_2^{-1})\E_F\phi_{\nu}(Y_1-Y_2)-2\E_F\phi_{\nu}(X_1-Y_1)+\frac{n_1+n_2}{n_1n_2}\\
&=\|\phi_X-\phi_Y\|_\nu^2-n_1^{-1}\E_F\phi_{\nu}(X_1-X_2)-n_2^{-1}\E_F(Y_1-Y_2)+\frac{n_1+n_2}{n_1n_2}\,.
\end{align*}
The above equation together with the fact that $\|\phi_{\nu}\|_\infty\leq 1$ imply that
\begin{align*}
\big|\E_F\{\ed_\nu^2(\hat{\mathbb P}^X,\hat{\mathbb P}^Y)\}-\|\phi_X-\phi_Y\|_\nu^2\big|\leq n_1^{-1}|\E_F\phi_{\nu}(X_1-X_2)|+n_2^{-1}|\E_F\phi_{\nu}(Y_1-Y_2)|+\frac{n_1+n_2}{n_1n_2}\leq cn_1^{-1}\,,
\end{align*}
where $c>0$ is an absolute constant.
Combining the above equation with \eqref{qbeta}, \eqref{p302} and \eqref{p303} completes the proof.

\subsubsection*{Lower bound of the minimax separation rate}


The proof follows similar ideas as the proof of Theorem~\ref{thm:lower:ind}.
Suppose $\{\eta_\ell\}_{\ell\geq1}$ is an orthogonal basis of $\H$. Let
\begin{align}\label{ccc}
C_0=\sum_{j=1}^\infty\lambda_j\l\eta_j,\cdot\r\eta_j\,,\qquad\text{where }\lambda_j=j^{-2}\,.
\end{align}
Let $\nu_0$ denote the mean-zero Gaussian measure with covariance operator $C$.
For $2\leq k\leq M$, let
\begin{align}\label{cn1k} C_{n_1;k}=C_0+\varrho_n\sqrt{\lambda_1\lambda_k}(\l\eta_1,\cdot\r\eta_k+\l\eta_k,\cdot\r\eta_1)\,.
\end{align}
Let $\zeta_1^*,\{\zeta_{\ell}\}_{\ell\geq1},\{\tilde\zeta_{\ell}\}_{\ell\geq1}$ be independent and identically distributed real-valued standard normal random variables, and let
\begin{align}\label{z}
&Z_{(0)}=\sum_{\ell=1}^\infty\sqrt{\lambda_\ell}\,\zeta_{\ell}\,\eta_\ell\,;\\
&Z_{(n_1;k)}=(\sqrt{\lambda_1}\eta_{1}+\sqrt{\lambda_k}\eta_k)\sqrt{\varrho_{n_1}}\zeta_1^*+(\sqrt{\lambda_1}\eta_{1}\zeta_1+\sqrt{\lambda_k}\eta_k\zeta_k)\sqrt{1-\varrho_{n_1}}+\sum_{\ell\geq2,\ell\neq k}^\infty\sqrt{\lambda_\ell}\,\zeta_{\ell}\,\eta_\ell\,,\notag
\end{align}
where $\varrho_n$ is to be specified later. Let $F_0$ and $F_{n_1;k}$ be the probability distributions of $Z_{(0)}$ and $Z_{(n_1;k)}$, respectively, and let $\phi_{(0)}$ and $\phi_{(n_1;k)}$ denotes their corresponding characteristic functionals.

It follows from Lemma~\ref{lem:klexpan} that $F_0$ is the mean-zero Gaussian distribution with covariance operator $C$ defined in \eqref{ccc}. For $F_{n_1;k}$, direct calculations yield that its characteristic functional is given by
\begin{align*}
\phi_{n_1;k}(w)&=\E\exp\bigg\{(\sqrt{\lambda_1}\l\eta_{1},w\r+\sqrt{\lambda_k}\l\eta_k,w\r)\sqrt{\varrho_{n_1}}\zeta_1^*+\sqrt{\lambda_1(1-\varrho_{n_1})}\l\eta_{1},w\r\zeta_1\\
&\qquad\qquad+\sqrt{\lambda_k(1-\varrho_{n_1})}\l\eta_k,w\r\zeta_k+\sum_{\ell\geq2,\ell\neq k}^\infty\sqrt{\lambda_\ell}\,\zeta_{\ell}\,\l\eta_\ell,w\r\bigg\}\\
&=\exp\bigg\{-\frac{1}{2}\varrho_{n_1}(\sqrt{\lambda_1}\l\eta_{1},w\r+\sqrt{\lambda_k}\l\eta_k,w\r)^2-\frac{1}{2}\lambda_1(1-\varrho_{n_1})\l\eta_{1},w\r^2\\
&\qquad\qquad-\frac{1}{2}\lambda_k(1-\varrho_{n_1})\l\eta_k,w\r^2-\frac{1}{2}\sum_{\ell\geq2,\ell\neq k}^\infty\lambda_\ell\l\eta_\ell,w\r^2\bigg\}\\
&=\exp\bigg(-\frac{1}{2}\sum_{\ell=1}^\infty\lambda_\ell\l\eta_\ell,w\r^2-\varrho_{n_1}\sqrt{\lambda_1\lambda_k}\l\eta_1,w\r\l\eta_k,w\r\bigg)=\exp(-2^{-1}\l C_{n_1;k}w,w\r)\,.
\end{align*}
The above equation implies that $F_{(n_1;k)}$ is the mean-zero Gaussian distribution with covariance operator $C_{n_1;k}$ defined in \eqref{cn1k}. Furthermore, it holds that
\begin{align}\label{phin1k}
\phi_{n_1;k}(w)=\phi_{(0)}(w)\exp(-\varrho_n\sqrt{\lambda_1\lambda_k}\l\eta_1,w\r\l\eta_k,w\r)\,.
\end{align}

Let $F_0$ denote the joint distribution of the sample $\T_{n,2}$ under the null hypothesis, where $Z_{1},Z_2,\ldots,Z_{n_1+n_2}$ are i.i.d.~copies of $Z_{(0)}$ in \eqref{z}, that is,
\begin{align*}
Z_{1},\ldots,Z_{n_1},Z_{n_1+1},\ldots,Z_{n_1+n_2}\iidsim\mathcal N(0,C_0)\,.
\end{align*}
For $2\leq k\leq M$, let $F_{n_1;k}$ denote the joint distribution of the sample $\T_{n,2}$, where $Z_{1},\ldots,Z_{n_1}$ are i.i.d.~copies of $Z_{(1)}$ and $Z_{n_1+1},\ldots,Z_{n_1+n_2}$ are i.i.d.~copies of $Z_{(n_1;k)}$, from two independent samples, that is,
\begin{align*}
Z_{1},\ldots,Z_{n_1}\iidsim\mathcal N(0,C_0),\qquad Z_{n_1+1},\ldots,Z_{n_1+n_2}\iidsim\mathcal N(0,C_{n_1;k})\,.
\end{align*}
Let $F_{n,2}$ denote the uniform probability measure on the set $\{F_{n_1,2},,\ldots,F_{n_1,M}\}$.

Following the proof of Theorem~\ref{thm:lower:ind}, we aim to prove
\begin{align}\label{el2}
\E_{\P_{F_0}}\bigg\{\frac{\d\P_{F_{n,2}}(\T_{n,2})}{\d\P_{F_0}(\T_{n,2})}\bigg\}^2\leq1+4(1-\alpha-\beta)^2\,,
\end{align}
so that a lower bound for minimax risk is given by $\beta_{n_1,\varrho}^*\geq\beta$.

Next, we compute the energy distance $\ed_{\nu_0}^2(Z_{(0)},Z_{(n_1;k)})$ with respect to~the reference probability measure $\nu_0$. It follows from \eqref{phin1k} and calculations similar to the ones to derive \eqref{ink0} and \eqref{p2} in the proof of Theorem~\ref{thm:lower:ind} that
\begin{align*}
\ed_{\nu_0}^2(Z_{(0)},Z_{(n_1;k)})&=\int_\H|\phi_{(0)}(w)-\phi_{(n_1;k)}(w)|^2\,\nu_0(\d w)\\
&=\int_\H\exp(-\l C_0w,w\r)\{\exp(-\varrho_n\sqrt{\lambda_1\lambda_k}\l\eta_1,w\r\l\eta_k,w\r)-1\}^2\,\nu_0(\d w)\\
&\geq(1-e^{-1/2})^2\varrho_n^2\lambda_1\lambda_k\int_\H\exp(-2\l C_0w,w\r)\l\eta_1,w\r^2\l\eta_k,w\r^2\,\nu_0(\d w)\\
&\geq(1-e^{-1/2})^2\varrho_n^2\times\frac{\lambda_1^2\lambda_k^2}{(4\lambda_1^2+1)(4\lambda_k^2+1)}\times\{\E{\phi_{\nu_0}}(2\tilde X)\}^2\geq c\varrho_n^2\,,
\end{align*}
where $\tilde X\sim \N(0,2C_0)$. We then proceed to verify \eqref{el2}. Observe that
$$C_0^{-1/2}C_{n_1;k}C_0^{-1/2}=\varrho_n(\l\cdot,\eta_1\r\eta_k+\l\cdot,\eta_k\r\eta_1)+\sum_{\ell=1}^\infty\l\eta_\ell,\cdot\r\eta_\ell=I+\varrho_n(\l\cdot,\eta_1\r\eta_k+\l\cdot,\eta_k\r\eta_1)\,.$$
In view of Lemma~\ref{lem:rd}, we compute the eigenvalues and eigenvectors of the following $\H$-operator:
\begin{align*}
S_{n_1;k}=I-C_0^{-1/2}C_{n_1;k}C_0^{-1/2}=-\varrho_n(\l\cdot,\eta_1\r\eta_k+\l\cdot,\eta_k\r\eta_1)\,.
\end{align*}
Observe that $S_{n_1;k}$ is a rank-2 operator on $\H$, and
\begin{align*}
&S_{n,k}(\eta_1+\eta_k)=-\varrho_n(\eta_1+\eta_k)\,;\qquad S_{n,k}(\eta_1-\eta_k)=\varrho_n(\eta_1-\eta_k)\,.
\end{align*}
This implies that the nonzero eigenvalues of $S_{n,k}$ are $\alpha_{k,1}=-\varrho_n$ and $\alpha_{k,2}=\varrho_n$ with corresponding eigenvectors $\xi_{k,1}=(\eta_1+\eta_k)/\sqrt{2}$ and $\xi_{k,2}=(\eta_1-\eta_k)/\sqrt{2}$. 
Therefore, applying Lemma~\ref{lem:rd}, we deduce that
\begin{align*}
&\frac{\d\P_{F_{n_1;k}}(z)}{\d\P_{F_0}(z)}=(1-\varrho_n^2)^{-1/2}\exp\bigg[-\frac{1}{2}\Big\{\frac{\alpha_{k,1}}{1-\alpha_{k,1}}\l z,C_{0}^{-1/2}\xi_{k,1}\r^2+\frac{\alpha_{k,2}}{1-\alpha_{k,2}}\l z,C_{0}^{-1/2}\xi_{k,2}\r^2\Big\}\bigg]\\
&=(1-\varrho_n^2)^{-1/2}\exp\bigg[-\frac{1}{4}\Big\{\frac{-\varrho_n}{1+\varrho_n}\l z,C_0^{-1/2}(\eta_1+\eta_k)\r^2+\frac{\varrho_n}{1-\varrho_n}\l z,C_0^{-1/2}(\eta_1-\eta_k)\r^2\Big\}\bigg]\\
&=(1-\varrho_n^2)^{-1/2}\exp\bigg\{\frac{\varrho_n}{4(1+\varrho_n)}(\lambda_1^{-1/2}\l z,\eta_1\r+\lambda_k^{-1/2}\l z,\eta_k\r)^2-\frac{\varrho_n}{4(1-\varrho_n)}(\lambda_1^{-1/2}\l z,\eta_1\r-\lambda_k^{-1/2}\l z,\eta_k\r)^2\bigg\}\,.
\end{align*}
Observe that by construction, we have
\begin{align*}
&\frac{\d\P_{F_{n,2}}(\T_{n,2})}{\d\P_{F_0}(\T_{n,2})}=(M-1)^{-1}\sum_{k=2}^{M}\prod_{j=1}^{n_2}\frac{\d\P_{F_{n_1;k}}(Z_{n_1+j})}{\d\P_{F_0}(Z_{n_1+j})}\\
&=(M-1)^{-1}(1-\varrho_n^2)^{-n_2/2}\sum_{k=2}^{M}\prod_{j=1}^{n_2}\exp\bigg\{\frac{\varrho_n}{4(1+\varrho_n)}(\lambda_1^{-1/2}\l Z_{n_1+j},\eta_1\r+\lambda_k^{-1/2}\l Z_{n_1+j},\eta_k\r)^2\\
&\hspace{7cm}-\frac{\varrho_n}{4(1-\varrho_n)}(\lambda_1^{-1/2}\l Z_{n_1+j},\eta_1\r-\lambda_k^{-1/2}\l Z_{n_1+j},\eta_k\r)^2\bigg\} \,.
\end{align*}
Therefore, we obtain from the above equation that
\begin{align*}
\E_{\P_{F_0}}\bigg\{\frac{\d\P_{F_{n,2}}(\T_{n,2})}{\d\P_{F_0}(\T_{n,2})}\bigg\}^2=J_{n,1}+J_{n,2}\,,
\end{align*}
where
\begin{align*}
&J_{n,1}=(M-1)^{-2}(1-\varrho_n^2)^{-n_2}\\
&\sum_{2\leq k_1\neq k_2\leq M}\bigg(\int_{\H}\exp\bigg[\frac{\varrho_n}{4(1+\varrho_n)}\Big\{(\lambda_1^{-1/2}\l z,\eta_1\r+\lambda_{k_1}^{-1/2}\l z,\eta_{k_1}\r)^2+(\lambda_1^{-1/2}\l z,\eta_1\r+\lambda_{k_2}^{-1/2}\l z,\eta_{k_2}\r)^2\Big\}\\
&\qquad~~~-\frac{\varrho_n}{4(1-\varrho_n)}\Big\{(\lambda_1^{-1/2}\l z,\eta_1\r-\lambda_{k_1}^{-1/2}\l z,\eta_{k_1}\r)^2+(\lambda_1^{-1/2}\l z,\eta_1\r-\lambda_{k_2}^{-1/2}\l z,\eta_{k_2}\r)^2\Big\}\bigg]\nu_0(\d z)\bigg)^{n_2}\,;\\
&J_{n,2}=(M-1)^{-2}(1-\varrho_n^2)^{-n_2}\sum_{k=1}^{M}\bigg[\int_{\H}\exp\bigg\{\frac{\varrho_n}{2(1+\varrho_n)}(\lambda_1^{-1/2}\l z,\eta_1\r+\lambda_k^{-1/2}\l z,\eta_k\r)^2\\
&\hspace{7cm}-\frac{\varrho_n}{2(1-\varrho_n)}(\lambda_1^{-1/2}\l z,\eta_1\r-\lambda_k^{-1/2}\l z,\eta_k\r)^2\bigg\}\nu_0(\d z)\bigg]^{n_2}\,.
\end{align*}

For the first term $J_{n,1}$, by writing $z_\ell=\l z,\eta_\ell\r$ for $\ell\geq1$,  following similar calculations that yield \eqref{jn1} in the proof of Theorem~\ref{thm:lower:ind}, we obtain
\begin{align*}
&J_{n,1}=(M-1)^{-2}(1-\varrho_n^2)^{-n_2}\\
&\sum_{2\leq k_1\neq k_2\leq M}\bigg(\int_{\mathbb R^3}\exp\bigg[\frac{\varrho_n}{4(1+\varrho_n)}\big\{(\lambda_1^{-1/2}z_1+\lambda_{k_1}^{-1/2}z_{k_1})^2+(\lambda_1^{-1/2}z_1+\lambda_{k_2}^{-1/2}z_{k_2})^2\big\}\\
&-\frac{\varrho_n}{4(1-\varrho_n)}\big\{(\lambda_1^{-1/2}z_1-\lambda_{k_1}^{-1/2}z_{k_1})^2+(\lambda_1^{-1/2}z_1-\lambda_{k_2}^{-1/2}z_{k_2})^2\big\}\bigg]\,\N_{\lambda_1}(\d z_1)\,\N_{\lambda_{k_1}}(\d z_{k_1})\,\N_{\lambda_{k_2}}(\d z_{k_2})\bigg)^{n_2}\\
&=(M-1)^{-2}(1-\varrho_n^2)^{-n_2}\sum_{2\leq k_1\neq k_2\leq M}\bigg\{\int_{\mathbb R^3}\exp\big([s_1,s_{k_1},s_{k_2}]\Sigma_2[s_1,s_{k_1},s_{k_2}]\trans\big) \N_{1}(\d s_1)\N_{1}(\d s_{k_1})\N_{1}(\d s_{k_2})\bigg\}^{n_2}\,,
\end{align*}
where
\begin{align*}
\Sigma_2=\left[\begin{matrix}
\frac{-\varrho_n^2}{1-\varrho_n^2} & \frac{\varrho_n}{2(1-\varrho_n^2)} &\frac{\varrho_n}{2(1-\varrho_n^2)} \\
\frac{\varrho_n}{2(1-\varrho_n^2)}&\frac{-\varrho_n^2}{2(1-\varrho_n^2)} & 0\\
\frac{\varrho_n}{2(1-\varrho_n^2)} & 0& \frac{-\varrho_n^2}{2(1-\varrho_n^2)} 
\end{matrix}\right]\,.
\end{align*}
Note that $\det(I_3-2\Sigma_2)=(1-\varrho_n^2)^{-2}$. Applying equation \eqref{formula}, we deduce that for $n$ large enough,
\begin{align}\label{j1}
J_{n,1}&=(M-2)(M-1)^{-1}(1-\varrho_n^2)^{-n_2}\det(I_3-2\Sigma_2)^{-n_2/2}\notag\\
&=(M-2)(M-1)^{-1}(1-\varrho_n^2)^{-n_2/2}\leq\exp(n_2\varrho_n^2/2)\{1+o(1)\}\,.
\end{align}
For the second term $J_{n,2}$, it follows from similar calculations that
\begin{align*}
J_{n,2}&=(M-1)^{-2}(1-\varrho_n^2)^{-n_2}\sum_{k=2}^{M}\bigg[\int_{\mathbb R^2}\exp\bigg\{\frac{\varrho_n}{2(1+\varrho_n)}(\lambda_1^{-1/2}z_1+\lambda_{k}^{-1/2}z_{k})^2\\
&\hspace{5.5cm}-\frac{\varrho_n}{2(1-\varrho_n)}(\lambda_1^{-1/2}z_1-\lambda_{k}^{-1/2}z_{k})^2\bigg\}\,\N_{ \lambda_1}(\d z_1)\,\N_{ \lambda_k}(\d z_k)\bigg]^{n_2}\\
&=(M-1)^{-1}(1-\varrho_n^2)^{-n_2}\bigg[\int_{\mathbb R^2}\exp\bigg\{\frac{\varrho_n(s+t)^2}{2(1+\varrho_n)}-\frac{\varrho_n(s-t)^2}{2(1-\varrho_n)}\bigg\}\,\N_{1}(\d s)\,\N_{1}(\d t)\bigg]^{n_2}\\
&=(M-1)^{-1}(1-\varrho_n^2)^{-n_2}\bigg\{\int_{\mathbb R^2}\exp([s,t]\Sigma_1[s,t]\trans)\,\N_{1}(\d s)\,\N_{1}(\d t)\bigg\}^{n_2}\,,
\end{align*}
where $\Sigma_1$ is defined in \eqref{sigma}.
Observing that $\det(I_2-2\Sigma_1)=1$. We therefore obtain from the above equation and \eqref{formula} that, for $n$ large enough,
\begin{align*}
J_{n,2}&=(M-1)^{-1}(1-\varrho_n^2)^{-n_2}\det(I_2-2\Sigma_1)^{-n/2}=(M-1)^{-1}(1-\varrho_n^2)^{-n_2}\notag\\
&=(M-1)^{-1}\exp(n_2\varrho_n^2)\{1+o(1)\}\,.
\end{align*}
Therefore, combining the above equation with \eqref{j1} yields that, for $n$ large enough,
\begin{align*}
\E_{\P_{F_0}}\bigg\{\frac{\d\P_{F_n}(\T_{n,2})}{\d\P_{F_0}(\T_{n,2})}\bigg\}^2&\leq\exp(n_2\varrho_n^2)+(M-1)^{-1}\exp(2n_2\varrho_n^2)\,.
\end{align*}
Therefore, taking $M=\lceil 4^{-1}(1-\alpha-\beta)^{-2}\rceil+2>4^{-1}(1-\alpha-\beta)^{-2}$ and $\varrho_n=(cn)^{-1/2}$ suffices \eqref{el2} and completes the proof, where 
\begin{align*}
c=\log\big\{2^{-1}\sqrt{M^2+4M+4(1-\alpha-\beta)^2M}-M/2\big\}>0\,.
\end{align*}

\subsection{Proof of Proposition~\ref{prop:me:two}}\label{app:prop:me:two}

We obtain from the proof of Theorem~\ref{thm:upper:ind} that, for some absolute constant $c>0$, let $\rho_0=c(n_1+n_2)^{-1/2}\big\{\log(\alpha^{-1})+\log(\beta^{-1})\big\}$. A sufficient condition for
\begin{align*}
\sup_{\nu\in\V_G(c)}\sup_{F\in\mathcal F_\nu(\rho_0)}\P_F\{\tilde\psi_{\nu,n,2}(\alpha)=0\}\leq\beta
\end{align*}
is that $\int_\H|\phi_{\tilde X}(w)-\phi_{\tilde Y}(w)|^2\nu(\d w)\geq\rho_0$.
It follows from the Cauchy-Schwarz inequality that
\begin{align}\label{rpre}
\int_\H|\phi_{\tilde X}(w)-\phi_{\tilde Y}(w)|^2\nu(\d w)&=\int_\H|\phi_X(w)-\phi_Y(w)|^2|\phi_U(u)|^2\nu(\d w)\notag\\
&\geq\bigg[\int_\H|\phi_X(w)-\phi_Y(w)|\,\nu(\d w)\bigg]^2\bigg[\int_\H|\phi_U(w)|^{-2}\,\nu(\d w)\bigg]^{-1}\,.
\end{align}

Next, we compute $\int_\H|\phi_U(w)|^{-2}\nu(\d w)$. By Proposition~1.13 in \cite{daprato2006}, we have that, for the mean-zero Gaussian measure $\nu$,
\begin{align}\label{R}
R(a):=\int_\H\exp(a\|w\|^2/2)\,\nu(\d w)=\bigg\{\prod_{k=1}^\infty(1-a\lambda_k)\bigg\}^{-1/2}\,,\qquad\quad a\in\mathbb R\,,
\end{align}
where $\{\lambda_k\}_{k=1}^\infty$ are the eigenvalues of the covariance operator of $\nu$.
Observing the fact that, for $w\in\H$, $(1+\|w\|)^{s}\leq c_U(1+\|w\|^s)$, in view of \eqref{Fu}, it holds that, for the real-valued function $R$ in \eqref{R},
\begin{align*}
\int_\H|\phi_U(w)|^{-2}\nu(\d w)\leq c_s \int_\H(1+\|w\|^{2s})\nu(\d w)\leq c\int_\H\|w\|^{2s}\nu(\d w)=2^s\frac{d^s}{d a^s}R(0)\,,
\end{align*}
where, for $s\in\mathbb R$, $\frac{d^s}{d a^s}$ denotes the $s$-order Riemann-Liouville fractional derivative. Observe the fact that for $t\in\mathbb R$,
\begin{align*}
\frac{d^s}{d a^s}(a^t)=\frac{t!}{(t-a)!}a^{t-s}\,.
\end{align*}
Then, it follows from direct calculations that
\begin{align*}
\int_\H|\phi_U(w)|^{-2}\nu(\d w)\leq c_s\{{\rm tr}(C_\nu)\}^{s}\,.
\end{align*}
where $c_s>0$ only depends on $s$. The proof is therefore complete in view of \eqref{rpre}.

\subsection{Proof of Proposition~\ref{prop:agg}}\label{app:prop:agg}

For (i), under the null hypothesis of independence, $\{\hat\dcov_\nu^2(X,Y,\tilde\Pi_b)\}_{1\leq b\leq B_2}$ are identically distributed as $\hat\dcov_\nu^2(X,Y)$, so that for any $F\in\mathcal P_0$, it holds that $\P_F\{\hat\psi_n(\alpha)\}\leq\alpha$.

For (ii), let $T_b(u)=\one\big[\max_{\nu\in\V_0}\{\hat\dcov_\nu^2(X,Y,\tilde\Pi_b)-\hat q_{\nu}(1-u)\}>0\big]$,
so that $\hat\psi_n(\alpha)=B_2^{-1}\sum_{b=1}^{B_2}T_b(\hat u_\alpha)$ in view of \eqref{hatpsi:ind}. 
Define the event
\begin{align*}
\mathcal E=\bigg\{\frac{1}{B_2}\sum_{b=1}^{B_2}T_b(u)-\E\{T_b(u)\}\leq\sqrt{(2B_2)^{-1}\log(2/\beta)}\bigg\}\,.
\end{align*}
Let $\boldsymbol{\Pi}_{B_1}=\{\Pi_1,\ldots,\Pi_{B_1}\}$. Conditional on $\T_n$ and $\boldsymbol{\Pi}_{B_1}$, the random variables $T_1(u),\ldots,T_{B_2}(u)$ are independent. By the conditional Hoeffding's inequality on $\T_n$ and $\boldsymbol{\Pi}_{B_1}$, we obtain that, for any $\boldsymbol{\Pi}_{B_1}$ and $\T_n$,
\begin{align*}
\P\bigg[\Big|\frac{1}{B_2}\sum_{b=1}^{B_2}T_b(u)-\E\{T_1(u)|\boldsymbol{\Pi}_{B_1},\T_n\}\Big|\geq \sqrt{(2B_2)^{-1}\log(2/\beta)} \,\big|\,\boldsymbol{\Pi}_{B_1},\T_n\bigg]\leq \beta\,,
\end{align*}
for any $\boldsymbol{\Pi}_{B_1}$ and $\T_n$, which implies that $\P(\mathcal E)\leq 1-\beta/2$.  For $F\in\mathcal P_0$, it holds that $\hat\dcov_\nu^2(X,Y,\tilde\Pi_b)$ is identically distributed as $\hat\dcov^2_\nu(X,Y)$.
By taking $u=\alpha/(2|\V_0|)$ and $B_2\geq2\alpha^{-2}\log(2/\beta)$, so that $\sqrt{(2B_2)^{-1}\log(2/\beta)}\leq \alpha/2$, we have, on the event $\mathcal E$, it holds that
\begin{align*}
&\frac{1}{B_2}\sum_{b=1}^{B_2}T_b(u)
%
\leq\P_F\Big[\max_{\nu\in\V_0}\{\hat\dcov_\nu^2(X,Y,\tilde\Pi_b)-\hat q_{\nu}(1-u)\}>0\Big]+\sqrt{(2B_2)^{-1}\log(2/\beta)}\\
&\leq\sum_{\nu\in\V_0} \P_F\big\{\hat\dcov^2_\nu(X,Y)>\hat q_{\nu}(1-u)\big\}+\sqrt{(2B_2)^{-1}\log(2/\beta)}\leq |\V_0|u+\sqrt{(2B_2)^{-1}\log(2/\beta)}\leq\alpha\,.
\end{align*}
In view of the definition of $\hat\psi_n(\alpha)$ in \eqref{hatpsi:ind}, we have $\hat u_\alpha\geq\alpha/(2|\V_0|)$.
We have, for any $F\in\mathcal P_1$, the type-II error is such that
\begin{align*}
&\P_{F}\{\hat\psi_n(\alpha)=0\}=\P_{F}\{\hat\psi_n(\alpha)=0;\mathcal E\}+\P_{F}\{\hat\psi_n(\alpha)=0;\mathcal E^c\}\leq\P_{F}\{\hat\psi_n(\alpha)=0;\mathcal E\}+\P_{F}(\mathcal E^c)\\
&\leq\P_{F}\{\hat\psi_n(\alpha)=0,\mathcal E\}+\beta/2=\P_{F}\Big[\max_{\nu\in\V_0}\{\hat\dcov_\nu^2(X,Y,\tilde\Pi_b)-\hat q_{\nu}(1-\hat u_\alpha)\}\leq0;\mathcal E\Big]+\beta/2\\
&\leq \min_{\nu\in\V_0}\P_{F}\{\hat\dcov_\nu^2(X,Y,\tilde\Pi_b)\leq\hat q_{\nu}(1-\hat u_\alpha);\mathcal E\}+\beta/2\leq \min_{\nu\in\V_0}\P_{F}\{\hat\dcov_\nu^2(X,Y,\tilde\Pi_b)\leq\hat q_{\nu}(1-\alpha/(2|\V_0|))\}+\beta/2\,.
\end{align*}
In view of the definition of $\hat\psi_{\nu,n}$ in \eqref{hatpsi}, we obtain from the above equation that $$\P_{F}\{\hat\psi_n(\alpha)=0\}\leq \beta/2+\P_F\{\hat\psi_{\nu,n}(\alpha/(2|\V_0|))\}\,,$$ so that the result follows in view of Theorem~\ref{thm:upper:ind} and the definition of uniform separation rate.

\section{Auxiliary theoretical details}\label{app:aux}

\subsection{Proof of Lemma~\ref{lem:a.5}}\label{proof:lem:a.5}

It suffices to show that
\begin{align*}
\hat\dcov _\nu^2(X,Y)&=\frac{1}{n^2}\sum_{j,k=1}^n\phi_{\nu}(X_j-X_k)\phi_{\nu}(Y_j-Y_k)+\frac{1}{n^4}\sum_{j,k,\ell,r=1}^n\phi_{\nu}(X_j-X_k)\phi_{\nu}(Y_\ell-Y_r)\\
&\quad-\frac{2}{n^3}\sum_{j,k,\ell=1}^n\phi_{\nu}(X_j-X_k)\phi_{\nu}(Y_j-Y_\ell)\,.
\end{align*}
Observe that $\hat\dcov_\nu^2(X,Y)=I_1+I_2-I_3-I_4$, where
\begin{align*}
I_1&=\int_{\H^2}\hat\phi_{X,Y}(w_1,w_2)\hat\phi_{X,Y}(-w_1,-w_2)\,\nu(\d w_1)\,\nu(\d w_2)\,,\\
I_2&=\int_{\H^2}\hat\phi_{X}(w_1)\hat\phi_{X}(-w_1)\hat\phi_{Y}(w_2)\hat\phi_{Y}(-w_2)\,\nu(\d w_1)\,\nu(\d w_2)\,,\\
I_3&=\int_{\H^2}\hat\phi_{X,Y}(w_1,w_2){\hat\phi_{X}(-w_1)\hat\phi_{Y}(-w_2)}\,\nu(\d w_1)\,\nu(\d w_2)\,,\\
I_4&=\int_{\H^2}\hat\phi_{X,Y}(-w_1,-w_2)\hat\phi_{X}(w_1)\hat\phi_{Y}(w_2)\,\nu(\d w_1)\,\nu(\d w_2)\,.
\end{align*}
For the first term $I_1$, we have
\begin{align*}
&I_1=\int_{\H^2}\bigg\{\frac{1}{n}\sum_{j=1}^n\exp(\i\l X_j,w_1\r+\i\l Y_j,w_2\r)\bigg\}\bigg\{\frac{1}{n}\sum_{k=1}^n\exp(-\i\l X_k,w_1\r-\i\l Y_k,w_2\r)\bigg\}\,\nu(\d w_1)\,\nu(\d w_2)\notag\\
&=\frac{1}{n^2}\sum_{j,k=1}^n\int_{\H^2}\exp(\i\l X_j-X_k,w_1\r+\i\l Y_j-Y_k,w_2\r)\,\nu(\d w_1)\,\nu(\d w_2)=\frac{1}{n^2}\sum_{j,k=1}^n\phi_{\nu}(X_j-X_k)\phi_{\nu}(Y_j-Y_k)\,.
\end{align*}
For the second term $I_2$, we have
{\small
\begin{align*}
I_2&=\frac{1}{n^4}\int_{\H^2}\bigg\{\sum_{j=1}^n\exp(\i\l X_j,w_1\r)\bigg\}
\bigg\{\sum_{k=1}^n\exp(-\i\l X_k,w_1\r)\bigg\}
\bigg\{\sum_{\ell=1}^n\exp(\i\l Y_j,w_1\r)\bigg\}
\bigg\{\sum_{r=1}^n\exp(-\i\l Y_r,w_1\r)\bigg\}\nu(\d w_1)\nu(\d w_2)\\
&=\frac{1}{n^4}\sum_{j,k,\ell,r=1}^n\int_{\H^2}\exp(\i\l w_1,X_j-X_k\r)\times\exp(\i\l w_2,Y_\ell-Y_r\r)\,\nu(\d w_1)\,\nu(\d w_2)=\frac{1}{n^4}\sum_{j,k,\ell,r=1}^n\phi_{\nu}(X_j-X_k)\phi_{\nu}(Y_\ell-Y_r)\,.
\end{align*}
}
For $I_3$ and $I_4$, under Assumption~\ref{a:sym}, we have
{\small
\begin{align*}
&I_3=I_4=\int_{\H^2}\bigg\{\frac{1}{n}\sum_{j=1}^n\exp(\i\l X_j,w_1\r+\i\l Y_j,w_2\r)\bigg\}\bigg\{\frac{1}{n}\sum_{k=1}^n\exp(-\i\l X_k,w_1\r)\bigg\}\bigg\{\frac{1}{n}\sum_{\ell=1}^n\exp(-\i\l Y_\ell,w_1\r)\bigg\}\,\nu(\d w_1)\,\nu(\d w_2)\\
&=\frac{1}{n^3}\sum_{j,k,\ell=1}^n\int_{\H^2}\exp(\i\l X_j-X_k,w_1\r+\i\l Y_j-Y_\ell,w_2\r)\,\nu(\d w_1)\,\nu(\d w_2)=\frac{1}{n^3}\sum_{j,k,\ell=1}^n\phi_{\nu}(X_j-X_k)\phi_{\nu}(Y_j-Y_\ell)\,.
\end{align*}
}

\subsection{Proof of Lemma~\ref{lem:dcovv}}\label{app:lem:dcovv}


In view of the definition of $h$ in \eqref{h}, let
\begin{align*}
&\tilde h\{(x_1,y_1),(x_2,y_2),(x_3,y_3),(x_4,y_4)\}\\
&=\phi_{\nu}(x_{1}-x_{2})\phi_{\nu}(y_{3}-y_{4})+\phi_{\nu}(x_{1}-x_{2})\phi_{\nu}(y_{1}-y_{2})-2\phi_{\nu}(x_{1}-x_{2})\phi_{\nu}(y_{1}-y_{3})\,,
\end{align*}
so that
\begin{align*}
h\{(x_1,y_1),(x_2,y_2),(x_3,y_3),(x_4,y_4)\}=\frac{1}{4!}\sum_{(k_1,k_2,k_3,k_4)\in I_4(4)}\tilde h\{(x_{k_1},y_{k_1}),(x_{k_2},y_{k_2}),(x_{k_3},y_{k_3}),(x_{k_4},y_{k_4})\}\,.
\end{align*}
For $h_1$, 
since $h$ in \eqref{h} and $d_\nu$ in \eqref{d} are symmetric and $(X^{(1)},Y^{(1)})$, $(X^{(2)},Y^{(2)})$, $(X^{(3)},Y^{(3)})$, $(X^{(4)},Y^{(4)})$ are i.i.d., Under the condition that $X\indep Y$,
\begin{align*}
&h_1(x,y)=\frac{1}{4}\sum_{k=1}^4\E\Big[\tilde h\{(X^{(1)},Y^{(1)}),(X^{(2)},Y^{(2)}),(X^{(3)},Y^{(3)}),(X^{(4)},Y^{(4)})\,|\,(X^{(k)},Y^{(k)})=(x,y)\}\Big]\\
&=\E\big\{\phi_{\nu}(x-X^{(2)})\phi_{\nu}(Y^{(3)}-Y^{(4)})+\phi_{\nu}(x-X^{(2)})\phi_{\nu}(y-Y^{(2)})-2\phi_{\nu}(x-X^{(2)})\phi_{\nu}(y-Y^{(3)})\big\}\\
& +\E\big\{\phi_{\nu}(X^{(1)}-x)\phi_{\nu}(Y^{(3)}-Y^{(4)})+\phi_{\nu}(X^{(1)}-x)\phi_{\nu}(Y^{(1)}-y)-2\phi_{\nu}(X^{(1)}-x)\phi_{\nu}(Y^{(1)}-Y^{(3)})\big\}\\
& +\E\big\{\phi_{\nu}(X^{(1)}-X^{(2)})\phi_{\nu}(y-Y^{(4)})+\phi_{\nu}(X^{(1)}-X^{(2)})\phi_{\nu}(Y^{(1)}-Y^{(2)})-2\phi_{\nu}(X^{(1)}-X^{(2)})\phi_{\nu}(Y^{(1)}-y)\big\}\\
& +\E\big\{\phi_{\nu}(X^{(1)}-X^{(2)})\phi_{\nu}(Y^{(3)}-y)+\phi_{\nu}(X^{(1)}-X^{(2)})\phi_{\nu}(Y^{(1)}-Y^{(2)})-2\phi_{\nu}(X^{(1)}-X^{(2)})\phi_{\nu}(Y^{(1)}-Y^{(3)})\big\}\\
&=0\,.
\end{align*}
Hence, we have obtained $\var\{h_1(X,Y)\}=0$.
Moreover, in view of the definition of $h_2$ in \eqref{h2}, we find
\begin{align}\label{t}
&h_2\{(x_1,y_1),(x_2,y_2)\}=\E\big[h\big\{(x_1,y_1),(x_2,y_2),(X^{(3)},Y^{(3)}),(X^{(4)},Y^{(4)})\big\}\big]\notag\\
&=\frac{1}{12}\sum_{1\leq k_1<k_2\leq4}\E\Big[\tilde h\{(X^{(1)},Y^{(1)}),\ldots,(X^{(4)},Y^{(4)})\,|\,(X^{(k_1)},Y^{(k_1)})=(x_1,y_1),(X^{(k_2)},Y^{(k_2)})=(x_2,y_2)\}\Big]\notag\\
&+\frac{1}{12}\sum_{1\leq k_1<k_2\leq4}\E\Big[\tilde h\{(X^{(1)},Y^{(1)}),\ldots,(X^{(4)},Y^{(4)})\,|\,(X^{(k_1)},Y^{(k_1)})=(x_2,y_2),(X^{(k_2)},Y^{(k_2)})=(x_1,y_1)\}\Big]\,.
\end{align}
Note that direct calculations yields
\begin{align*}
&\sum_{1\leq k_1<k_2\leq4}\E\Big[\tilde h\{(X^{(1)},Y^{(1)}),\ldots,(X^{(4)},Y^{(4)})\,|\,(X^{(k_1)},Y^{(k_1)})=(x_1,y_1),(X^{(k_2)},Y^{(k_2)})=(x_2,y_2)\}\Big]\\
%
&=\frac{1}{12}\big\{\phi_{\nu}(x_1-x_2)\E \phi_{\nu}(Y^{(3)}-Y^{(4)})+\phi_{\nu}(x_1-x_2)\phi_{\nu}(y_1-y_2)-2\phi_{\nu}(x_1-x_2)\E \phi_{\nu}(y_1-Y^{(3)})\\
& +2\E \phi_{\nu}(x_1-X^{(2)})\E \phi_{\nu}(y_2-Y^{(4)})-2\E \phi_{\nu}(x_1-X^{(2)})\phi_{\nu}(y_1-y_2)+2\E \phi_{\nu}(x_1-X^{(2)})\E \phi_{\nu}(y_1-Y^{(2)})\\
& -2\E  \phi_{\nu}(X^{(1)}-x_1)\E \phi_{\nu}(Y^{(1)}-Y^{(3)})+\E \phi_{\nu}(X^{(1)}-X^{(2)})\phi_{\nu}(y_1-y_2)+\E \phi_{\nu}(X^{(1)}-X^{(2)})\E \phi_{\nu}(Y^{(1)}-Y^{(2)})\\
& -2\E \phi_{\nu}(X^{(1)}-X^{(2)})\E \phi_{\nu}(Y^{(1)}-y_1)\big\}\,.
\end{align*}
Combining the above equation and \eqref{t}, direct calculations yield
\begin{align*}
h_2\{(x_1,y_1),(x_2,y_2)\}=6^{-1}g_X(x_1,x_2)\,g_Y(y_1,y_2)\,.
\end{align*}
Note that $\E g_X(X_1,X_2) =\E g_Y(Y_1,Y_2) =0$, so that since $X\indep Y$, we obtain
\begin{align*}
\E h_2\{(x_1,y_1),(x_2,y_2)\}=\E g_X(X_1,X_2) \E g_Y(Y_1,Y_2) =0\,.
\end{align*}
Therefore, in view of (ii) of Theorem~\ref{thm:dvar},
\begin{align*}
&\var\big[ h_2\{(X^{(1)},Y^{(1)}),(X^{(2)},Y^{(2)})\}\big]=\E\big[ h_2\{(X^{(1)},Y^{(1)}),(X^{(2)},Y^{(2)})\}\big]^2\\
&=36^{-1}\E\{g_X(X_1,X_2)\}^2\E\{g_Y(Y_1,Y_2)\}^2=36^{-1}\dvar_\nu^2(X)\,\dvar_\nu^2(Y)\,.
\end{align*}
Furthermore, in view of (iii) of Theorem~\ref{thm:dvar}, $\dvar_\nu^2(X)\dvar_\nu^2(Y)\neq0$ since $X$ and $Y$ cannot be degenerate since $X\indep Y$, which proves the claim.

\baselineskip=11pt
\vspace{1cm}

\noindent

\begingroup
\renewcommand{\section}[2]{\subsection#1{#2}}
{\centering

}
\endgroup

\end{document}